\documentclass[a4paper,10pt,notitlepage,reqno]{article}
\usepackage[utf8]{inputenc}
\usepackage[english]{babel}
\usepackage{amssymb,amsthm,bm} 
\usepackage{mathtools}
\usepackage{mathrsfs}
\usepackage{enumitem}
\usepackage{authblk}
\usepackage{cite}
\usepackage{xcolor}
\usepackage{tikz}
\usepackage{stackrel}
\usepackage{geometry}
\usepackage{hyperref}
\theoremstyle{plain} 
\newtheorem{theorem}{Theorem}[section]
\newtheorem*{theorem*}{Theorem}
\newtheorem{lemma}{Lemma}[section]
\newtheorem{proposition}{Proposition}[section]
\newtheorem{corollary}[theorem]{Corollary}
\newtheorem*{corollary*}{Corollary}

\theoremstyle{definition}

\newtheorem*{definition}{Definition}

\theoremstyle{remark}
\newtheorem{remark}{Remark}[section]

\numberwithin{equation}{section}

\usepackage{color}
\usepackage[normalem]{ulem}
\definecolor{DarkGreen}{rgb}{0,0.5,0.1} 

\newcommand\soutD{\bgroup\markoverwith
{\textcolor{DarkGreen}{\rule[.5ex]{2pt}{1pt}}}\ULon}

\definecolor{lorenzo}{rgb}{0.2, 0.65, 0.75}

\newcommand{\Hm}[1]{\leavevmode{\marginpar{\tiny%
$\hbox to 0mm{\hspace*{-1.5mm}$\leftarrow$\hss}%
\vcenter{\vrule depth 0.1mm height 0.1mm width \the\marginparwidth}%
\hbox to
0mm{\hss$\rightarrow$\hspace*{-1.5mm}}$\\\relax\raggedright #1}}}

\definecolor{Darkgblue}{rgb}{0.3,0.3,0.5}

\newcommand{\R}{\mathbb{R}}
\newcommand{\C}{\mathbb{C}}
\DeclareMathOperator{\Div}{div}

\DeclareMathOperator\arctanh{arctanh}

\DeclareMathOperator{\realpart}{Re}
\renewcommand{\Re}{\realpart}

\DeclareMathOperator{\imaginarypart}{Im}
\renewcommand{\Im}{\imaginarypart}

\title{A unified approach to Hardy-type inequalities with Bessel pairs}

\author[1]{Lucrezia Cossetti}
\affil[1]{{Ikerbasque \& Universidad del Pa\'is Vasco/Euskal Herriko Unibertsitatea, UPV/EHU, \newline Aptdo. 644, 48080, Bilbao, Spain; lucrezia.cossetti@ehu.eus}}

\author[2]{Lorenzo D'Arca}
\affil[2]{{Department of Mathematics “Guido Castelnuovo”, Sapienza University of Rome,\newline Piazzale Aldo Moro 5, Roma
00185, Italy; lorenzo.darca@uniroma1.it}}

\begin{document}
\maketitle

\begin{abstract}
	\noindent	
In this paper, we provide suitable characterisations of pairs of weights $(V,W),$ known as Bessel pairs, that ensure the validity of weighted Hardy-type inequalities. The abstract approach adopted here makes it possible to establish such inequalities also going beyond the classical Euclidean setting and also within a more general $L^p$ framework. As a byproduct of our method, we obtain explicit expressions for the maximizing functions and, in certain specific situations, we show that the associated constants are sharp. We emphasise that our approach unifies, generalises and improves several existing results in the literature.
\end{abstract}

\tableofcontents

\section{Introduction}
In the pioneering work~\cite{GM2011}, Ghoussoub and Moradifam introduced the notion of Bessel pairs in relation to the validity of suitable weighted Hardy-type inequalities, more precisely they established necessary and sufficient conditions on a pair of positive, radial functions $V$ and $W$ on a ball $B_R$ with radius $R$ in $\R^N,$ under which the following Hardy inequality holds 
\begin{equation}\label{eq:weighted-Hardy-L2}
	\int_{B_R} V(x) |\nabla u|^2\, dx \geq \int_{B_R} W(x) |u|^2\, dx,
	\quad \forall\,  u \in C^\infty (B_R).
\end{equation}
More specifically, they showed that~\eqref{eq:weighted-Hardy-L2} holds true if and only if the ordinary differential equation
\begin{equation}\label{eq:ODE}
	y''(r) + \Big(\frac{N-1}{r} + \frac{V_r(r)}{V(r)} \Big)y'(r) + \frac{W(r)}{V(r)} y(r)=0
\end{equation}
admits a positive solution on the interval $(0,R].$ 
This characterisation makes a very useful connection between the validity of weighted Hardy-type inequalities in the spirit of~\eqref{eq:weighted-Hardy-L2} and the oscillatory behaviour of certain ordinary differential equations like~\eqref{eq:ODE}, which allowed the authors to improve, extend and unify many available results on Hardy-type inequalities (refer to~\cite{GM2011} and references therein). 
This connection was further developed by the same authors in the monograph~\cite{GM2013-book}, where they addressed, among other issues, the possibility to establish an analogous correlation for Hardy-type inequalities for more general uniformly elliptic operators. 
More precisely, they posed the problem of developing suitable characterisations for a pair of functions $(V,W)$ (radial or not) in order to guarantee the validity of the Hardy-type inequality 
\begin{equation}\label{eq:general}
	\int_\Omega V(x) |\nabla u|_{A}^2\, dx \geq \int_\Omega W(x) |u|^2\, dx, 
\end{equation} 
where $|\xi|_{A}:=\langle A(x) \xi, \xi \rangle,$ for $\xi\in \R^N$ and where $A(x)$ is a given symmetric, uniformly positive definite $N\times N$ matrix. In~\cite[Ch.9]{GM2013-book}, inequality~\eqref{eq:general} was proved for the pair $(V,W)=(1, \frac{1}{4} \frac{|\nabla E|_{\!A}^2}{E^2} ),$ where $E$ is a positive solution to $-\Div(A \nabla E)\, dx= \mu$ on $\Omega,$ where $\mu$ is any non-negative non-zero finite measure on $\Omega.$ An $L^p$- version of~\eqref{eq:general}, with $p\neq 2,$ was provided in the same work.

In the recent paper~\cite{RS2025}, the authors generalised the results obtained in~\cite{GM2013-book}. More precisely, they proved the inequality~\eqref{eq:general} for a more general class of positive, \emph{radial} Bessel pairs $(V,W).$ Some examples in a non-Euclidean setting are provided too (refer to~\cite[Thm.2.2]{RS2025} for more details). 

Further results in a non-Euclidean context were obtained in~\cite{GJR24}. Here, the notion of Bessel pairs was developed to derive Hardy and Rellich identities and inequalities for Baouendi--Grushin operators. For this type of operators refined weighted $L^p$-Hardy and Caffarelli--Kohn--Nirenberg inequalities have also been established in~\cite{YZ25, SYZ25, AYZ25}. Improved Hardy-type inequalities in the Heisenberg group in the presence of magnetic fields have been obtained in~\cite{CFKP23}.

\medskip
The aim of our paper is to generalise and unify the aforementioned results by establishing this type of inequalities within a broad abstract setting, from which the particular cases follow as (almost) direct consequences. 

Our approach makes full use of the classical technique of \emph{developing a square} and integrating by parts, a method traditionally employed, among other proofs, to establish the classical Hardy inequality in $L^2.$ 

More precisely, we will make essential use of identity~\eqref{eq: identita_fondamentale} below, introduced for real-valued functions in~\cite[Prop. 2.1]{DA} and here extended to complex-valued functions, and of its vector-valued counterpart~\eqref{eq: identita_fondamentale_campivettoriali}. 

Once identities~\eqref{eq: identita_fondamentale} and~\eqref{eq: identita_fondamentale_campivettoriali} are established, abstract Hardy-type inequalities (see Theorem~\ref{thm: Hardy lungo Z} and Theorem~\ref{thm: Hardy_completa}) follow from an appropriate and clever choice of the functions involved.

This identity-based strategy was fundamental also in~\cite{YZ25, SYZ25, AYZ25}. We stress that our approach is reminiscent of that introduced by D' Ambrosio~\cite{DA05} and is also closely related to the supersolution method, for which we refer the reader to the comprehensive review by Cazacu~\cite{C220}. 

\medskip
In order to state our main achievements we need to introduce the mathematical setting we will work in.

We will consider $1\leq h \leq N$ vector fields in $\R^N$ of the form
\begin{equation}\label{eq: definizione_X_i}
X_i = \sum_{j=1}^N \sigma_{ij}(x) \frac{\partial}{\partial x_j}, \qquad i = 1, \dots, h,
\end{equation}
where \(\sigma_{ij} \in C(\mathbb{R}^N)\), with \(\frac{\partial}{\partial x_j} \sigma_{ij}\) also belonging to \(C(\mathbb{R}^N)\). Under these assumptions, the formal adjoint of \(X_i\) is well-defined and is given by
\begin{equation}\label{eq: definizione_aggiunto_X_i}
X_i^* = -\sum_{j=1}^N \frac{\partial}{\partial x_j} \big(\sigma_{ij}(x) \cdot \big).
\end{equation}
With these definitions at hands we shall define the horizontal gradient, the Laplacian and the $p$-Laplacian associated to the vector fields defined in~\eqref{eq: definizione_X_i} as follows:
\begin{itemize}
\item The \emph{horizontal gradient} $\nabla_\mathcal{L} = (X_1, \dots, X_h)$ can be written as  
\begin{equation}\label{eq: nabla_L come gradiente euclideo}
\nabla_\mathcal{L} = \sigma \nabla,
\end{equation}
where \(\sigma(x)\in\mathcal{M}_{h \times N}\) is the matrix whose entries are $(\sigma(x))_{ij} = \sigma_{ij}(x),$
and  \(\nabla\) represents the standard Euclidean gradient in \(\mathbb{R}^N\). Related to this operator we have the \emph{horizontal divergence} operator defined as
\begin{equation}\label{eq: definizione_divergenza_L}
\operatorname{div}_\mathcal{L}(\cdot) \coloneqq -\nabla_\mathcal{L}^* \cdot = \operatorname{div}(\sigma^T \cdot),
\end{equation}
where $\nabla_\mathcal{L}^* = (X_1^*, \dots, X_h^*).$
\item The \emph{Laplacian} is defined as $\mathcal{L} = -\nabla_\mathcal{L}^* \cdot \nabla_\mathcal{L}.$
\item The \emph{$p$-Laplacian} is defined as $\mathcal{L}_p u = -\nabla_\mathcal{L}^* \cdot \big(|\nabla_\mathcal{L} u|^{p-2} \nabla_\mathcal{L} u \big).$
\end{itemize}

Now we are in position to state our main results.
The first theorem examines Hardy-type inequalities along a fixed vector field \(Z\). This case is especially important, as it allows us to recover all radial Hardy and Poincaré inequalities by appropriately selecting \(Z\) (see Section~\ref{section:d-radial}).
\begin{theorem}
\label{thm: Hardy lungo Z}
Let \(\Omega \subseteq \mathbb{R}^N\) be an open set. Assume \(p\geq2\) and consider a vector field \(Z \colon \Omega \subseteq \mathbb{R}^N \to \mathbb{R}^h\) and two functions \(V(x) \geq 0\) and \(W(x)\) such that \(W(x) \in L^1_{\text{loc}}(\Omega)\) and \(V(x) |Z|^p \in L^1_{\text{loc}}(\Omega)\). Let \(\lambda > 0\) be a fixed positive constant. Assume that there exists a real-valued function \(h \in C^1(\Omega)\), \(h \neq 0\), that solves the following equation in the distributional sense:
\begin{equation}
\label{eq: hardy_lungo_Z_eqdiff_in_D'}
-\operatorname{div}_\mathcal{L} \big( V |\nabla_\mathcal{L} h \cdot Z|^{p-2} (\nabla_\mathcal{L} h \cdot Z) Z \big) = \lambda W |h|^{p-2} h \quad \text{in } \mathcal{D}'(\Omega),
\end{equation}
which means
\begin{equation}
\label{eq: hardy_lungo_Z_eqdiff_in_D'_integrale}
\int_\Omega  V |\nabla_\mathcal{L} h \cdot Z|^{p-2} (\nabla_\mathcal{L} h \cdot Z)\,\nabla_\mathcal{L} \phi \cdot Z  \, dx = \lambda \int_\Omega  W |h|^{p-2} h\,\phi \, dx \quad \forall \phi \in C^\infty_c(\Omega;\C).
\end{equation}
Then for all $u \in C^\infty_c(\Omega;\C)$ one has
\begin{equation}
\label{eq: hardy_lungo_Z}
\int_\Omega V(x) |\nabla_\mathcal{L} u \cdot Z|^p \, dx \geq \lambda \int_\Omega W(x) |u|^p \, dx.
\end{equation}
Moreover, the function \(u=h\), if admissible, is the candidate maximising function that achieves the equality.
\end{theorem}
\begin{remark}
We explicitly remark that, as will become clear from the proof, we must assume \(V(x) \geq 0\), although no sign assumption has been made on the function \(W(x)\). Moreover, from the assumptions on \(V\), \(W\), and \(Z\) in Theorem~\ref{thm: Hardy lungo Z}, it follows that 
$ V(x) |\nabla_\mathcal{L} h \cdot Z|^{p-2} (\nabla_\mathcal{L} h \cdot Z)\, Z \in L^1_{\text{loc}}(\Omega, \mathbb{C}^h)$
and similarly 
\(
W(x) |h|^{p-2} h \in L^1_{\text{loc}}(\Omega,\C).
\)
Thus, equation \eqref{eq: hardy_lungo_Z_eqdiff_in_D'} is well-defined in the distributional sense.
\end{remark}

In some situations, a single direction determined by \(Z\) may not suffice to fully characterize a Hardy-type inequality. In these cases, it becomes necessary to account for contributions from all directions. The second theorem addresses this broader context, establishing a complete Hardy inequality that incorporates every relevant direction.

\begin{theorem}
\label{thm: Hardy_completa}
Let \(\Omega \subseteq \mathbb{R}^N\) be an open set. Assume \(p\geq2\) and let \(V(x) \geq 0\) and \(W(x)\) be two functions in \(L^1_{\text{loc}}(\Omega)\). Let \(\lambda > 0\) be a fixed positive constant. Assume there exists a real-valued function \(h \in C^1(\Omega\), \(h \neq 0\), that solves the equation
\begin{equation}
\label{eq: eqdiff_completa_in_D'}
-\operatorname{div}_\mathcal{L} \big( V |\nabla_\mathcal{L} h|^{p-2} \nabla_\mathcal{L} h \big) = \lambda W |h|^{p-2} h \quad \text{in } \mathcal{D}'(\Omega,\C).
\end{equation}
Then for all $u \in C^\infty_c(\Omega)$ one has
\begin{equation}
\label{eq: hardy_completa}
\int_\Omega V(x) |\nabla_\mathcal{L} u|^p \, dx \geq \lambda \int_\Omega W(x) |u|^p \, dx.
\end{equation}
Moreover, the function \(u=h\), if admissible, is the candidate maximising function that achieves the equality.
\end{theorem}

\subsection*{Notations}
We collect below some useful notations and properties that will be used throughout the paper.
\begin{itemize} 
\item The Laplacian $\mathcal{L}$ and the $p$-Laplacian $\mathcal{L}_p$ can be expressed in divergence form. More specifically one has 
\begin{equation*}\label{eq: L operatore di divergenza}
\mathcal{L} = \operatorname{div}_\mathcal{L}(\nabla_\mathcal{L}) = \operatorname{div}(A \nabla),
\end{equation*}
and
\begin{equation*}\label{eq: Lp operatore divergenza}
\mathcal{L}_p u = \operatorname{div}_\mathcal{L} \big(|\nabla u|^{p-2} \nabla_\mathcal{L} u \big) = \operatorname{div} \big(|\nabla_\mathcal{L} u|^{p-2} A \nabla u \big),
\end{equation*}
respectively. Here  $A(x)$ is  the symmetric \(N \times N\) matrix defined as \(A(x) = \sigma^T \sigma\).
\item Let $f$ and $g$ two scalar functions, the Leibniz rule for the horizontal gradient states as follows
\begin{equation}\label{eq:Leibniz}
	\nabla_{\mathcal{L}}(fg)=\nabla_{\mathcal{L}}(f) g + f \nabla_{\mathcal{L}}(g).
\end{equation} 
\item Let \(w\) be a function and \(Z\) a vector field, both assumed to be sufficiently regular. Throughout the rest of this work, the following two formulas are fundamental and frequently applied, often implicitly
\begin{equation*}\label{eq: divergenza_L funzione x campo}
\operatorname{div}_\mathcal{L}(\omega Z) = -\nabla_\mathcal{L}^* (\omega Z) = \nabla_\mathcal{L} \omega \cdot Z + \omega \operatorname{div}_\mathcal{L}(Z),
\end{equation*}
\begin{equation*}\label{eq: integrazione per parti}
\int_\Omega \operatorname{div}_\mathcal{L}(Z) \, dx = \int_\Omega \operatorname{div}(\sigma^T Z) \, dx = \int_{\partial \Omega} \sigma^T Z \cdot \nu \, d\mathcal{H}_{N-1} = \int_{\partial \Omega} Z \cdot \sigma \nu \, d\mathcal{H}_{N-1} = \int_{\partial \Omega} Z \cdot \nu_\mathcal{L} \, d\mathcal{H}_{N-1}.
\end{equation*}
Here, \(\nu\) is the outward unit normal to \(\Omega\), and \(\nu_\mathcal{L}\) is defined as \(\nu_\mathcal{L} \coloneqq \sigma \nu\).
\end{itemize}
%

The organisation of the paper is as follows. In Section~\ref{section:proof-main-results} we shall give the proofs of our main results Theorem~\ref{thm: Hardy lungo Z} and Theorem~\ref{thm: Hardy_completa}. The proof of Theorem~\ref{thm: Hardy lungo Z} strongly relies on a scalar algebraic identity (see~\eqref{eq: identita_fondamentale}) which, as a byproduct, provides the following positivity
\begin{equation}\label{eq:scalar-positivity}
	\|f\|_{L^p(\Omega;\C)}^p + (p-1) \|g\|_{L^p(\Omega;\C)}^p -p\Re (|g|^{p-2}g,f)_{L^2(\Omega; \C)}\geq 0,
	\qquad f,g\in L^p(\Omega;\C).
\end{equation}
Once~\eqref{eq:scalar-positivity} is established, then Theorem~\ref{thm: Hardy lungo Z} is obtained through a targeted choice of $f$ and $g$ in~\eqref{eq:scalar-positivity}. The proof of Theorem~\ref{thm: Hardy_completa} follows a similar strategy, the main difference is that the positivity result~\eqref{eq:scalar-positivity} has to be replaced by an analogous result for vector-valued functions (see~\eqref{eq: identita_fondamentale_campivettoriali} below). 
Section~\ref{section:d-radial}, Section~\ref{section:hardy-special} and Section~\ref{section:examples-subelliptic} are devoted to provide some interesting applications of the general results of Section~\ref{section:proof-main-results}. Many of those applications unify, generalise or improve several already existing results.

\section{Proof of main results}
\label{section:proof-main-results}
In order to prove Theorem~\ref{thm: Hardy lungo Z}, we will use the following fundamental algebraic identity, that is a generalisation to complex-valued functions of an identity already obtained for real-valued functions by the second author in~\cite[Prop. 2.1]{DA}.

\begin{proposition}\label{prop:identita_fondamentale}
Let \(\Omega \subseteq \mathbb{R}^N\) be an open set. For every \(f, g \in L^p(\Omega; \C)\) with \(p \geq 2\), the following holds:
\begin{equation}
\label{eq: identita_fondamentale}
\|w(p, f, g)\|^2_{L^2(\Omega;\C)} + \|\widetilde{w}(p, f, g)\|^2_{L^2(\Omega;\C)}  = \|f\|^p_{L^p(\Omega;\C)} + (p - 1)\|g\|^p_{L^p(\Omega;\C)} - p\Re(|g|^{p-2}g, f)_{L^2(\Omega;\C)},
\end{equation}

\textit{where \(w\) and $\widetilde{w}$ are defined as}
\begin{equation}
\label{eq: definizione_w}
w(p, f, g)^2 = p(p-1) \int_0^1 s |s g + (1-s)f|^{p-4} \left[\Re (f-g)(s\overline{g} + (1-s) \overline{f})\right]^2 \, ds,
\end{equation}
and 
\begin{equation}\label{eq: definizione_tildew}
\widetilde{w}(p, f, g)^2 = p\int_0^1 s |s g + (1-s)f|^{p-4} \left[\Im(f\,\overline{g})\right]^2 \, ds.
\end{equation}
\end{proposition}
\begin{remark}
From~\eqref{eq: identita_fondamentale} in Prop.~\ref{prop:identita_fondamentale} one immediately has the lower bound
\begin{equation*}
	 \|f\|^p_{L^p(\Omega;\C)} + (p - 1)\|g\|^p_{L^p(\Omega;\C)} - p\Re(|g|^{p-2}g, f)_{L^2(\Omega;\C)}\geq 0.
\end{equation*}
\end{remark}
\begin{remark}\label{rmk:maximising}
Observe that from the definition of $w$ in~\eqref{eq: definizione_w} and $\widetilde{w}$ in~\eqref{eq: definizione_tildew} it is easy to see that \(\|w(p, f, g)\|^2_{L^2(\Omega;\C)} + \|\widetilde{w}(p, f, g)\|^2_{L^2(\Omega;\C)}= 0\) if and only if \(f = g\) almost everywhere in \(\Omega\).
\end{remark}
\begin{proof}
	The proof of~\eqref{eq: identita_fondamentale} can be obtained generalising the proof provided in~\cite[Prop. 2.1]{DA} to complex-valued functions. Alternatively, the identity~\eqref{eq: identita_fondamentale} can be obtained as a particular case of~\eqref{eq: identita_fondamentale_campivettoriali} in Prop.~\ref{prop:identita_fondamentale_campivettoriali} (\emph{cfr.} Rmk.~\ref{rmk:campivettoriali_to_scalare}).
\end{proof}
Now we are in position to prove our first result Theorem~\ref{thm: Hardy lungo Z}.
\begin{proof}[Proof of Theorem \ref{thm: Hardy lungo Z}]
Given \(u \in C^\infty_c(\Omega,\C)\), we define the auxiliary functions
\begin{equation}\label{eq:f-g}
f= V^{\frac{1}{p}} \nabla_\mathcal{L} u \cdot Z \in L^p(\Omega),
\qquad 
 g = \frac{\nabla_\mathcal{L} h \cdot Z}{h} \,u\, V^{\frac{1}{p}} \in L^p(\Omega).
\end{equation}
Notice that the assumptions \(V(x) \geq 0\) and \(h \neq 0\) in \(\Omega\) are crucial for $f$ and $g$ above to be well defined.

We want to plug $f$ and $g$ defined in~\eqref{eq:f-g} in~\eqref{eq: identita_fondamentale}. We shall compute first the scalar product
\begin{equation*}
	\mathcal{I}:= p \Re(|g|^{p-2}g, f)_{L^2(\Omega)}.
\end{equation*}  
Using that $\nabla_{\mathcal{L}} |u|^p=p |u|^{p-2} \Re (u \nabla_{\mathcal{L}}\overline{u}),$ one has
\begin{equation*}
\begin{split}
\mathcal{I}&=p \int_\Omega \frac{V(x)|\nabla_\mathcal{L} h \cdot Z|^{p-2} \nabla_\mathcal{L} h \cdot Z}{|h|^{p-2} h} |u|^{p-2} \Re (u \nabla_\mathcal{L} \overline{u}) \cdot Z \, dx\\
&=\int_\Omega \nabla_\mathcal{L} |u|^p \cdot \frac{V(x) |\nabla_\mathcal{L} h \cdot Z|^{p-2} \nabla_\mathcal{L} h \cdot Z}{|h|^{p-2} h} Z \, dx\\
&= \int_\Omega \nabla_\mathcal{L} \left( \frac{|u|^p}{|h|^{p-2} h} \right) \cdot \left( V(x) |\nabla_\mathcal{L} h \cdot Z|^{p-2} \nabla_\mathcal{L} h \cdot Z \right) Z \, dx
+ (p-1) \int_\Omega |u|^p \frac{|\nabla_\mathcal{L} h \cdot Z|^p}{|h|^p} V(x) \, dx,
\end{split}
\end{equation*}
where in the last identity we have used the Leibniz rule~\eqref{eq:Leibniz}.

Using~\eqref{eq: hardy_lungo_Z_eqdiff_in_D'_integrale} with the choice $\phi=\frac{|u|^p}{|h|^{p-2} h},$ we deduce
\[
\mathcal{I} = \lambda \int_\Omega |u|^p W(x) \, dx 
+ (p-1) \int_\Omega |u|^p \frac{|\nabla_\mathcal{L} h \cdot Z|^p}{|h|^p} V(x)\, dx.
\]
By combining all terms in the fundamental identity \eqref{eq: identita_fondamentale}, the Hardy-Poincaré inequality naturally follows,
\begin{equation}\label{eq:last}
0 \leq \| w(p, f, g) \|^2_{L^2}   + \| \widetilde{w}(p, f, g) \|^2_{L^2}
= \int_\Omega V(x) |\nabla_\mathcal{L} u \cdot Z|^p \, dx 
- \lambda \int_\Omega W(x) |u|^p  \, dx.
\end{equation}
From~\eqref{eq:last} it is also evident that a function \(u\) achieves equality if and only \(\|w(p, f, g)\|^2_{L^2} + \|\tilde{w}(p, f, g)\|^2_{L^2}\) is zero. As observed in Remark~\ref{rmk:maximising}, this happens if and only if \(f = g\) almost everywhere in \(\Omega\), which gives
\[
f - g = V(x)^{\frac{1}{p}} \nabla_\mathcal{L} u \cdot Z - \frac{\nabla_\mathcal{L} h \cdot Z}{h} \, u \, V(x)^{\frac{1}{p}}
= V(x)^{\frac{1}{p}} h(x) \nabla_\mathcal{L} \left( \frac{u}{h} \right) \cdot Z = 0 \quad \text{a.e. in } \Omega.
\]
Since \(u=h\) is a solution to this equation, the proof of the theorem is complete.
\end{proof}

The proof of Theorem~\ref{thm: Hardy_completa} follows the same steps as the proof of Theorem~\ref{thm: Hardy lungo Z}. However, in this case, we cannot rely on identity~\eqref{eq: identita_fondamentale}, as it is valid only for scalar functions. Instead, we will need the following  vector-valued counterpart of~\eqref{eq: identita_fondamentale}.

\begin{proposition}
\label{prop:identita_fondamentale_campivettoriali}
Let \(\Omega \subseteq \mathbb{R}^N\) be an open set, and let \(f, g \in L^p(\Omega, \mathbb{C}^h)\) with \(p \geq 2\). The following identity holds:
\begin{multline}
\label{eq: identita_fondamentale_campivettoriali}
\int_\Omega |f|^p \, dx + (p-1) \int_\Omega |g|^p \, dx 
- p \Re (|g|^{p-2}g,f)_{L^2(\Omega, \C^h)} \\
=
\frac{1}{p-1} \|\mathrm{w}(p, f, g)(f - g)\|_{L^2(\Omega,\C^h)}^2 + \|\widetilde{\mathrm{w}}(p, f, g)(f - g)\|_{L^2(\Omega,\C^h)}^2,
\end{multline}
where \(\mathrm{w}(p, f, g)\) and \(\widetilde{\mathrm{w}}(p, f, g)\) are defined as 
\begin{equation*}
	|\mathrm{w}(p,f,g)|^2=p(p-1) \int_0^1 s |s g + (1-s) f|^{p-2} \, ds,
\end{equation*}
and
\begin{equation*}
	|\widetilde{\mathrm{w}}(p,f,g)(f-g)|^2=p(p-2) \int_0^1 s |s g + (1-s) f|^{p-4} \big[\Re \langle f - g, s g + (1-s) f \rangle_{\C^h} \big]^2 \, ds.
\end{equation*}
\end{proposition}
%
\begin{remark}\label{rmk:campivettoriali_to_scalare}
When \(h = 1\), identity~\eqref{eq: identita_fondamentale_campivettoriali} reduces to~\eqref{eq: identita_fondamentale}. Indeed, in this case, one has 
\begin{equation*}
\begin{split}
	\frac{1}{p-1}\|\mathrm{w}(p, f, g)(f - g)\|_{L^2(\Omega,\C^h)}^2
	&=p \int_{\Omega} |f-g|^2 \int_0^1 s |sg+(1-s)f|^{p-2}\, ds\, dx\\
	&=p \int_{\Omega} \int_0^1 s |sg+(1-s)f|^{p-4} |(f-g)(s\overline{g}+(1-s)\overline{f})|^2\, ds\, dx\\
	&=p \int_{\Omega} \int_0^1 s |sg+(1-s)f|^{p-4} [\Re ((f-g)(s\overline{g}+(1-s)\overline{f}))]^2\, ds\, dx\\
	&\phantom{=}+ p \int_{\Omega} \int_0^1 s |sg+(1-s)f|^{p-4} [\Im ((f-g)(s\overline{g} +(1-s)\overline{f}))]^2\, ds\, dx.
\end{split}
\end{equation*}
Then the thesis follows from summing this last identity to $\|\widetilde{\mathrm{w}}(p, f, g)(f - g)\|_{L^2(\Omega,\C)}^2.$ 
\end{remark}

The proof of Proposition~\ref{prop:identita_fondamentale_campivettoriali} will be given in Appendix~\ref{app:auxiliary}.

\medskip
We shall continue with the proof of Theorem~\ref{thm: Hardy_completa}.
\begin{proof}[Proof of Theorem~\ref{thm: Hardy_completa}]
Given \(u \in C^\infty_c(\Omega,\C)\), we define the following auxiliary vector-valued functions
\begin{equation*}
	f= V^{\frac{1}{p}} \nabla_\mathcal{L} u \in L^p(\Omega, \mathbb{C}^h),
	\qquad 
	g= \frac{\nabla_\mathcal{L} h}{h} V^{\frac{1}{p}} u \in L^p(\Omega, \mathbb{C}^h).
\end{equation*}
We proceed now in a similar way as in the proof of Theorem~\ref{thm: Hardy lungo Z}. We plug $f$ and $g$ defined above in~\eqref{eq: identita_fondamentale_campivettoriali}, using again that $\nabla_{\mathcal{L}} |u|^p=p |u|^{p-2} \Re (u \nabla_{\mathcal{L}}\overline{u})$ and then using~\eqref{eq: eqdiff_completa_in_D'}, the scalar product p\(\Re(|g|^{p-2} g, f)_{L^2(\Omega, \C^h)}\) can be written as follows  
\begin{equation*}
p\Re (|g|^{p-2} g, f)_{L^2(\Omega, \C^h)} = \lambda \int_\Omega W(x) |u|^p \, dx + (p-1) \int_\Omega |u|^p \frac{|\nabla_\mathcal{L} h|^p}{|h|^p} V(x) \, dx.
\end{equation*}
Plugging the latter into~\eqref{eq: identita_fondamentale_campivettoriali} gives
\begin{equation}\label{eq:preliminary-completa}
\int_\Omega V(x) |\nabla_\mathcal{L} u|^p \, dx - \lambda \int_\Omega W(x) |u|^p \, dx 
=\frac{1}{p-1} \|\mathrm{w}(p, f, g)(f - g)\|_{L^2(\Omega,\C^h)}^2 + \|\widetilde{\mathrm{w}}(p, f, g)(f - g)\|_{L^2(\Omega,\C^h)}^2
\geq 0.
\end{equation}
Hence, inequality~\eqref{eq: hardy_completa} holds.

Notice now that equality holds if and only if both terms on the right-hand side of~\eqref{eq:preliminary-completa} vanish. In particular, \(u\) must satisfy that \(\|\mathrm{w}(p, f, g)(f - g)\|_{L^2}^2 = 0\), which implies \(f = g\) almost everywhere in \(\Omega\). In other words,
\[
f - g = V(x)^{\frac{1}{p}} h(x) \nabla_\mathcal{L} \left( \frac{u}{h} \right) = 0 \quad \text{a.e. in } \Omega.
\]
Since \(u=h\) solves the equation above, then the theorem follows.
\end{proof}

\section{Application I: \(d\)-radial Hardy-type inequalities}
\label{section:d-radial}
The main result of this section is a Hardy-Poincaré type inequality when the weights \(V\) and \(W\) are \emph{radial} with respect to a suitable function \(d\). We will explore specific choices for \(V\) and \(W\), offering concrete examples that include both familiar and previously unexplored inequalities. 

\medskip
We shall fix now the mathematical setting for this section. 

Let \(\Omega \subseteq \mathbb{R}^N\) be an open set, and let \(p \geq 2\). Consider \(d \colon \Omega \to \mathbb{R}\) a non-negative, continuous function such that \(\nabla_\mathcal{L} d \neq 0\) almost everywhere in \(\Omega\). 

Assume that there exists \(x_0 \in \Omega\) such that \(d \in C^\infty(\Omega \setminus \{x_0\})\), and define
\[
R := \sup_{x \in \Omega} d(x), 
\qquad \overline{r} := \inf_{x \in \Omega} d(x),
\]
where \(R\) may take the value \(+\infty\). 

Let  \(V(r) \geq 0\) and \(W(r)\) be two functions defined on the interval \(I=(\overline{r}, R)\). 

Consider \(\alpha, \beta \in \mathbb{R}\) and \(\lambda \geq 0\) being fixed constants, then we shall make the following assumptions:
\begin{itemize}

\item[\((\mathrm{H}.0)\)] The functions \(V(d(x))\) and \(W(d(x)) |\nabla_\mathcal{L} d|^p\) belong to \(L^1_{\text{loc}}(\Omega)\). Furthermore, 
\(W(r) r^{-(\beta - 1)(p - 1)} \in C(I)\).

\item[\((\mathrm{H}.1)\)] There exists a function \(\varphi \neq 0\), \(\varphi \in C^1(I)\) such that 
\((V(r) r^{-(\beta - 1)(p - 1)} |\varphi'|^{p-2} \varphi'(r)) \in C^1(I)\), and \(\varphi\) satisfies the differential equation
\begin{equation}
\label{eq: equazione_differenziale_per_varphi}
\left( \frac{V(r)}{r^{(\beta - 1)(p - 1)}} |\varphi'|^{p-2} \varphi' \right)' 
+ \lambda \frac{W(r)}{r^{(\beta - 1)(p - 1)}} |\varphi|^{p-2} \varphi = 0, \quad \text{for } r \in I.
\end{equation}

\item[\((\mathrm{H}.2)\)] The function \(d\) satisfies
\[
\begin{cases}
\displaystyle\mathcal{L}_p d^\beta = \alpha \delta_{x_0} & \text{in } \mathcal{D}'(\Omega), \quad \text{if } \beta \neq 0, \\[5pt]
\displaystyle\mathcal{L}_p(-\ln d) = \alpha \delta_{x_0} & \text{in } \mathcal{D}'(\Omega), \quad \text{if } \beta = 0.
\end{cases}
\]
\item[\((\mathrm{H}.3)\)] The functions \(\varphi(d)\) and 
\(\displaystyle
\frac{V(d)}{d^{(\beta - 1)(p - 1)}} |\varphi'(d)|^{p-2} \varphi'(d) 
\) belong to \(C^1(\Omega)\). Additionally, we ask the condition
\begin{equation}
\label{eq: condizione_compatibilita_alfa}
\alpha \cdot \left. \frac{V(d)}{d^{(\beta - 1)(p - 1)}} |\varphi'(d)|^{p-2} \varphi'(d) \right|_{x_0} = 0
\end{equation}
to be satisfied.
\end{itemize}
Observe that hypothesis (H.3) requires either \(\alpha = 0\) or 
\(\displaystyle
\left(\frac{V(d)}{d^{(\beta - 1)(p - 1)}} |\varphi'(d)|^{p-2} \varphi'(d)\right)(x_0) = 0.
\)

\medskip
Now we are in position to state the first result of this section.
\begin{theorem}
\label{thm: d-radial_Hardy_poincare_inequality}
Assume \emph{(H.0)--(H.3)}. Then for every \(u \in C^\infty_c(\Omega),\) the following inequality holds
\begin{equation}
\label{eq: d-radial_hardy_poincare_inequality}
\int_\Omega V(d) \left| \nabla_\mathcal{L} u \cdot \frac{\nabla_\mathcal{L} d}{|\nabla_\mathcal{L} d|} \right|^p \, dx 
\geq 
\lambda \int_\Omega W(d) \, |u|^p \, |\nabla_\mathcal{L} d|^p \, dx.
\end{equation}
Moreover, the candidate maximising function is \(u=\varphi(d(x))\).
\end{theorem}

\begin{proof}
Inequality~\eqref{eq: d-radial_hardy_poincare_inequality} would simply follows from~\eqref{eq: hardy_lungo_Z} in Theorem~\ref{thm: Hardy lungo Z} choosing \( Z = \frac{\nabla_\mathcal{L} d}{|\nabla_\mathcal{L} d|} \). thus, we need to verify that hypotheses of Theorem~\ref{thm: Hardy lungo Z} are satisfied. 
It is easy to see that we only need to check the existence of a solution \( h \in C^1(\Omega) \) to equation~\eqref{eq: hardy_lungo_Z_eqdiff_in_D'}. We claim that such solution is given by \(h(x) = \varphi(d(x))\) with $\varphi$ as in (H.1) and (H.3). Therefore we want to show that for every \(\phi \in C^\infty_c(\Omega)\) equation 
\begin{equation}
\label{eq: equazione_da_verificare_distribuzionale}
\int_\Omega 
V(d) \, |\varphi'(d)|^{p-2} \varphi'(d) \, |\nabla_\mathcal{L} d|^{p-2} \, \nabla_\mathcal{L} d \cdot \nabla_\mathcal{L} \phi \, dx 
=
\lambda 
\int_\Omega 
W(d) \, |\nabla_\mathcal{L} d|^p \, |\varphi(d)|^{p-2} \varphi(d) \, \phi(x) \, dx.
\end{equation}
Let us consider the left-hand side of~\eqref{eq: equazione_da_verificare_distribuzionale}, namely 
\begin{equation*}
	\mathcal{I}= \int_\Omega 
V(d) \, |\varphi'(d)|^{p-2} \varphi'(d)\, |\nabla_\mathcal{L} d|^{p-2} \, \nabla_\mathcal{L} d \cdot \nabla_\mathcal{L} \phi \, dx.
\end{equation*}
Multiplying and dividing by $d^{(\beta -1)(p-1)}$ and using Leibniz rule~\eqref{eq:Leibniz} one has 
\begin{multline*}
	\mathcal{I}=
	\int_\Omega d^{(\beta-1)(p-1)} |\nabla_\mathcal{L} d|^{p-2} \nabla_\mathcal{L} d \cdot \nabla_\mathcal{L} \left(\phi\,\, \frac{V(d) \, |\varphi'(d)|^{p-2} \varphi'(d)}{d^{(\beta-1)(p-1)}} \right) dx\\
	- \int_\Omega \phi\, \, d^{(\beta-1)(p-1)} |\nabla_\mathcal{L} d|^{p-2} \nabla_\mathcal{L} d \cdot \nabla_\mathcal{L} \left( \frac{V(d) \, |\varphi'(d)|^{p-2} \varphi'(d)}{d^{(\beta-1)(p-1)}} \right) dx=\mathcal{I}_1 + \mathcal{I}_2.
\end{multline*}
Observe that, using (H.1), $\mathcal{I}_2$ can be rewritten as
\begin{equation*}
\mathcal{I}_2= \lambda \int_\Omega W(d) |\nabla_{\mathcal{L}} d|^p |\varphi(d)|^{p-1} \varphi(d)\, \phi\, dx.
\end{equation*} 
Thus we need to show that $\mathcal{I}_1=0.$ If $\beta\neq 0,$ we observe that
\begin{equation*}
	d^{(\beta-1)(p-1)} |\nabla_{\mathcal{L}}d|^{p-2} \nabla_{\mathcal{L}}d
	= \frac{1}{|\beta|^{p-2}\beta} |\nabla_{\mathcal{L}}d^{\beta}|^{p-2} \nabla_{\mathcal{L}}d^{\beta}.
\end{equation*}
Now integrating by parts and using hypotheses (H.2) and (H.3) one has
\[
\mathcal{I}_1 = -\frac{\phi(x_0)}{|\beta|^{p-2} \beta} \, \alpha \cdot \left( 
\frac{V(d)}{d^{(\beta - 1)(p - 1)}} |\varphi'(d)|^{p-2} \varphi'(d) 
\right)(x_0) = 0.
\]
When \(\beta = 0\), using that 
\begin{equation*}
	d^{-(p-1)} |\nabla_{\mathcal{L}}d|^{p-2} \nabla_{\mathcal{L}}d
	= |\nabla_{\mathcal{L}}(-\ln(d)|^{p-2} \nabla_{\mathcal{L}}(-\ln(d)),
\end{equation*}
integrating by parts and using (H.2) and (H.3) gives
\[
\mathcal{I}_1 = -\phi(x_0) \, \alpha \cdot \left( 
\frac{V(d)}{d^{-(p - 1)}} |\varphi'(d)|^{p-2} \varphi'(d) 
\right)(x_0) = 0.
\]
This concludes the proof.
\end{proof}

\subsection{Classical weighted Hardy inequalities}
In order to recover the more classical form of Hardy inequalities, we shall make further assumptions on the function $d.$ We assume that \(d \in C^\infty(\Omega)\), \(d(x) \neq 0\) for all \(x \in \Omega\) and that there exists \(\beta \in \mathbb{R}\) such that 
\[
\begin{cases} 
\mathcal{L}_p d^\beta = 0 & \text{in } \Omega, \quad \text{if } \beta \neq 0, \\[5pt]
\mathcal{L}_p(-\ln d) = 0 & \text{in } \Omega, \quad \text{if } \beta = 0.
\end{cases}
\]

Now we are in position to prove the following result.
\begin{theorem}
\label{thm: d_radial_weighted_hardy_inequality}
Let \(\theta \in \mathbb{R}\) be a parameter, then for all \(u \in C^\infty_c(\Omega)\) one has
\begin{equation}
\label{eq: d_radial_weighted_hardy_inequality}
\int_\Omega \left| \nabla_\mathcal{L} u \cdot \frac{\nabla_\mathcal{L} d}{|\nabla_\mathcal{L} d|} \right|^p \frac{1}{d^{p(\theta - 1)}} \, dx 
\geq \left| \frac{\beta(p-1) + p(\theta - 1)}{p} \right|^p 
\int_\Omega \frac{|u|^p}{d^{p\theta}} |\nabla_\mathcal{L} d|^p \, dx.
\end{equation}
Moreover, the candidate maximising function is \(u=d^{\frac{\beta(p-1) + p(\theta - 1)}{p}}.\)
\end{theorem}
\begin{proof}
We want to see~\eqref{eq: d_radial_weighted_hardy_inequality} as a special case of~\eqref{eq: d-radial_hardy_poincare_inequality} in Theorem~\ref{thm: d-radial_Hardy_poincare_inequality}.

Consider \(V(r) = r^{-p(\theta - 1)}\) and \(W(r) = r^{-p\theta}\). With this choice, the differential equation~\eqref{eq: equazione_differenziale_per_varphi} in hypothesis (H.1) becomes 
\begin{equation}
\label{eq:classical-preliminary}
\big( r^{-(\beta - 1)(p - 1) - p(\theta - 1)} |\varphi'|^{p-2} \varphi' \big)' 
+ \lambda r^{-(\beta - 1)(p - 1) - p\theta} |\varphi|^{p-2} \varphi = 0.  
\end{equation}
Assume the solution $\varphi$ to be of power-type, namely \(\varphi(r) = r^\gamma\), for a suitable $\gamma$ to be chosen later. Then $\varphi$ satisfies~\eqref{eq:classical-preliminary} if and only if 
\[  
\lambda = -|\gamma|^{p-2} \gamma \big((\gamma - 1)(p - 1) - (\beta - 1)(p - 1) - p(\theta - 1)\big).  
\]  
Maximizing this expression with respect to \(\gamma\), we find that the optimal \(\gamma\) is given by  
\(  \displaystyle
\gamma = \frac{\beta(p - 1) + p(\theta - 1)}{p},  
\) 
and the corresponding \(\lambda\) is  
\[  
\lambda = \left| \frac{\beta(p - 1) + p(\theta - 1)}{p} \right|^p.  
\]  
This concludes the proof.
\end{proof}

\subsubsection{Hardy inequalities related to \(p\)-Harmonic distance functions}
In this section we will provide some application of the previous result Theorem~\ref{thm: d_radial_weighted_hardy_inequality} to specific situations. 

We shall consider first the case of \(d\) being a \(p\)-harmonic function, this requirement corresponds to the parameter choice \(\beta = 1\) in (H.2). The corresponding result is stated below.

\begin{theorem}
\label{thm: hardy_inequality_for_p_armonic_distance_function}
Let \(d \in C^\infty(\Omega)\), with \(d(x) \neq 0\) in $\Omega$ and such that \(\mathcal{L}_p d = 0\) in \(\Omega\). 
Let \(\theta \in \mathbb{R}\) be a parameter, then for all \(u \in C^\infty_c(\Omega)\) one has
\begin{equation}
\label{eq: hardy_inequality_for_p_armonic_distance_function}
\int_\Omega \left| \nabla_\mathcal{L} u \cdot \frac{\nabla_\mathcal{L} d}{|\nabla_\mathcal{L} d|} \right|^p \frac{1}{d^{p(\theta - 1)}} \, dx 
\geq \left| \frac{p\theta - 1}{p} \right|^p  
\int_\Omega \frac{|u|^p}{d^{p\theta}} |\nabla_\mathcal{L} d|^p \, dx.  
\end{equation}  
Moreover, the candidate maximising function is \(u=d^{\frac{p\theta - 1}{p}}\).
\end{theorem}
 
As a direct consequence of this result we obtain the following Hardy inequality in \(\mathbb{R}^2\). We emphasise that, in the specific case \(\theta = 1\) and if the directional gradient is replaced by the full gradient, the resulting estimate already appeared in \cite{DA05}.

\begin{corollary}
Let \(\Omega = \left(-\frac{\pi}{2}, \frac{\pi}{2}\right) \times \mathbb{R}\) and assume \(\nabla_\mathcal{L} = \nabla = (\partial_x, \partial_y)\). Let \(\theta \in \mathbb{R}\), then for all \(u \in C^\infty_c(\Omega)\) one has
\begin{equation}
\label{eq: hardy_in_R2_gradiente_perp}
\int_\Omega 
\frac{|\nabla u \cdot (-\sin x, \cos x)|^2}{\left(\cos x\right)^{2(\theta - 1)}} 
\, \frac{1}{e^{2(\theta - 1)y}} 
\, dx \, dy  
\geq 
\left(\frac{2\theta - 1}{2}\right)^2  
\int_\Omega 
\frac{|u|^2}{\left(\cos x\right)^{2\theta}} 
\, \frac{1}{e^{2(\theta - 1)y}} 
\, dx \, dy.
\end{equation}
Moreover, the constant is sharp.  
\end{corollary}

\begin{proof}
Inequality \eqref{eq: hardy_in_R2_gradiente_perp} directly follows from ~\eqref{eq: hardy_inequality_for_p_armonic_distance_function} in Theorem \ref{thm: hardy_inequality_for_p_armonic_distance_function} by setting \( p = 2 \) and choosing \( d = e^y \cos x \). Thus we only need to prove that the constant is sharp.
According to Theorem~\ref{thm: hardy_inequality_for_p_armonic_distance_function}, the maximising function is given by \( d(x, y)^{\theta - \frac{1}{2}} = e^{(\theta - \frac{1}{2})y} \cos^{\theta - \frac{1}{2}} x \). Although this function is not admissible, this problem can be addressed by introducing a suitable cut-off function.
Let \(\eta \in C^\infty_c(\mathbb{R})\), and let \(f_\varepsilon \in C^\infty_c\big(-\frac{\pi}{2}, \frac{\pi}{2}\big),\) $0 \leq f_\varepsilon(x) \leq 1$ such that
\[
f_\varepsilon(x) =
\begin{cases} 
 1  \qquad -\frac{\pi}{2}< |x| \leq \frac{\pi}{2} \frac{1}{(2\varepsilon + 1)}, \\
 0  \quad \frac{\pi}{2} \frac{1}{(\varepsilon + 1)} \leq |x| < \frac{\pi}{2},
\end{cases} 
\]
moreover assume 
\[
|f'_\varepsilon(x)| \leq c \frac{(\varepsilon + 1)(2\varepsilon + 1)}{\varepsilon} \quad \text{in}\quad \frac{\pi}{2} \frac{1}{(2\varepsilon + 1)} < |x| < \frac{\pi}{2} \frac{1}{(\varepsilon + 1)},
\]
for some constant \(c > 0\).
Define the test function 
\[
u_\varepsilon(x, y) = \cos^{\theta - \frac{1}{2}}(x) f_\varepsilon(x) e^{(\theta - \frac{1}{2})y} \eta(y) \in C^\infty_c(\Omega).
\]
Straightforward computation show that
\begin{multline*}
	\frac{|\nabla u_\varepsilon \cdot (-\sin x, \cos x)|^2}{\cos^{2(\theta - 1)}(x)}  \frac{1}{e^{2(\theta - 1)y}}
	= \left(\frac{2\theta - 1}{2}\right)^2 \frac{f_\varepsilon^2(x)}{\cos(x)} e^y \eta^2(y)\\
+ \cos(x) \sin^2(x) f_\varepsilon'^2(x) e^y \eta^2(y)
+ \cos^3(x) f_\varepsilon^2(x) e^y \eta'^2(y)
- (2\theta - 1) \sin(x) f_\varepsilon(x) f_\varepsilon'(x) e^y \eta^2(y)\\
+ (2\theta - 1) \cos(x) f_\varepsilon^2(x) e^y \eta(y) \eta'(y)
- 2 \cos^2(x) \sin(x) f_\varepsilon(x) f_\varepsilon'(x) e^y \eta(y) \eta'(y).
\end{multline*}
We observe that once integrating the previous identity over \(\Omega\), all terms on the right-hand side of the identity except the first one are bounded by a constant independent of \(\varepsilon\) that we will indicate as $\mathcal{O}(1).$

Furthermore, we easily have that
\[
\int_{-\frac{\pi}{2}}^{\frac{\pi}{2}} \frac{f_\varepsilon^2(x)}{\cos x} \, dx 
\geq \int_{-\frac{\pi}{2} \frac{1}{(2\varepsilon + 1)}}^{\frac{\pi}{2} \frac{1}{(2\varepsilon + 1)}} \frac{1}{\cos x} \, dx 
= 2 \arctanh\left(\sin\left(\frac{\pi}{2} \frac{1}{(2\varepsilon + 1)}\right)\right).
\]
Putting those altogether we conclude
\[
\left(\frac{2\theta - 1}{2}\right)^2 
\leq \frac{\displaystyle\int_\Omega 
\frac{|\nabla u_\varepsilon \cdot (-\sin x, \cos x)|^2}{\left(\cos x\right)^{2(\theta - 1)}}  
\, \frac{1}{e^{2(\theta - 1)y}} 
\, dx \, dy}
{\displaystyle\int_\Omega 
\frac{|u_\varepsilon|^2}{\left(\cos x\right)^{2\theta}}  
\, \frac{1}{e^{2(\theta - 1)y}} 
\, dx \, dy}
\leq \left(\frac{2\theta - 1}{2}\right)^2 
+ \frac{\mathcal{O}(1)}{\arctanh\left(\sin\left(\frac{\pi}{2} \frac{1}{(2\varepsilon + 1)}\right)\right)},
\]
which gives the thesis when $\varepsilon$ goes to $0^+.$
\end{proof}

\subsubsection{Sharp Hardy inequalities related to fundamental solutions}
\label{sezione: Sharp Hardy inequalities related to fundamental solution}
In this section we want to show that under suitable stronger hypotheses on \(d\)  and  \(\mathcal{L}\) one can prove that the constant in the corresponding identity is sharp.  

Specifically, let \(d \in C^\infty(\mathbb{R}^N \setminus \{0\}) \cap C(\mathbb{R}^N)\), with \(d(x) \geq  0\) and \(d(x) = 0\) if and only if \(x = 0\). We assume the existence of a family of dilations \(\delta_\lambda\colon \mathbb{R}^N \to \mathbb{R}^N\), namely
\[
\delta_\lambda(x_1, \ldots, x_N) = (\lambda^{\beta_1} x_1, \ldots, \lambda^{\beta_N} x_N), \quad \beta_j > 0, 
\quad j\in \{1,2,\dots, N\},
\]
with respect to which both the function \(d\) and the vector field \(\nabla_\mathcal{L}\) are homogeneous of degree one. Specifically, they satisfy
\[
\begin{cases} 
d(\delta_\lambda(x)) = \lambda d(x), \\[5pt]
\nabla_\mathcal{L}[u(\delta_\lambda(x))] = \lambda (\nabla_\mathcal{L} u)(\delta_\lambda(x)).
\end{cases}
\]
Moreover, we further assume that
\[
\begin{cases} 
\mathcal{L}_p \, d^{\frac{p - Q}{p - 1}} = 0 \quad \text{in } \mathbb{R}^N \setminus \{0\}, & \text{if } p \neq Q, \\[5pt]
\mathcal{L}_p(-\ln d) = 0 \quad \text{in } \mathbb{R}^N \setminus \{0\}, & \text{if } p = Q,
\end{cases}
\]
where \(Q\) denotes the homogeneous dimension associated with the family of dilations \(\delta_\lambda\) which is given by
\(
Q = \beta_1 + \beta_2 + \cdots + \beta_N.
\)

\medskip
We shall show the following result.

\begin{theorem}
\label{thm: sharp_d_radial_weighted_hardy_inequaity}
Let \(\theta \in \mathbb{R}\), and let \(\Omega \subseteq \mathbb{R}^N \setminus \{0\}\), then for all $u \in C^\infty_c(\Omega)$ one has
\begin{equation}
\label{eq: sharp_d_radial_weighted_hardy_inequaity}
\int_\Omega \left| \nabla_\mathcal{L} u \cdot \frac{\nabla_\mathcal{L} d}{|\nabla_\mathcal{L} d|} \right|^p \frac{1}{d^{p(\theta - 1)}} \, dx 
\geq \left| \frac{Q - p \theta}{p} \right|^p \int_\Omega \frac{|u|^p}{d^{p \theta}} |\nabla_\mathcal{L} d|^p \, dx.
\end{equation}
Moreover, if \(\Omega \cup \{0\}\) contains a neighborhood of the origin, then the constant \(\left| \frac{Q - p \theta}{p} \right|^p\) is sharp.
\end{theorem}

\begin{remark}
If \(\Omega = \mathbb{R}^N \setminus \{0\}\), this result is already found in \cite{DA}.
\end{remark}

The proof of the theorem relies on the following technical results obtained in~\cite{DA} and~\cite{DA2}. We refer to these works for the corresponding proofs and further details.

The first lemma establishes that the function $d$ behaves like a proper distance functions: its ‘‘balls" are compact.
\begin{lemma}[Lemma 2.2 in \cite{DA2}]
\label{lemma: Lemma 2.2 in DA2}
The set \(B_R^d=\{x \in \mathbb{R}^N \mid d(x) \leq R\}\) is compact for every \(R > 0\).
\end{lemma}
The second result tells that integration over $d$-balls and spheres exhibits exact scaling laws analogous to those of the corresponding euclidean sets.
\begin{lemma}[Lemma 2.3 in \cite{DA2}]
\label{lemma: Lemma 2.3 in DA2}
Let \(\alpha \geq 0\) and define \(\lambda_\alpha := \int_{B_1^d} |\nabla_\mathcal{L} d|^\alpha \, dx\). Then, we have
\begin{equation}\label{eq:d-ball-sphere}
\int_{B_R^d} |\nabla_\mathcal{L} d|^\alpha \, dx = \lambda_\alpha R^Q, 
\qquad \text{and} \qquad 
\int_{S_R^d} \frac{|\nabla_\mathcal{L} d|^\alpha}{|\nabla d|} \, d\mathcal{H}^{N-1} = Q \lambda_\alpha R^{Q-1},
\end{equation}
where $B_R^d:=\{x\in \R^N\mid d(x)\leq R\}$ and $S_R^d:=\{x\in \R^N\mid d(x)= R\}.$ 
\end{lemma}


We will also need the following inequality.
\begin{lemma}[Lemma 2.1 in~\cite{DA}]
\label{lemma: Lemma 2.1 in DA}
For every \(p \geq 2\), there exists a positive constant \(c_p > 0\) such that for every \(a, b \in \mathbb{R}\),
\begin{equation}\label{eq:p-power-ineq}
|a + b|^p \leq |a|^p + c_p \left(|a|^{p-1} |b| + |b|^p\right).
\end{equation}
\end{lemma}

Now we are in position to prove Theorem~\ref{thm: sharp_d_radial_weighted_hardy_inequaity}.

\begin{proof}[Proof of Theorem~\ref{thm: sharp_d_radial_weighted_hardy_inequaity}]
Inequality~\eqref{eq: sharp_d_radial_weighted_hardy_inequaity} follows directly from inequality~\eqref{eq: d_radial_weighted_hardy_inequality} in Theorem~\ref{thm: d_radial_weighted_hardy_inequality}. Therefore we only need to show that the constant in~\eqref{eq: sharp_d_radial_weighted_hardy_inequaity} is sharp. 

\medskip
We first consider the case $\Omega = \mathbb{R}^N \setminus \{0\},$ the general case will be addressed later.

Let us consider \(u_\varepsilon(x) = d(x)^{\frac{p \theta - Q}{p}} g_\varepsilon(d(x)) \in C^\infty_c(\mathbb{R}^N \setminus \{0\})\), where $g_\varepsilon$ represents a cut-off function explicitly defined in~\eqref{eq: definizione g_veps}. The function $u_\varepsilon$ represents the candidate maximising function for~\eqref{eq: sharp_d_radial_weighted_hardy_inequaity} (ref. Theorem~\ref{thm: d_radial_weighted_hardy_inequality}). 

Using coarea formula and the second identity in~\eqref{eq:d-ball-sphere} of Lemma~\ref{lemma: Lemma 2.3 in DA2}, one has 
\begin{equation*}
\begin{split}
\int_{\mathbb{R}^N} \frac{|u_\varepsilon|^p}{d^{p \theta}} |\nabla_\mathcal{L} d|^p \, dx 
&= \int_{\mathbb{R}^N} d(x)^{-Q} |g_\varepsilon(d(x))|^p |\nabla_\mathcal{L} d|^p \, dx
= \int_0^{+\infty} r^{-Q} |g_\varepsilon(r)|^p \int_{\{d = r\}} \frac{|\nabla_\mathcal{L} d|^p}{|\nabla d|} \, d\mathcal{H}_{N-1} \, dr\\
&=\lambda_p Q \left\{-\ln(4 \varepsilon^2) + \mathcal{O}(1)\right\},
\end{split}
\end{equation*}
where in the last identity we have used~\eqref{eq:g_eps1}.

By applying inequality~\eqref{eq:p-power-ineq} in Lemma~\ref{lemma: Lemma 2.1 in DA} and, again, the second identity in~\eqref{eq:d-ball-sphere}, one gets
\begin{equation*}
\begin{split}
\int_{\mathbb{R}^N} \left| \nabla_\mathcal{L} u_\varepsilon \cdot \frac{\nabla_\mathcal{L} d}{|\nabla_\mathcal{L} d|} \right|^p \frac{1}{d^{p(\theta - 1)}} \, dx 
&= \int_{\mathbb{R}^N} \left| 
\left(\frac{p \theta - Q}{p}\right) d^{-\frac{Q}{p}} g_\varepsilon(d) |\nabla_\mathcal{L} d| 
+ d^{1 - \frac{Q}{p}} g'_\varepsilon(d) |\nabla_\mathcal{L} d| 
\right|^p \, dx\\ 
&
\leq \left| \frac{Q - p \theta}{p} \right|^p \lambda_p Q \int_0^{+\infty} r^{-1} |g_\varepsilon(r)|^p \, dr\\ 
&\quad + b_1 \int_0^{+\infty} |g_\varepsilon(r)|^{p-1} |g'_\varepsilon(r)| \, dr
+ b_2 \int_0^{+\infty} r^{p-1} |g'_\varepsilon(r)|^p \, dr,
\end{split}
\end{equation*}
for some positive constants \( b_1 \) and \( b_2 \), independent of \( \varepsilon \). This and~\eqref{eq:g_eps1}-\eqref{eq:g_eps3} give
\[
\int_{\mathbb{R}^N} \left| \nabla_\mathcal{L} u_\varepsilon \cdot \frac{\nabla_\mathcal{L} d}{|\nabla_\mathcal{L} d|} \right|^p \frac{1}{d^{p(\theta - 1)}} \, dx 
\leq \left| \frac{Q - p \theta}{p} \right|^p \lambda_p Q \left\{-\ln(4 \varepsilon^2) + \mathcal{O}(1)\right\}.
\]
Therefore
\[
\left| \frac{Q - p \theta}{p} \right|^p \leq 
\frac{\displaystyle \int_{\mathbb{R}^N} \left| \nabla_\mathcal{L} u_\varepsilon \cdot \frac{\nabla_\mathcal{L} d}{|\nabla_\mathcal{L} d|} \right|^p \frac{1}{d^{p(\theta - 1)}} \, dx}
{\displaystyle \int_{\mathbb{R}^N} \frac{|u_\varepsilon|^p}{d^{p \theta}} |\nabla_\mathcal{L} d|^p \, dx}
\leq 
\frac{\displaystyle -\left| \frac{Q - p \theta}{p} \right|^p \ln(4 \varepsilon^2) + \mathcal{O}(1)}
{\displaystyle -\ln(4 \varepsilon^2) + \mathcal{O}(1)}.
\]
Letting $\varepsilon$ go to $0^+$, gives the claim in the case \(\Omega = \mathbb{R}^N \setminus \{0\}\).

\medskip
We consider now the more general case of $\Omega$ such that \(\Omega \cup \{0\}\) contains a neighborhood of the origin. The proof in this case will strongly rely on the invariance under dilations of~\eqref{eq: sharp_d_radial_weighted_hardy_inequaity}. 
Let \(c(\Omega)\) denote the optimal constant in~\eqref{eq: sharp_d_radial_weighted_hardy_inequaity}.  
By hypothesis, there exists \(r > 0\) sufficiently small such that
\[
\label{eq: catene_inclusioni_sharp_hardy_foundamntal_solution}
B^d_r \setminus \{0\} \subseteq \Omega \subseteq \mathbb{R}^N \setminus \{0\},
\]
as a result
\[
\left| \frac{Q - p \theta}{p} \right|^p = c(\mathbb{R}^N \setminus \{0\}) \leq c(\Omega) \leq c(B^d_r \setminus \{0\}).
\]
Let $u\in C^\infty_c(\R^N\setminus \{0\}),$ then there exists $R>0$ sufficiently large such that $u\in C^\infty_c(B_R^d \setminus \{0\}),$ this implies $c(B_R^d \setminus \{0\})\leq c(\R^N\setminus \{0\}).$
The previous facts together give 
\[
\left| \frac{Q - p \theta}{p} \right|^p = c(\mathbb{R}^N \setminus \{0\}) \leq c(\Omega) \leq c(B^d_r \setminus \{0\}) = c(B^d_R \setminus \{0\}) \leq c(\mathbb{R}^N \setminus \{0\}) = \left| \frac{Q - p \theta}{p} \right|^p,
\]
where in the previous we have used that, due to the invariance of inequality~\eqref{eq: sharp_d_radial_weighted_hardy_inequaity} under dilations, one has
\[
c(B^d_R \setminus \{0\}) = c(B^d_1 \setminus \{0\}), \quad \forall\, R > 0.
\]
This concludes the proof in the general case.
\end{proof}

\begin{remark}[Criticality of the Hardy Inequality]
We shall discuss here the criticality (\emph{ref}.~\cite{Pinchover07}) of inequality~\eqref{eq: sharp_d_radial_weighted_hardy_inequaity} in the special case $\theta=1.$ As a first result we show that the operator 
\begin{equation*}
	\mathcal{L}_p - \Bigl|\tfrac{Q-p}{p}\Bigr|^p \frac{1}{d^p} |\nabla_\mathcal{L} d|^p
\end{equation*}
is critical in $\R^N.$ This result was already proved in~\cite{CKNL24} for the standard $p$-Laplacian (in~\cite{CKNL24} they impose the restriction  $p<d,$ in our case this is not needed as we are working in \( C_c^\infty(\mathbb{R}^N \setminus \{0\}) ,\) so no integrability conditions at the origin are needed). 
\begin{proposition}\label{prop:criticality}
Let $W \in L^1_{\mathrm{loc}}(\mathbb{R}^N \setminus \{0\})$ a non-negative potential.   
Assume that for all $u\in C^\infty_c(\R^N\setminus \{0\})$ one has
\begin{equation}\label{eq:plus_W}
\int_{\mathbb{R}^N} 
\Big| \nabla_\mathcal{L} u \cdot \tfrac{\nabla_\mathcal{L} d}{|\nabla_\mathcal{L} d|} \Big|^p dx
\geq
 \Bigl|\frac{Q-p}{p}\Bigr|^p  \int_{\mathbb{R}^N} \frac{|u|^p}{d^p} |\nabla_\mathcal{L} d|^p \, dx
+ \int_{\mathbb{R}^N} W \, |u|^p |\nabla_\mathcal{L} d|^p \, dx,
\end{equation}
then $W=0$ for almost every $x\in \R^N.$
\end{proposition}

\begin{proof}
Fix $R>0.$ We define a Lipschitz cit-off function $\psi_R$ as follows: $\psi_R(r)=0$ if $0 \leq r < \tfrac{1}{R^2}$ and if $r\geq R^2,$ $\psi_R(r)=1$ if $\tfrac{1}{R} \leq r \leq R,$ moreover $\psi_R(r)=2 + \tfrac{\ln r}{\ln R}$ if $\tfrac{1}{R^2} \leq r \leq \tfrac{1}{R}$ and $\psi_R(r)=2 - \tfrac{\ln r}{\ln R}$ if $ R \leq r \leq  R^2.$ 

A direct computation shows that
\begin{equation*}
\int_0^{+\infty} r \, |\psi'_R(r)|^2 \, dr 
= \frac{2}{\ln R},
\qquad 
\text{and}
\qquad 
\int_0^{+\infty} \psi_R(r)\,\psi'_R(r)\, dr = 0.
\end{equation*}

We consider the test function
\begin{equation*}
u_R = d^{-\frac{Q-p}{p}}\,\psi_R(d),
\end{equation*}
with $\psi_R$ as above. An explicit computation shows that
\begin{equation}\label{eq:bound-lnR}
\int_{\mathbb{R}^N}\Bigl|\nabla_{\mathcal L}u_R\cdot \frac{\nabla_{\mathcal L} d}{|\nabla_{\mathcal L} d|}\Bigr|^p\,dx
-
\left|\frac{Q-p}{p}\right|^p
\int_{\mathbb{R}^N}\frac{|u_R|^p}{d^p}\,|\nabla_{\mathcal L} d|^p\,dx
= \mathcal{O} \Bigl(\frac{1}{\ln R}\Bigr)\xrightarrow[R\to+\infty]{} 0.
\end{equation}
Using Fatou’s lemma, applying~\eqref{eq:plus_W} with $u=u_R$ and using bound~\eqref{eq:bound-lnR} one gets
\begin{multline*}
0 \le \int_{\mathbb{R}^N} W\, d^{\,p-Q}\,|\nabla_{\mathcal L} d|^p\,dx
\le
\liminf_{R\to+\infty}\int_{\mathbb{R}^N} W\,|u_R|^p\,|\nabla_{\mathcal L} d|^p\,dx\\
\le
\liminf_{R\to+\infty} \Biggl[
\int_{\mathbb{R}^N}\Bigl|\nabla_{\mathcal L}u_R\cdot \frac{\nabla_{\mathcal L} d}{|\nabla_{\mathcal L} d|}\Bigr|^p\,dx
-\left|\frac{Q-p}{p}\right|^p
\int_{\mathbb{R}^N}\frac{|u_R|^p}{d^p}\,|\nabla_{\mathcal L} d|^p\,dx
\Biggr]
= 0.
\end{multline*}
Hence $W=0$ almost everywhere in $\mathbb{R}^N.$
\end{proof}

Despite the result in Proposition~\ref{prop:criticality}, one can show that $W$ can be strictly larger than zero on a subset of $\R^N$ of positive measure. 
In order to prove this result we will make use of the following lower bound which is a generalisation to complex-valued functions of a previous result obtained in~\cite[Thm. 3.1.]{DA2}.
\begin{lemma}
Let $D\subset \R^N$ and let $p\geq 2.$ then there exist a constant $c_p\in [2^{-p}, p 2^{-p}]$ such that for any $f,g\in L^p(D; \C)$ one has 
\begin{equation}\label{eq:Lorenzo2}
	\|f\|_{L^p(D; \C)}^p + (p-1) \|g\|_{L^p(D;\C)}^p -p\Re (|g|^{p-2} g, f)_{L^2(D;\C)}
	\geq  c_p \|f-g\|_{L^p(D;\C)}^p.
\end{equation}
\end{lemma}
\begin{proof}
We use identity~\eqref{eq: identita_fondamentale_campivettoriali} in Proposition~\ref{prop:identita_fondamentale_campivettoriali} with $h=1$.  
From~\cite[Prop.~3.1]{DA2}, there exists $c_p \in [2^{-p},\, p\,2^{-p}]$ such that
\[
|\mathrm{w}(p,f,g)|^2 \;\ge\; c_p\,|f-g|^{p-2}
\]
and, by similar arguments to those used in~\cite[Prop.~3.1]{DA2}, one obtains the analogous lower bound
\[
|\widetilde{\mathrm{w}}(p,f,g)(f-g)|^2 \;\ge\; c_p\,\frac{p-2}{p-1}\,|f-g|^{p}
\]
From this, the claim follows.
\end{proof}
%
Now we are in position to state and prove the aforementioned result.
\begin{proposition}\label{prop: controesempio citicita hardy}
Let $Q \ge 1$ and $p \ge 2$, then there exists
$W \in L^1_{\mathrm{loc}}(\mathbb{R}^N \setminus \{0\}),$ with $W>\Big|\frac{Q-p}{p}\Big|^p \frac{1}{d^p}$ on a set of positive measure, such that for all $ u \in C_c^\infty(\mathbb{R}^N\setminus\{0\})$ one has
\begin{equation*}
\int_{\mathbb{R}^N}\Bigl|\nabla_\mathcal{L} u \cdot 
\frac{\nabla_\mathcal{L} d}{|\nabla_\mathcal{L} d|}\Bigr|^p dx
\ge
\int_{\mathbb{R}^N} W\,|u|^p\,|\nabla_\mathcal{L} d|^p dx.
\end{equation*}
\end{proposition}

\begin{proof}
Using~\eqref{eq: identita_fondamentale} with the choice 
\begin{equation*}
f = \nabla_\mathcal{L} u \cdot \frac{\nabla_\mathcal{L} d}{|\nabla_\mathcal{L} d|},
\qquad \text{and} \qquad
g = -\left(\frac{Q-p}{p}\right)\frac{u}{d}\,|\nabla_\mathcal{L} d|,
\end{equation*}
integrating by parts and using~\eqref{eq:Lorenzo2} we have
\begin{equation}\label{eq:former} 
\int_{\mathbb{R}^N}\Bigl|\nabla_\mathcal{L} u \cdot 
\frac{\nabla_\mathcal{L} d}{|\nabla_\mathcal{L} d|}\Bigr|^p dx
-
\left|\frac{Q-p}{p}\right|^p 
\int_{\mathbb{R}^N}\frac{|u|^p}{d^p}\,|\nabla_\mathcal{L} d|^p dx
\geq  c_p \int_{\mathbb{R}^N} d^{-(Q-p)} 
\Bigl|\nabla_\mathcal{L}\phi\cdot \frac{\nabla_\mathcal{L} d}{|\nabla_\mathcal{L} d|}\Bigr|^p dx,
\end{equation}
where we have defined $\phi = d^{\frac{Q-p}{p}}u$.
Moreover, from Theorem~\ref{thm: d-radial_Hardy_poincare_inequality}
with $\widetilde V=d^{-(Q-p)}$ and $\widetilde W=d^{-(Q-p)}\tfrac{1-d}{d}$, since the function $\varphi(r)=e^{-r}$ solves the associated differential equation
\[
\bigl( r^{Q-1}\,\widetilde V\,|\varphi'|^{p-2}\varphi' \bigr)' 
+ (p-1)r^{Q-1}\widetilde W(r)|\varphi|^{p-2}\varphi = 0,
\]
one has 
\begin{equation}\label{eq:latter}
\int_{\mathbb{R}^N} d^{-(Q-p)} 
\Bigl|\nabla_\mathcal{L}\phi\cdot \frac{\nabla_\mathcal{L} d}{|\nabla_\mathcal{L} d|}\Bigr|^p dx
\ge
(p-1)\int_{\mathbb{R}^N} d^{-(Q-p)}\frac{1-d}{d}\,|\phi|^p\,|\nabla_\mathcal{L} d|^p dx,
\end{equation}
for every $\phi \in C_c^\infty(\mathbb{R}^N\setminus\{0\}).$

Putting~\eqref{eq:former} and~\eqref{eq:latter} together, we conclude that
\[
\int_{\mathbb{R}^N}\Bigl|\nabla_\mathcal{L} u \cdot 
\frac{\nabla_\mathcal{L} d}{|\nabla_\mathcal{L} d|}\Bigr|^p dx
-
\left|\frac{Q-p}{p}\right|^p 
\int_{\mathbb{R}^N}\frac{|u|^p}{d^p}\,|\nabla_\mathcal{L} d|^p dx
\ge
c_p(p-1) \int_{\mathbb{R}^N}\frac{1-d}{d}\,|u|^p\,|\nabla_\mathcal{L} d|^p dx.
\]
Therefore the claim holds with
\[
W = \left|\frac{Q-p}{p}\right|^p \frac{1}{d^p} 
+ c_p(p-1)\left(\frac{1-d}{d}\right),
\]
observing that $W>\left|\frac{Q-p}{p}\right|^p \frac{1}{d^p}$ when $d<1.$
\end{proof}

\end{remark}

\subsection{Cylindrical Hardy inequalities}
\label{sezione: sharp_cylindrical_hardy_inequality}
In this section, we consider a special form of the matrix $\sigma$ defining $\nabla_{\mathcal{L}}$ (see~\eqref{eq: nabla_L come gradiente euclideo}). More precisely, assume that there exists \(m \in \mathbb{N}\), with \(1 \leq m \leq h\), such that \(\sigma \in \mathcal{M}_{h \times N}\) takes the form
\begin{equation}
\label{eq: matrice sigma cilindrica}
\sigma = 
\begin{pmatrix} 
\mathcal{I}_m & \sigma_1 \\ 
0 & \sigma_2 
\end{pmatrix},
\end{equation}
where \(\sigma_1 \in \mathcal{M}_{m \times (N-m)}\), \(\sigma_2 \in \mathcal{M}_{(h-m) \times (N-m)}\) and \(\mathcal{I}_m\) denotes the identity matrix in \(\mathbb{R}^m\).

Moreover, we impose the condition
\begin{equation}
\label{eq: condizione su sigma_1 per vandermonde0}
\frac{\partial}{\partial y_j} (\sigma_1)_{ij} = 0, \quad \forall\, i = 1, \dots, m, \quad \forall\, j = 1, \dots, N-m.
\end{equation}
Explicit relevant examples of operators of this form can be found in Section~\ref{section:examples-subelliptic}.

We shall consider \(\mathbb{R}^N=\mathbb{R}^m \times \mathbb{R}^{N-m}\), with coordinates \((x, y)\), where \(x \in \mathbb{R}^m\) and \(y \in \mathbb{R}^{N-m}\). 

In addition, we assume that there exist a family of dilations
\begin{equation*}
\delta_\lambda(x,y) = (\lambda^{\beta_1} x_1, \ldots, \lambda^{\beta_m} x_m, \lambda^{\eta_1} y_1, \ldots, \lambda^{\eta_{N-m}} y_{N-m}),
\end{equation*}
with respect to which \(\nabla_\mathcal{L}\) is homogeneous of degree one. From the specific structure of \(\sigma\) in \eqref{eq: matrice sigma cilindrica}, it follows that \(\beta_1 = \cdots = \beta_m = 1\), namely
\begin{equation}\label{eq:eta-dil}
\delta_\lambda(x_1, \ldots, x_m, y_1, \ldots, y_{N-m}) = (\lambda x_1, \ldots, \lambda x_m, \lambda^{\eta_1} y_1, \ldots, \lambda^{\eta_{N-m}} y_{N-m}).
\end{equation}
From~\eqref{eq: matrice sigma cilindrica} and~\eqref{eq: condizione su sigma_1 per vandermonde0} we see that \(\mathcal{L}_p |x| = \Delta_p |x|\), where \(\Delta_p\) denotes the standard Euclidean \(p\)-Laplacian. Furthermore one has
\begin{equation}
\label{eq: L_p e funzione |x|}
\begin{cases}
\displaystyle\mathcal{L}_p |x|^{\frac{p-m}{p-1}} = 0 & \text{in } \mathbb{R}^m \setminus \{0\} \times \mathbb{R}^{N-m}, \quad \text{if } p \neq m, \\[5pt]
\displaystyle\mathcal{L}_p(-\ln |x|) = 0 & \text{in } \mathbb{R}^m \setminus \{0\} \times \mathbb{R}^{N-m}, \quad \text{if } p = m.
\end{cases}
\end{equation}

Given $r>0,$ and $l\in \mathbb{N},$ we shall denote $B_r^l$ as the Euclidean ball in $\R^l$ centered at the origin and with radius $r.$
Now we are in position to state the main result of the section.
\begin{theorem}
\label{thm: sharp_cilindrical_hardy_inequality}
Let \(\sigma\) be as in~\eqref{eq: matrice sigma cilindrica}, fix \(\theta \in \mathbb{R}\), and assume \(\Omega \subseteq \mathbb{R}^m \setminus \{0\} \times \mathbb{R}^{N-m}\), then for all \(u \in C^\infty_c(\Omega)\) one has 
\begin{equation}
\label{eq: sharp_cilindrical_hardy_inequality}
\int_\Omega \left| \nabla_\mathcal{L} u \cdot \left( \frac{x}{|x|}, 0_{h-m} \right) \right|^p \frac{1}{|x|^{p(\theta - 1)}} \, dx \, dy 
\geq \left| \frac{m - p \theta}{p} \right|^p \int_\Omega \frac{|u(x, y)|^p}{|x|^{p \theta}} \, dx \, dy.
\end{equation}
Moreover, if \(B^m_r \setminus \{0\} \times B^{N-m}_r \subseteq \Omega\) for some \(r > 0\) and \(x \cdot (\sigma_1 y) = 0\), then the constant \(\big| \frac{m - p \theta}{p} \big|^p\) is sharp.
\end{theorem}

\begin{remark}
This result was already obtained by D'Ambrosio in~\cite{DA05}, nevertheless in this work the general condition \(x \cdot (\sigma_1 y) = 0\), which ensures the sharpness of the constant, was not explicitly identified.
\end{remark}

\begin{proof}[Proof of Theorem~\ref{thm: sharp_cilindrical_hardy_inequality}]
Inequality~\eqref{eq: sharp_cilindrical_hardy_inequality} follows directly from Theorem~\ref{thm: d_radial_weighted_hardy_inequality} by choosing \(d(x, y) = |x|\). Moreover, it is known that the candidate maximizing function is \( |x|^{\frac{p \theta - m}{p}} \).

To establish the sharpness of the constant, we consider first the specific case \(\Omega = \mathbb{R}^m \setminus \{0\} \times \mathbb{R}^{N-m}\). We define a family of approximating functions 
\[
u_\varepsilon(x, y) = |x|^{\frac{p \theta - m}{p}} g_\varepsilon(|x|) \varphi(|y|)\in C^\infty_c(\mathbb{R}^m \setminus \{0\} \times \mathbb{R}^{N-m}),
\]
where \(g_\varepsilon(r)\) is the cut-off function defined in~\eqref{eq: definizione g_veps}, and \(\varphi(r) \) is any function in \(C^\infty_c(\mathbb{R}_+)\).

Observe that \(
\nabla_\mathcal{L} |y| = \left( \sigma_1 \frac{y}{|y|}, \sigma_2 \frac{y}{|y|} \right),
\)
therefore, since we are assuming $x \cdot (\sigma_1 y)=0,$ we have 
\begin{equation}\label{eq:der-phi-zero}
\nabla_\mathcal{L} |y| \cdot \left( \frac{x}{|x|}, 0_{h-m} \right) = \frac{1}{|x| |y|} x \cdot (\sigma_1 y) = 0.
\end{equation}

Using~\eqref{eq:der-phi-zero} and~\eqref{eq:p-power-ineq} in Lemma~\ref{lemma: Lemma 2.1 in DA}, we find
\begin{equation*}
\begin{split}
\int_{\mathbb{R}^N} \left| \nabla_\mathcal{L} u_\varepsilon \cdot \left( \frac{x}{|x|}, 0_{h-m} \right) \right|^p \frac{1}{|x|^{p(\theta - 1)}} \, dx \, dy 
&= \int_{\mathbb{R}^N} \left| \left( \frac{p \theta - m}{p} \right) |x|^{-\frac{m}{p}} g_\varepsilon(|x|) \varphi(|y|) + |x|^{1 - \frac{m}{p}} g'_\varepsilon(|x|) \varphi(|y|) \right|^p \, dx \, dy \vspace{0.2cm}\\
& \leq \left| \frac{p \theta - m}{p} \right|^p \int_{\mathbb{R}^{N-m}} |\varphi|^p \, dy \int_{\mathbb{R}^m} |x|^{-m} |g_\varepsilon(|x|)|^p \, dx\\
 &\phantom{\leq } + \left| \frac{p \theta - m}{p} \right|^{p-1} c_p \int_{\mathbb{R}^{N-m}} |\varphi|^p \, dy \int_{\mathbb{R}^m} |x|^{1-m} |g_\varepsilon(|x|)|^{p-1} |g'_\varepsilon(|x|)| \, dx\\
 &\phantom{\leq }
 + c_p \int_{\mathbb{R}^{N-m}} |\varphi|^p \, dy \int_{\mathbb{R}^m} |x|^{p-m} |g'_\varepsilon(|x|)|^p \, dx.
\end{split}
\end{equation*}
Using spherical coordinates in \(\mathbb{R}^m\), we have
\begin{equation*}
\begin{split}
\int_{\mathbb{R}^N} \left| \nabla_\mathcal{L} u_\varepsilon \cdot \left( \frac{x}{|x|}, 0_{h-m} \right) \right|^p \frac{1}{|x|^{p(\theta - 1)}} \, dx \, dy
&\leq \left| \frac{p \theta - m}{p} \right|^p |\mathbb{S}^{m-1}| 
\int_{\mathbb{R}^{N-m}} |\varphi|^p \, dy \left \{\int_0^{+\infty} r^{-1} |g_\varepsilon(r)|^p \, dr\right.\\
&\left.\phantom{\leq }+ b_1 \int_0^{+\infty}  |g_\varepsilon(r)|^{p-1} |g'_\varepsilon(r)| \, dr 
+ b_2 \int_0^{+\infty} r^{p-1} |g'_\varepsilon(r)|^p \, dr \right\},
\end{split}
\end{equation*}
where \(b_1\) and \(b_2\) are suitable positive constants independent of \( \varepsilon \).

Using~\eqref{eq:g_eps1}-\eqref{eq:g_eps3}, we have
\[
\int_{\mathbb{R}^N} \left| \nabla_\mathcal{L} u_\varepsilon \cdot \left( \frac{x}{|x|}, 0_{h-m} \right) \right|^p \frac{1}{|x|^{p(\theta - 1)}} \, dx \, dy = \left| \frac{m - p \theta}{p} \right|^p |\mathbb{S}^{m-1}| \int_{\mathbb{R}^{N-m}} |\varphi|^p \, dy \left\{ -\ln(4 \varepsilon^2)  + \mathcal{O}(1)\right\}.
\]
Similarly, it is straightforward to verify that
\[
\int_{\mathbb{R}^N} \frac{|u_\varepsilon(x, y)|^p}{|x|^{p \theta}} \, dx \, dy = |\mathbb{S}^{m-1}| \int_{\mathbb{R}^{N-m}} |\varphi|^p \, dy \left\{ -\ln(4 \varepsilon^2) + \mathcal{O}(1) \right\}.
\]
Altogether gives
\[
\left| \frac{m - p \theta}{p} \right|^p 
\leq 
\frac{\displaystyle \int_{\mathbb{R}^N} \left| \nabla_\mathcal{L} u_\varepsilon \cdot \left( \frac{x}{|x|}, 0_{h-m} \right) \right|^p \frac{1}{|x|^{p(\theta - 1)}} \, dx \, dy }{\displaystyle \int_{\mathbb{R}^N} \frac{|u_\varepsilon(x, y)|^p}{|x|^{p \theta}} \, dx \, dy}
\leq 
\frac{\displaystyle - \left| \frac{m - p \theta}{p} \right|^p \ln(4 \varepsilon^2) + \mathcal{O}(1)}{\displaystyle -\ln(4 \varepsilon^2) + \mathcal{O}(1)},
\]
and the claim follows by letting $\varepsilon$ go to $ 0^+.$ 

To consider the general case, let \(c(\Omega)\) denote the optimal constant in \eqref{eq: sharp_cilindrical_hardy_inequality} and let \(r > 0\) be such that \(B^m_r \setminus \{0\} \times B^{N-m}_r \subseteq \Omega\), and define \(B^1_{s^{\eta_i}} = \{ |y_i| \leq s^{\eta_i} \}\). 
We can further assume that there exists \(s > 0\), sufficiently small, such that 
\[
B^*_s := B^m_s \setminus \{0\} \times B^1_{s^{\eta_1}} \times \cdots \times B^1_{s^{\eta_{N-m}}} \subseteq \Omega \subseteq \mathbb{R}^m \setminus \{0\} \times \mathbb{R}^{N-m},
\]
where $\eta_1, \cdots, \eta_{N-m}$ where introduced in~\eqref{eq:eta-dil}.
We have
\[
\left| \frac{m - p \theta}{p} \right|^p = c(\mathbb{R}^m \setminus \{0\} \times \mathbb{R}^{N-m}) \leq c(\Omega) \leq c(B^*_s).
\]
Given \(u \in C^\infty_c(\mathbb{R}^m \setminus \{0\} \times \mathbb{R}^{N-m})\), then, in particular, $u\in C^\infty_c(B^*_R),$ for a sufficiently large $R>0.$ Therefore $c(B^*_R)\leq c(\mathbb{R}^m \setminus \{0\} \times \mathbb{R}^{N-m}).$ The previous fact together give 
\begin{equation*}
	\left| \frac{m - p \theta}{p} \right|^p = c(\mathbb{R}^m \setminus \{0\} \times \mathbb{R}^{N-m}) \leq c(\Omega) \leq c(B^*_s)
	=c(B^*_R)
	\leq c(\mathbb{R}^m \setminus \{0\} \times \mathbb{R}^{N-m})
	=\left| \frac{m - p \theta}{p} \right|^p.
\end{equation*} 
Here we have used that, since \(\delta_\lambda(B^*_1) = B^*_\lambda\), then \(c(B^*_s) = c(B^*_1)\) for all \(s > 0.\)
\end{proof}

\subsection{Logarithmic Hardy inequalities}
In this section, we simplify, unify, and extend some results from D’Ambrosio \cite{DA05} and Kalaman and Yessirkegenov \cite{KY}.

\subsubsection{Sharp Logarithmic Hardy inequalities}
We shall make similar assumptions as in Section~\ref{sezione: Sharp Hardy inequalities related to fundamental solution}.
Specifically, let \(d \in C^\infty(\mathbb{R}^N \setminus \{0\}) \cap C(\mathbb{R}^N)\), with \(d(x) \geq  0\) and \(d(x) = 0\) if and only if \(x = 0\). We assume the existence of a family of dilations \(\delta_\lambda\colon \mathbb{R}^N \to \mathbb{R}^N\), namely
\[
\delta_\lambda(x_1, \ldots, x_N) = (\lambda^{\beta_1} x_1, \ldots, \lambda^{\beta_N} x_N), \quad \beta_j > 0, 
\quad j\in \{1,2,\dots, N\},
\]
with respect to which both the function \(d\) and the vector field \(\nabla_\mathcal{L}\) are homogeneous of degree one. 

Consider a domain \(\Omega \subseteq \mathbb{R}^N \setminus \{0\}\) and suppose \(R := \sup_\Omega d(x) < +\infty\).
Moreover, we further assume that
\[
\begin{cases} 
\mathcal{L}_p \, d^{\frac{p - Q}{p - 1}} = 0 \quad \text{in } \Omega, & \text{if } p \neq Q, \\[5pt] 
\mathcal{L}_p(-\ln d) = 0 \quad \text{in } \Omega, & \text{if } p = Q,
\end{cases}
\]
where \(Q\), as usual, denotes the homogeneous dimension associated with the family of dilations \(\delta_\lambda\) which is given by
\(
Q = \beta_1 + \beta_2 + \cdots + \beta_N.
\)

\medskip
Now we can state the main result of this section.
\begin{theorem}
\label{thm: sharp_logaritmic_hardy_inequality}
Let \(\theta \in \mathbb{R}\), and let \(\Omega\) be as above, then for all \(u \in C^\infty_c(\Omega)\) one has
\begin{equation}
\label{eq: sharp_logaritmic_hardy_inequality}
\int_\Omega \left| \nabla_\mathcal{L} u \cdot \frac{\nabla_\mathcal{L} d}{|\nabla_\mathcal{L} d|} \right|^p \left( \ln \frac{R}{d} \right)^{\theta + p} \frac{1}{d^{Q-p}} \, dx 
\geq \left| \frac{\theta + 1}{p} \right|^p \int_\Omega \frac{|u|^p}{d^Q} \left( \ln \frac{R}{d} \right)^\theta |\nabla_\mathcal{L} d|^p \, dx.
\end{equation}
Moreover, if \(\Omega \cup \{0\}\) contains a neighborhood of the origin, then the constant is sharp.
\end{theorem}

\begin{remark}
In the classical Hardy inequality~\eqref{eq: sharp_d_radial_weighted_hardy_inequaity}, the optimal constant vanishes when \(Q = p \theta\), in this case the inequality is said to be critical. However, the previous result shows that, introducing a logarithmic weight, produces a non-trivial inequality. More specifically, for \(Q = p \theta\) one obtains from~\eqref{eq: sharp_logaritmic_hardy_inequality}
\[
\int_\Omega \left| \nabla_\mathcal{L} u \cdot \frac{\nabla_\mathcal{L} d}{|\nabla_\mathcal{L} d|} \right|^p \frac{1}{d^{p(\theta - 1)}} \left( \ln \frac{R}{d} \right)^{\theta + p} \, dx  
\geq \left| \frac{Q + p}{p^2} \right|^p \int_\Omega \frac{|u|^p}{d^{p \theta}} \left( \ln \frac{R}{d} \right)^\theta |\nabla_\mathcal{L} d|^p \, dx.
\]
\end{remark}

\begin{proof}[Proof of Theorem~\ref{thm: sharp_logaritmic_hardy_inequality}]
Inequality~\eqref{eq: sharp_logaritmic_hardy_inequality} simply follows from~\eqref{eq: d-radial_hardy_poincare_inequality} in Theorem~\ref{thm: d-radial_Hardy_poincare_inequality}. Indeed, the hypotheses of Theorem~\ref{thm: d-radial_Hardy_poincare_inequality} are  satisfied,  as a matter of fact, a nontrivial solution to~\eqref{eq: equazione_differenziale_per_varphi} in (H.1) is given by
\begin{equation*}
	\varphi(r) = \left(\ln\frac{R}{r}\right)^{\!-\frac{\theta + 1}{p}}, 
	\qquad
\lambda = \left|\frac{\theta + 1}{p}\right|^p.
\end{equation*}
Moreover, the candidate maximizing function is given by \( \left(\ln\frac{R}{d}\right)^{-\frac{\theta + 1}{p}}\). We only need that the constant in~\eqref{eq: sharp_logaritmic_hardy_inequality} is sharp.

We consider first the case \(\Omega = B^d_R \setminus \{0\}.\)

Let \(u_\varepsilon(x)\) be the sequence of test functions defined as 
\[
u_\varepsilon(x) = \left(\ln\frac{R}{d}\right)^{-\frac{\theta+1}{p}} g_\varepsilon\left(\ln\frac{R}{d}\right), \quad u_\varepsilon \in C^\infty_c(B^d_R \setminus \{0\}),
\]
where $g_\varepsilon$ is the cut-off function introduced in~\eqref{eq: definizione g_veps}.
Using the coarea formula, Lemma~\ref{lemma: Lemma 2.3 in DA2} and the asymptotic~\eqref{eq:g_eps-log1} from Proposition~\ref{prop: asintotico_g_veps_logaritmico}, we obtain
\[
\int_{B^d_R} \frac{|u_\varepsilon|^p}{d^Q} \left(\ln \frac{R}{d}\right)^\theta |\nabla_\mathcal{L} d|^p \, dx = \lambda_p Q \big(\!-\ln(4 \varepsilon^2) + \mathcal{O}(1)\big).
\]
Similarly, using~\eqref{eq:p-power-ineq} in Lemma~\ref{lemma: Lemma 2.1 in DA} and the asyptotics~\eqref{eq:g_eps-log1}-\eqref{eq:g_eps-log3} in Proposition~\ref{prop: asintotico_g_veps_logaritmico}, we find
\[
\int_{B^d_R} \left|\nabla_\mathcal{L} u_\varepsilon \cdot \frac{\nabla_\mathcal{L} d}{|\nabla_\mathcal{L} d|}\right|^p \left(\ln \frac{R}{d}\right)^{\theta + p} \frac{1}{d^{Q-p}} \, dx 
\leq \left|\frac{\theta+1}{p}\right|^p \lambda_p Q \big(-\ln(4 \varepsilon^2) + \mathcal{O}(1)\big).
\]
Altogether one gets 
\[
\left|\frac{\theta+1}{p}\right|^p \leq \frac{\displaystyle \int_{B^d_R} \left|\nabla_\mathcal{L} u_\varepsilon \cdot \frac{\nabla_\mathcal{L} d}{|\nabla_\mathcal{L} d|}\right|^p \left(\ln \frac{R}{d}\right)^{\theta + p} \frac{1}{d^{Q-p}} \, dx}{\displaystyle \int_{B^d_R} \frac{|u_\varepsilon|^p}{d^Q} \left(\ln \frac{R}{d}\right)^\theta |\nabla_\mathcal{L} d|^p \, dx } 
\leq \left|\frac{\theta+1}{p}\right|^p\frac{-\ln(4 \varepsilon^2) + \mathcal{O}(1)}{-\ln(4 \varepsilon^2) + \mathcal{O}(1)}.
\]
Then the claim follows taking the limit $\varepsilon$ to $0^+.$
To handle the general case, let \(c(\Omega)\) denote the optimal constant in \eqref{eq: sharp_logaritmic_hardy_inequality}. By the homogeneity of the inequality with respect to the family of dilations \(\delta_\lambda\), we have 
\[
c(B^d_R \setminus \{0\}) = c(B^d_1 \setminus \{0\}), \quad \forall R > 0.
\]
Let \(R = \sup_\Omega d(x)\). From the inclusion \(\Omega \subseteq B^d_R \setminus \{0\}\), it follows that 
\(
c(B^d_1 \setminus \{0\}) = c(B^d_R \setminus \{0\}) \leq c(\Omega).
\)
Conversely, let \(r > 0\) be sufficiently small such that \(B^d_r \setminus \{0\} \subseteq \Omega\). Then 
\(
 c(\Omega) \leq c(B^d_r \setminus \{0\}) = c(B^d_1 \setminus \{0\}).
\)
\end{proof}

\subsubsection{Cylindrical Logarithmic Hardy inequalities}
\label{sezione: Cylindrical Logarithmic Hardy Inequality}

We shall show the following result.

\begin{theorem}
\label{thm: Cylindrical logaritmic hardy inequality}
Let \(R > 0\), fix \(\,\theta \in \mathbb{R}\), and assume \(\Omega = B^m_R \setminus \{0\} \times \mathbb{R}^{N-m}\), then for all \(u \in C^\infty_c(\Omega)\) one has
\begin{equation}
\label{eq: cylindrical logaritmic hardy inequality}
\int_\Omega \left| \nabla_\mathcal{L} u \cdot \left(\frac{x}{|x|}, 0_{h-m}\right) \right|^p \frac{1}{|x|^{m-p}} \left(\ln\left(\frac{R}{|x|}\right)\right)^{\theta+p} \, dx \, dy 
\geq \left|\frac{\theta+1}{p}\right|^p \int_\Omega \frac{|u|^p}{|x|^m} \left(\ln\left(\frac{R}{|x|}\right)\right)^\theta \, dx \, dy.
\end{equation}
Moreover, if \(x \cdot (\sigma_1 y) = 0\), the constant \(\left|\frac{\theta+1}{p}\right|^p\) is sharp.
\end{theorem}

\begin{proof}
Inequality~\eqref{eq: cylindrical logaritmic hardy inequality} follows directly from~\eqref{eq: cylindrical logaritmic hardy inequality} of Theorem~\ref{thm: d-radial_Hardy_poincare_inequality} by choosing \(d(x, y) = |x| \) and
\begin{equation*}
	\varphi(r) = \left(\ln \frac{R}{r}\right)^{-\frac{\theta+1}{p}},
	\quad \text{and} \quad 
	\lambda = \left|\frac{\theta+1}{p}\right|^p
\end{equation*} 
as the solution of the associated differential equation~\eqref{eq: equazione_differenziale_per_varphi}. Moreover $\varphi(|x|)$ as defined above is also the candidate maximiser for~\eqref{eq: cylindrical logaritmic hardy inequality}.
To prove the sharpness of the constant, we consider the family of test functions
\[
u_\varepsilon = \left(\ln \frac{R}{|x|}\right)^{-\frac{\theta+1}{p}} g_\varepsilon\left(\ln \frac{R}{|x|}\right) \phi(|y|), \quad u_\varepsilon \in C^\infty_c(\Omega),
\]
where \(g_\varepsilon\) is defined in \eqref{eq: definizione g_veps}, and \(\phi(r)\) is an arbitrary function in \(C^\infty_c(\mathbb{R}_+)\).

Using spherical coordinates in $\R^m$ and~\eqref{eq:g_eps-log1} in Proposition~\ref{prop: asintotico_g_veps_logaritmico}, we deduce that  
\[
\displaystyle \int_\Omega \frac{|u_\varepsilon|^p}{|x|^m} \left(\ln\left(\frac{R}{|x|}\right)\right)^\theta \, dx \, dy = |\mathbb{S}^{m-1}| \int_{\mathbb{R}^{N-m}} |\phi(|y|)|^p \, dy \left(-\ln(4 \varepsilon^2) + \mathcal{O}(1)\right).
\]
Similarly, using~\eqref{eq:p-power-ineq} in Lemma~\ref{lemma: Lemma 2.1 in DA} and the assumption \(x \cdot (\sigma_1 y) = 0\),  and also the asyptotics~\eqref{eq:g_eps-log1}-\eqref{eq:g_eps-log3} in Proposition~\ref{prop: asintotico_g_veps_logaritmico}, we find
\begin{multline*}
\int_\Omega \left| \nabla_\mathcal{L} u_\varepsilon \cdot \left(\frac{x}{|x|}, 0_{N-m}\right) \right|^p \frac{1}{|x|^{m-p}} \left(\ln\left(\frac{R}{|x|}\right)\right)^{\theta+p} \, dx \, dy \\
\leq |\mathbb{S}^{m-1}| \int_{\mathbb{R}^{N-m}} |\phi(|y|)|^p \, dy \left(-\left|\frac{\theta+1}{p}\right|^p \ln(4 \varepsilon^2) + \mathcal{O}(1)\right).
\end{multline*}
Altogether gives
\[
\left|\frac{\theta+1}{p}\right|^p \leq 
\frac{\displaystyle \int_\Omega \left| \nabla_\mathcal{L} u_\varepsilon \cdot \left(\frac{x}{|x|}, 0_{N-m}\right) \right|^p \frac{1}{|x|^{m-p}} \left(\ln\left(\frac{R}{|x|}\right)\right)^{\theta+p} \, dx \, dy}
{\displaystyle \int_\Omega \frac{|u_\varepsilon|^p}{|x|^m} \left(\ln\left(\frac{R}{|x|}\right)\right)^\theta \, dx \, dy}
\leq 
\frac{\displaystyle -\left|\frac{\theta+1}{p}\right|^p \ln(4 \varepsilon^2) + \mathcal{O}(1)}{\displaystyle -\ln(4 \varepsilon^2) + \mathcal{O}(1)},
\]
and the sharpeness of the constant follows taking the limit as $\varepsilon$ goes to $0^+.$
\end{proof}

\subsection{Gaussian Hardy inequalities}
\label{sezione: Gaussian-type Hardy inequality}
In this section we shall work under assumptions closely related to those in Section~\ref{sezione: Sharp Hardy inequalities related to fundamental solution}, but here we will only consider the case \(\Omega = \mathbb{R}^N\) and \(\Omega = \mathbb{R}^N \setminus \{0\}\). For clarity, we briefly summarize the relevant conditions. Let \(d \in C^\infty(\mathbb{R}^N \setminus \{0\}) \cap C(\mathbb{R}^N)\), with \(d(x) \geq 0\) and \(d(x) = 0\) if and only if \(x = 0\).  
Additionally, suppose there exists a family of dilations defined as  
\(
\delta_\lambda(x_1, \ldots, x_N) = (\lambda^{\beta_1} x_1, \ldots, \lambda^{\beta_N} x_N),\) \(\beta_j \geq 1,
\)  
under which both the function \(d\) and the vector field \(\nabla_\mathcal{L}\) are homogeneous of degree one.  

Finally, if \(0 \in \Omega\), the function \(d\) is required to satisfy
\[
\begin{cases} 
\mathcal{L}_p \, d^{\frac{p - Q}{p - 1}} = l_p \delta_0 \quad \text{in } \mathcal{D}'(\Omega), & \text{if } p \neq Q, \\[5pt] 
\mathcal{L}_p(-\ln d) = l_p \delta_0 \quad \text{in } \mathcal{D}'(\Omega), & \text{if } p = Q, 
\end{cases}
\]
for some constant \(l_p \in \mathbb{R}\). Conversely, if \(0 \notin \Omega\), we simply set \(l_p = 0\).
Here \(Q = \beta_1 + \beta_2 + \cdots + \beta_N\) denotes the homogeneous dimension associated with the family of dilations. Note also that \(Q \geq 1\), as we have assumed here that \(\beta_j \geq 1\).

\medskip
Now we are in position to state the two main results of this section.
\begin{theorem}
\label{thm: hardy gaussiana tipo A}
Let \(\alpha\) and \(\beta\) be constants such that \(\beta > 0\) and \(\alpha \geq 2\), then for all \(u \in C^\infty_c(\mathbb{R}^N)\) one has
\begin{multline}
\label{eq: hardy gaussiana tipo A}
\int_{\mathbb{R}^N} \left| \nabla_\mathcal{L} u \cdot \frac{\nabla_\mathcal{L} d}{|\nabla_\mathcal{L} d|} \right|^p \displaystyle{ e^{-\tfrac{ d^{ \alpha}}{ \beta}}} dx 
\geq
\left( \frac{\alpha}{p \beta} \right)^p \int_{\mathbb{R}^N} |u|^p d^{p( \alpha - 1)} |\nabla_\mathcal{L} d|^p \displaystyle{ e^{-\tfrac{ d^{ \alpha}}{ \beta}}} dx  \\
- \left( \frac{\alpha}{p \beta} \right)^{p-1} \left(\alpha  (p - 1) + Q - p\right) \int_{\mathbb{R}^N} |u|^p d^{ \alpha (p - 1) - p} |\nabla_\mathcal{L} d|^p \displaystyle{ e^{-\tfrac{ d^{ \alpha}}{ \beta}}} dx.
\end{multline}
The candidate maximizing function is \(\displaystyle{ e^{\frac{ d^{ \alpha}}{p \beta}}}.\)
Moreover, this inequality is sharp in the sense that it does not hold if the right-hand side is multiplied by any constant \(c > 1\). 
\end{theorem}
\begin{remark}
Notice that the constant \((\alpha (p-1) + Q - p)\) in~\eqref{eq: hardy gaussiana tipo A} is positive.
\end{remark}

\begin{remark}[Standard Gaussian measure]
\label{rmk:standard-gaussian}
To obtain Hardy inequalities with respect to the standard Gaussian measure \(d\gamma = (2\pi)^{-\frac{N}{2}} e^{-\frac{d^2}{2}} \, dx\), it suffices to set \(\alpha = \beta = 2\) in Theorem~\ref{thm: hardy gaussiana tipo A}. Thus, one gets the following result. 

\begin{theorem}\label{thm:gen-BCT}
	Let $d\gamma$ be the standard Gaussian measure, \emph{i.e.} $d\gamma = (2\pi)^{-\frac{N}{2}} e^{-\frac{d^2}{2}} \, dx,$ then for  all $u \in C^\infty_c(\mathbb{R}^N)$ one has
	\begin{equation*}
		\int_{\mathbb{R}^N} \left| \nabla_\mathcal{L} u \cdot \frac{\nabla_\mathcal{L} d}{|\nabla_\mathcal{L} d|} \right|^p d\gamma 
\geq 
\frac{1}{p^p} \int_{\mathbb{R}^N} |u|^p d^p |\nabla_\mathcal{L} d|^p d\gamma 
- \frac{Q + p - 2}{p^{p-1}} \int_{\mathbb{R}^N} |u|^p d^{p-2} |\nabla_\mathcal{L} d|^p d\gamma
	\end{equation*}
	This inequality is sharp, as stated in Theorem~\ref{thm: hardy gaussiana tipo A}, with the maximizing function \(e^{\frac{d^2}{2p}}\).
\end{theorem}

We emphasise that Theorem~\ref{thm:gen-BCT} generalises to the $L^p$ framework and the subelliptic setting the following inequality proved by Brandolini, Chiacchio and Trombetti in~\cite{BCT}
\begin{equation*}
\int_{\mathbb{R}^N} \left( |\nabla u|^2 + \frac{N}{2} |u|^2 \right) e^{-\frac{|x|^2}{2}} \, dx 
\geq \frac{1}{4} \int_{\mathbb{R}^N} |x|^2 |u|^2 e^{-\frac{|x|^2}{2}} \, dx, \quad \forall u \in W^{1,2}(\mathbb{R}^N, e^{-\frac{|x|^2}{2}}dx).
\end{equation*}
\end{remark}

\medskip
The next result is a power-type weighted version of Theorem~\ref{thm: hardy gaussiana tipo A} above.

\begin{theorem}
\label{thm: hardy gaussiana tipo B}
Fix \(\theta \in \mathbb{R}\) and let \(\alpha, \beta\) be constants such that $\beta>0$ and \(\alpha \geq 2,\) then for all \(u \in C^\infty_c(\mathbb{R}^N \setminus \{0\})\) one has
\begin{multline}
\label{eq: hardy gaussiana tipo B}
\int_{\mathbb{R}^N} \left| \nabla_\mathcal{L} u \cdot \frac{\nabla_\mathcal{L} d}{|\nabla_\mathcal{L} d|} \right|^p \frac{1}{d^{p(\theta-1)}} \displaystyle{ e^{-\tfrac{ d^{ \alpha}}{ \beta}}} \, dx 
\geq
\left| \frac{Q - p \theta}{p} \right|^p \int_{\mathbb{R}^N} \frac{|u|^p}{d^{p \theta}} |\nabla_\mathcal{L} d|^p \displaystyle{ e^{-\tfrac{ d^{ \alpha}}{ \beta}}} \, dx\\ 
- \frac{\alpha}{\beta} \left(\frac{Q - p \theta}{p}\right) \left| \frac{Q - p \theta}{p} \right|^{p-2} \int_{\mathbb{R}^N} \frac{|u|^p}{d^{p\theta - \alpha}} |\nabla_\mathcal{L} d|^p \displaystyle{ e^{-\tfrac{ d^{ \alpha}}{ \beta}}} \, dx.
\end{multline}
The candidate maximizing function is \(d^{-\frac{Q - p \theta}{p}}\).
Furthermore, if \(Q \neq p \theta\), the inequality is sharp in the sense that it does not hold if the right-hand side is multiplied by any constant \(c > 1\). 
\end{theorem}

\begin{remark}
Notice that if \(Q - p \theta < 0\) in Theorem~\ref{thm: hardy gaussiana tipo B}, both terms on the right-hand side of~\eqref{eq: hardy gaussiana tipo B} are positive. Thus, this shows that replacing the classical Lebesgue measure \(dx\) in Hardy inequality~\eqref{eq: sharp_d_radial_weighted_hardy_inequaity} with the Gaussian measure \( e^{- d^{ \alpha}/ \beta} \, dx\) improves the inequality.
\end{remark}

\begin{remark}
Similarly to Remark~\ref{rmk:standard-gaussian} one has the following weighted version of Theorem~\ref{thm:gen-BCT}.
\begin{theorem}
Let $d\gamma$ be the standard Gaussian measure, \emph{i.e.} $d\gamma = (2\pi)^{-\frac{N}{2}} e^{-\frac{d^2}{2}} \, dx,$ then for  all $u \in C^\infty_c(\mathbb{R}^N)$ one has
\[
\int_{\mathbb{R}^N} \left| \nabla_\mathcal{L} u \cdot \frac{\nabla_\mathcal{L} d}{|\nabla_\mathcal{L} d|} \right|^p \frac{1}{d^{p(\theta-1)}} d\gamma 
\geq 
\left| \frac{Q - p \theta}{p} \right|^p \int_{\mathbb{R}^N} \frac{|u|^p}{d^{p \theta}} |\nabla_\mathcal{L} d|^p d\gamma 
- \frac{Q - p \theta}{p} \left| \frac{Q - p \theta}{p} \right|^{p-2} \int_{\mathbb{R}^N} \frac{|u|^p}{d^{p \theta - 2}} |\nabla_\mathcal{L} d|^p d\gamma.
\]
\end{theorem}
The candidate maximizing function is \(d^{-\frac{Q - p \theta}{p}}\).
Furthermore, if \(Q \neq p \theta\), the inequality is sharp, as stated in Theorem~\ref{thm: hardy gaussiana tipo B}.
\end{remark}

Now we begin with the proofs of Theorem~\ref{thm: hardy gaussiana tipo A} and Theorem~\ref{thm: hardy gaussiana tipo B}.

\begin{proof}[Proof of Theorem~\ref{thm: hardy gaussiana tipo A}]
Inequality~\eqref{eq: hardy gaussiana tipo A} follows directly from Theorem~\ref{thm: d-radial_Hardy_poincare_inequality} by choosing 
\[
V(r) = e^{-\tfrac{ r^{ \alpha}}{ \beta}}, \qquad \text{and} \qquad
W(r) = \left(r^{p( \alpha - 1)} - \frac{p \beta}{\alpha} \big( \alpha (p-1) + Q - p\big) r^{ \alpha (p-1) - p}\right) e^{-\tfrac{ r^{ \alpha}}{ \beta}},
\]
and observing that a solution to the corresponding differential equation~\eqref{eq: equazione_differenziale_per_varphi} is given by
\begin{equation*}
	\varphi(r) = e^{\tfrac{ r^{ \alpha}}{ p \beta}},
	\qquad 
	\lambda = \left(\frac{\alpha}{p \beta}\right)^p.
\end{equation*}
Moreover, since \(0 \in \Omega = \mathbb{R}^N\), we need to verify the compatibility condition~\eqref{eq: condizione_compatibilita_alfa} in (H.3) too. A straightforward calculation shows that this condition is always satisfied.

To show sharpness of the constant in~\eqref{eq: hardy gaussiana tipo A}, we consider the following sequence of approximating functions 
\[
u_\varepsilon(x) = \displaystyle{ e^{\tfrac{ d^{ \alpha}}{ p\beta}}} g_\varepsilon(d) \in C^\infty_c(\mathbb{R}^N \setminus \{0\}),
\]
where \(g_\varepsilon\) as in~\eqref{eq: definizione g_veps}.
We define 
\begin{equation}\label{eq:h-eps}
	h(\varepsilon):= \int_0^{+\infty} r^{Q - 1 +  \alpha (p-1)} |g_\varepsilon(r)|^p \, dr.
\end{equation}
Notice that $h(\varepsilon)$ as defined in~\eqref{eq:h-eps} goes to infinity as $\varepsilon$ goes to $0^+.$

Using the coarea formula, together with Lemma~\ref{lemma: Lemma 2.3 in DA2} and Lemma~\ref{lemma: Lemma 2.1 in DA}, we obtain

\begin{multline*}
\int_{\mathbb{R}^N} 
\left| \nabla_\mathcal{L} u_\varepsilon \cdot 
\frac{\nabla_\mathcal{L} d}{|\nabla_\mathcal{L} d|} \right|^p 
e^{-\tfrac{ d^{\alpha}}{\beta}} \, dx 
\leq \lambda_p \Bigg\{
\left(\frac{\alpha}{p \beta}\right)^p
h(\varepsilon) 
+ c \int_0^{+\infty} r^{Q - 1 + (p - 1)(\alpha - 1)} 
|g_\varepsilon(r)|^{p - 1} |g'_\varepsilon(r)| \, dr \\
+ c \int_0^{+\infty} r^{Q - 1} |g'_\varepsilon(r)|^p \, dr
\Bigg\},
\end{multline*}
where \( c > 0 \) denotes a generic constant that may change from line to line.

Using the definition of $g_\varepsilon$ one simply has
\begin{equation*}
\begin{split}
\int_0^{+\infty} r^{Q - 1 + (p - 1)(\alpha - 1)} 
|g_\varepsilon(r)|^{p - 1} |g'_\varepsilon(r)| \, dr 
&\leq \frac{c}{\varepsilon} 
\int_\varepsilon^{2\varepsilon} r^{Q - 1 + (p - 1)(\alpha - 1)} \, dr
+ c \varepsilon 
\int_{\frac{1}{2\varepsilon}}^{\frac{1}{\varepsilon}} 
r^{Q - 1 + (p - 1)(\alpha - 1)} \, dr \\[4pt]
&\leq c \!\left( 1 + 
\frac{1}{\varepsilon^{Q + (p - 1)(\alpha - 1) - 1}} \right).
\end{split}
\end{equation*}
Moreover,
\[
h(\varepsilon) \geq \int_{2\varepsilon}^{\frac{1}{2\varepsilon}} r^{Q - 1 + p(\alpha - 1)} \, dr
= \frac{c}{\varepsilon^{Q + p(\alpha - 1)}}.
\]
These two estimates together give 
\[
\frac{
\displaystyle \int_0^{+\infty} r^{Q - 1 + (p - 1)(\alpha - 1)} |g_\varepsilon(r)|^{p - 1} |g'_\varepsilon(r)| \, dr
}{
h(\varepsilon)
}
\leq c \, \varepsilon^{\alpha} \longrightarrow 0 
\quad \text{as } \varepsilon \to 0^+.
\]
Similarly
\[
\frac{
\displaystyle \int_0^{+\infty} r^{Q - 1} |g'_\varepsilon(r)|^p \, dr
}{
h(\varepsilon)
}
\leq c \left( \varepsilon^{2Q + p(\alpha - 2)} + \varepsilon^{p \alpha} \right)
\longrightarrow 0 
\quad \text{as } \varepsilon \to 0^+.
\]
Here we used \(\alpha \ge 2\).
Consequently, one gets
\[
\int_{\mathbb{R}^N} 
\left| \nabla_\mathcal{L} u_\varepsilon \cdot \frac{\nabla_\mathcal{L} d}{|\nabla_\mathcal{L} d|} \right|^p 
e^{-\tfrac{ d^{ \alpha}}{ \beta}} \, dx 
\leq \lambda_p \left\{
\left(\frac{\alpha}{p \beta}\right)^p h(\varepsilon) + o\big(h(\varepsilon)\big)
\right\}.
\]

Proceeding in the same way, we obtain
\[
\int_{\mathbb{R}^N} |u_\varepsilon|^p d^{p( \alpha - 1)} |\nabla_\mathcal{L} d|^p \displaystyle{ e^{-\tfrac{ d^{ \alpha}}{ \beta}}} \, dx 
- \left( \frac{p \beta}{\alpha} \right) \left(\alpha  (p - 1) + Q - p\right) \int_{\mathbb{R}^N} |u_\varepsilon|^p d^{ \alpha (p - 1) - p} |\nabla_\mathcal{L} d|^p \displaystyle{ e^{-\tfrac{ d^{ \alpha}}{ \beta}}} \, dx =
\]
\[
= \lambda_p \left\{h(\varepsilon) + o(h(\varepsilon))\right\}.
\]
Taking the ratio of these two expressions and letting \(\varepsilon \to 0^+\), we conclude.
\end{proof}

\begin{proof}[Proof of Theorem \ref{thm: hardy gaussiana tipo B}]
If \(Q = p \theta\), the inequality is trivially satisfied. Thus, we assume \(Q \neq p \theta\). Choose the functions 
\[
V(r) = \frac{1}{r^{p(\theta - 1)}} \displaystyle{ e^{-\tfrac{ r^{ \alpha}}{ \beta}}}, 
\qquad \text{and} \qquad
W(r) = \left(\frac{1}{r^{p \theta}} - \frac{\alpha}{\beta} \left(\frac{p}{Q - p \theta} \right)\frac{1}{r^{p\theta - \alpha}}\right) \displaystyle{ e^{-\tfrac{ r^{ \alpha}}{ \beta}}}
\]
in Theorem~\ref{thm: d-radial_Hardy_poincare_inequality}. The solution to the corresponding differential equation~\eqref{eq: equazione_differenziale_per_varphi} is given by
\begin{equation*}
	\varphi(r) = r^{-\frac{Q - p \theta}{p}},
	\qquad	
	\lambda = \left|\frac{Q - p \theta}{p}\right|^p.
\end{equation*}
This establishes inequality~\eqref{eq: hardy gaussiana tipo B}.  
To prove sharpeness of the constant in~\eqref{eq: hardy gaussiana tipo B}, we introduce the family of test functions
\[
u_\varepsilon(x) = d^{-\frac{Q - p \theta}{p}} g_\varepsilon(d) \in C^\infty_c(\mathbb{R}^N \setminus \{0\}),
\]
where \(g_\varepsilon\) is defined in \eqref{eq: definizione g_veps}. We define 
\[
h(\varepsilon) := \int_0^{+\infty} \frac{1}{r} |g_\varepsilon(r)|^p \displaystyle{ e^{-\tfrac{ r^{ \alpha}}{ \beta}}} \, dr.
\]
It is straightforward to show that $h(\varepsilon)$ goes to infinity as $\varepsilon$ goes to $0^+$.

Using the coarea formula and applying Lemma~\ref{lemma: Lemma 2.3 in DA2} and Lemma~\ref{lemma: Lemma 2.1 in DA}, we obtain the following estimates
\begin{equation*}
\begin{split}
\int_{\mathbb{R}^N} \left| \nabla_\mathcal{L} u_\varepsilon \cdot \frac{\nabla_\mathcal{L} d}{|\nabla_\mathcal{L} d|} \right|^p \frac{1}{d^{p(\theta - 1)}} \displaystyle{ e^{-\tfrac{ d^{ \alpha}}{ \beta}}} \, dx 
&\leq \lambda_p\left\{ \left| \frac{Q - p \theta}{p} \right|^p h(\varepsilon) + \mathcal{O}(1)\right\},\\
\int_{\mathbb{R}^N} \frac{|u_\varepsilon|^p}{d^{p \theta}} |\nabla_\mathcal{L} d|^p \displaystyle{ e^{-\tfrac{ d^{ \alpha}}{ \beta}}} \, dx 
&= \lambda_p\, h(\varepsilon),
\\
\int_{\mathbb{R}^N} \frac{|u_\varepsilon|^p}{d^{p \theta - \alpha}} |\nabla_\mathcal{L} d|^p \displaystyle{ e^{-\tfrac{ d^{ \alpha}}{ \beta}}} \, dx 
&= \mathcal{O}(1).
\end{split}
\end{equation*}
Altogether we have
\begin{equation*}
	\left|\frac{Q - p \theta}{p}\right|^p 
	\leq 
\frac{ \displaystyle
 \int_{\mathbb{R}^N} \left| \nabla_\mathcal{L} u_\varepsilon \cdot \frac{\nabla_\mathcal{L} d}{|\nabla_\mathcal{L} d|} \right|^p \frac{1}{d^{p(\theta-1)}}  e^{-\tfrac{ d^{ \alpha}}{ \beta}} \, dx
}{
 \displaystyle\int_{\mathbb{R}^N} \frac{|u_\varepsilon|^p}{d^{p \theta}} |\nabla_\mathcal{L} d|^p \displaystyle{ e^{-\tfrac{ d^{ \alpha}}{ \beta}}} \, dx 
- \frac{\alpha}{\beta} \left(\frac{Q - p \theta}{p}\right) \int_{\mathbb{R}^N} \frac{|u_\varepsilon|^p}{d^{p \theta - \alpha}} |\nabla_\mathcal{L} d|^p  e^{-\tfrac{ d^{ \alpha}}{ \beta}} \, dx}
\leq 
\frac{ \lambda_p\left\{ \left| \frac{Q - p \theta}{p} \right|^p h(\varepsilon) + \mathcal{O}(1)\right\}}{ \lambda_p h(\varepsilon) + \mathcal{O}(1)}.
\end{equation*} 
Then the claim follows by taking the limit as $\varepsilon$ goes to $0^+.$ 
\end{proof}

\subsection{Hardy inequalities in an annulus}
\label{sezione: Hardy inequality in an Annulus}
In the previous sections (refer, for example, to Theorem~\ref{thm: sharp_d_radial_weighted_hardy_inequaity} in Section~\ref{sezione: Sharp Hardy inequalities related to fundamental solution}), we showed that various Hardy-type inequalities achieve the sharp constant when the underlying domain contains a neighborhood of the origin. It is therefore natural to investigate how the optimal constant changes (improves) when one assumes that there exists \(s > 0\) such that \(B^N_s \subseteq \Omega^c\), where, as usual, \(B^N_s\) denotes the Euclidean ball of radius \(s\) centered at the origin in \(\mathbb{R}^N\).

As a prototypical example, we examine the case where the domain \(\Omega\) is an annular region defined as
\[
\Omega_{a,b} = \{x \in \mathbb{R}^N : a < d(x) < b\},
\]
with \(0<a<b.\) This kind of results have been investigated in the case \(p = 2\) and in the Euclidean setting in the recent work by Gesztesy and Pang~\cite{GP}.

\medskip
We adopt here the same assumptions as in Section \ref{sezione: Sharp Hardy inequalities related to fundamental solution}. 

For fixed real parameters $\lambda$ and $\theta$, we consider the following boundary value problem
\begin{equation}
\label{eq: problema autovalori A}
\begin{cases} 
\left(r^{Q-1 - p(\theta-1)} |\varphi'|^{p-2} \varphi'\right)' + \lambda r^{Q-1 - p\theta} |\varphi|^{p-2} \varphi = 0, & \text{for } a < r < b, \\ 
\varphi(a) = \varphi(b) = 0.
\end{cases}
\end{equation}
We say that a real-valued function \(\varphi\) is a solution to~\eqref{eq: problema autovalori A} if \(\varphi \in C^1([a,b])\), \(r^{Q-1-p(\theta-1)} |\varphi'|^{p-2} \varphi' \in C^1([a,b])\), and \(\varphi\) satisfies the differential equation~\eqref{eq: problema autovalori A} pointwise for every \(r \in (a,b)\). 

Moreover, if for a given \(\lambda \in \mathbb{R}\) there exists a non-trivial solution \(\varphi\) of~\eqref{eq: problema autovalori A}, then \(\lambda\) is referred to as an eigenvalue, and \(\varphi\) is called the eigenfunction associated with \(\lambda\).

It is well known (see, for instance, Kusano and Naito~\cite{KN}) that the eigenvalue problem \eqref{eq: problema autovalori A} admits a countable sequence of positive, simple eigenvalues
\[
0 < \lambda_1 < \lambda_2 < \ldots < \lambda_n \to +\infty \quad \text{as } n \to +\infty,
\]
with no other eigenvalues. 
Notice that the positivity of the eigenvalues can be easily deduced by multiplying the differential equation by \(\varphi\) and integrating by parts.
Furthermore, the corresponding eigenfunctions \(\varphi_n\) have exactly \(n-1\) zeros within the interval \((a,b)\). In particular, $\varphi_1$ does not have zeroes in $(a,b).$

The next crucial lemma provides a lower bound for the size of the first eigenvalue.

\begin{lemma}
\label{lemma: lambda_1 e piu grande costante classica}
Let \(\lambda_1\) be the first eigenvalue of~\eqref{eq: problema autovalori A}. Then, the following lower bound holds
\[
\lambda_1 > \left|\frac{Q - p \theta}{p}\right|^p.
\]
\end{lemma}

\begin{proof}
We proceed by contradiction and assume that \(\lambda_1 \leq \left|\frac{Q - p\theta}{p}\right|^p\). Consider the following two problems
\begin{equation}
\label{eq: problema autovalori A per lambda_1}
\begin{cases}  
\left(r^{Q-1 - p(\theta-1)} |\varphi_1'|^{p-2} \varphi_1'\right)' + \lambda_1 r^{Q-1 - p\theta} |\varphi_1|^{p-2} \varphi_1 = 0, & \text{for } a < r < b, \\  
\varphi_1(a) = \varphi_1(b) = 0,
\end{cases}
\end{equation}
and
\begin{equation}
\label{eq: problema autovalori A su tutto R}
\left(r^{Q-1 - p(\theta-1)} |\phi'|^{p-2} \phi'\right)' + \left|\frac{Q - p\theta}{p}\right|^p r^{Q-1 - p\theta} |\phi|^{p-2} \phi = 0, \quad \text{for } r \in (a,b).
\end{equation}
Without loss of generality, we can assume that \(\varphi_1 > 0\) in \((a,b)\).

Moreover, a non-vanishing solution of~\eqref{eq: problema autovalori A su tutto R} is given by \(\phi = r^{-\frac{Q - p\theta}{p}} > 0\). 

Using identity~\eqref{eq: identita_fondamentale} with the choices 
\[
\Omega = (a,b), 
\qquad f = r^{\frac{Q-1}{p} - (\theta-1)} \varphi_1', 
\qquad g = r^{\frac{Q-1}{p} - (\theta-1)} \varphi_1 \frac{\phi'}{\phi},
\]
and performing an integration by parts in  the mixed product \((|g|^{p-2} g, f)\), we obtain
\[
||w(p,f,g)(f-g)||^2_{L^2(a,b)} = \int_a^b r^{Q-1-p(\theta-1)} |\varphi_1'|^p \, dr - \left|\frac{Q-p\theta}{p}\right|^p \int_a^b r^{Q-1-p\theta} |\varphi_1|^p \, dr.
\]
Using the identity 
\[
\int_a^b r^{Q-1-p(\theta-1)} |\varphi_1'|^p \, dr = \lambda_1 \int_a^b r^{Q-1-p\theta} |\varphi_1|^p \, dr,
\]
obtained from~\eqref{eq: problema autovalori A per lambda_1} multiplying by $\varphi_1,$ integrating over $(a,b)$ and then integrating by parts, we find 
\[
0 \leq ||w(p,f,g)(f-g)||^2_{L^2(a,b)} = \left(\lambda_1 - \left|\frac{Q-p\theta}{p}\right|^p\right)  \int_a^b r^{Q-1-p\theta} |\varphi_1|^p \, dr \leq 0.
\]
From this, it follows that the \(L^2\)-norm must be zero, i.e., \(f = g\). This implies that \(\varphi_1 = c\phi\) for some constant \(c \neq 0\). However, this leads to a contradiction, as \(\phi\) does not satisfy the boundary conditions imposed on \(\varphi_1\). This concludes the proof.
\end{proof}

Now we are in position to prove the main result of this section. 

\begin{theorem}
\label{thm: sharp hardy inequality in an annulus}
Let \(\theta \in \mathbb{R},\) then for all \(u \in C^\infty_c(\Omega_{a,b})\) one has
\begin{equation}
\label{eq: sharp hardy inequality in an annulus}
\int_{\Omega_{a,b}} \left|\nabla_\mathcal{L} u \cdot \frac{\nabla_\mathcal{L} d}{|\nabla_\mathcal{L} d|} \right|^p \frac{1}{d^{p(\theta-1)}} \, dx \geq \lambda_1 \int_{\Omega_{a,b}} \frac{|u|^p}{d^{p\theta}} |\nabla_\mathcal{L} d|^p \, dx,
\end{equation}
where \(\lambda_1\) denotes the first eigenvalue of~\eqref{eq: problema autovalori A}. Moreover, this constant is sharp, and the candidate maximizing function is the corresponding eigenfunction \(\varphi_1\).
\end{theorem}

\begin{remark}
	This result highlights the close connection between the best constants of Hardy-Poincaré-type inequalities and the first eigenvalue in a related spectral problem. This relation has been intensively explored with a particular emphasis on the role played by the underline domain geometry (see, for instance,~\cite{Krejcirik08,KZ10,EKK08,BK25}).
\end{remark}

\begin{remark}
Being \(\lambda_1\) strictly greater than the sharp constant of a classical Hardy inequality on a domain that contains a neighborhood of the origin, namely \(\left|\frac{Q - p\theta}{p}\right|^p\) (Lemma~\ref{lemma: lambda_1 e piu grande costante classica}), inequality~\eqref{eq: sharp hardy inequality in an annulus} shows that in a the annulus $\Omega_{a,b}$ we indeed have an improvement of the best constant. 
As a matter of fact, the inequality~\eqref{eq: sharp hardy inequality in an annulus} is never critical, meaning that \(\lambda_1 \neq 0\), even in the case \(Q = p\theta\).
\end{remark}

\begin{proof}
The validity of the inequality~\eqref{eq: sharp hardy inequality in an annulus} follows directly from Theorem~\ref{thm: d-radial_Hardy_poincare_inequality}, which also shows that \(\varphi_1\) is the candidate maximizing function. 
Since \(\varphi_1(d) \in C^1(\Omega_{a,b})\) and \(\varphi_1(d(x)) = 0\) on \(\partial \Omega_{a,b}\), a standard approximation with smooth, compactly supported functions in \(\Omega_{a,b}\) shows that the constant \(\lambda_1\) is sharp.
\end{proof}

\begin{remark}
In the special case \(p = 2\), the constant \(\lambda_1\) and the function \(\varphi_1\) can be computed explicitly. The differential equation in \eqref{eq: problema autovalori A} simplifies to the second-order Cauchy-Euler equation
\[
r^2 \varphi''(r) + (Q - 2\theta + 1) r \varphi'(r) + \lambda_1 \varphi(r) = 0, \quad \text{for } r \in (a,b),
\]
whose general solution is given by
\[
\varphi(r) = b_1 r^{-\frac{Q - 2\theta}{2}} \sin\left(C \ln \frac{r}{a}\right) + b_2 r^{-\frac{Q - 2\theta}{2}} \cos\left(C \ln \frac{r}{a}\right),
\qquad C = \sqrt{\lambda_1 - \left(\frac{Q - 2\theta}{2}\right)^2},
\quad b_1, b_2 \in \mathbb{R}. 
\]
The boundary conditions uniquely determine the specific solution, namely
\begin{equation*}
\varphi_1(r) = r^{-\frac{Q - 2\theta}{2}} \sin\left(\frac{\pi}{\ln \frac{b}{a}} \ln \frac{r}{a}\right),
\qquad 
\lambda_1 = \left(\frac{Q - 2\theta}{2}\right)^2 + \left(\frac{\pi}{\ln \frac{b}{a}}\right)^2.
\end{equation*}
Thus, one has the following result.

\begin{corollary}
\label{cor:sharp hardy inequality in an annulus}
Let  $d=2.$ Let \(\theta \in \mathbb{R},\) then for all \(u \in C^\infty_c(\Omega_{a,b})\) one has
\begin{equation}
\int_{\Omega_{a,b}} \left|\nabla_\mathcal{L} u \cdot \frac{\nabla_\mathcal{L} d}{|\nabla_\mathcal{L} d|} \right|^2 \frac{1}{d^{2(\theta-1)}} \, dx \geq \left[\left(\frac{Q - 2\theta}{2}\right)^2 + \left(\frac{\pi}{\ln \frac{b}{a}}\right)^2\right] \int_{\Omega_{a,b}} \frac{|u|^2}{d^{2\theta}} |\nabla_\mathcal{L} d|^2 \, dx.
\end{equation}
Moreover, the constant is sharp, and the candidate maximizing function is 
\[
d(x)^{-\frac{Q - 2\theta}{2}} \sin\left(\frac{\pi}{\ln \frac{b}{a}} \ln \frac{d(x)}{a}\right).
\]
\end{corollary}
\end{remark}

\section{Application II: Hardy inequality for special domains and antisymmetric functions}
\label{section:hardy-special}
Assume \(N \geq 2\) and fix $h\leq N,$ consider the splitting \(\mathbb{R}^N = \mathbb{R}^m \times \mathbb{R}^{N-m}\), with \(2 \leq m \leq h,\) with coordinates $(x,y)$ such that \(x \in \mathbb{R}^m\) and \(y \in \mathbb{R}^{N-m}\).
Assume that the matrix \(\sigma \in \mathcal{M}_{h \times N}\), which defines the $h$ components-vector \(\nabla_\mathcal{L}\) in~\eqref{eq: nabla_L come gradiente euclideo}, has the following block structure
\begin{equation}
\sigma = 
\begin{pmatrix} 
\mathcal{I}_m & \sigma_1 \\ 
0 & \sigma_2 
\end{pmatrix},
\end{equation}
where \(\sigma_1 \in \mathcal{M}_{m \times (N-m)}\), \(\sigma_2 \in \mathcal{M}_{(h-m) \times (N-m)}\), and \(\mathcal{I}_m\) denotes the \(m \times m\) identity matrix. 
Moreover, we impose the condition
\begin{equation}
\label{eq: condizione su sigma_1 per vandermonde}
\frac{\partial}{\partial y_j} (\sigma_1)_{ij} = 0, \quad \forall i = 1, \dots, m, \quad \forall j = 1, \dots, N-m.
\end{equation}
In this section we investigate a different setting compared with the previous sections, namely we consider the case in which the origin is no longer inside \(\Omega\) but instead lies on its boundary. More precisely, we consider a domain of the form \(
\Omega = \Omega_m \times \mathbb{R}^{N-m} \subseteq \mathbb{R}^N,
\)
where  
\[
\Omega_m = \{(x_1, \dots, x_m) \in \mathbb{R}^m \mid x_1 < x_2 < \dots < x_m\}.
\]  
Let \(\nu_m: \mathbb{R}^N \to \mathbb{R}\) be defined as  
\[
\nu_m(x_1, \dots, x_m, y_1, \dots, y_{N-m}) = \prod_{i < j \leq m} (x_j - x_i).
\]
It is immediate to observe that \(\nu_m > 0\) in \(\Omega\) and \(\nu_m = 0\) on \(\partial \Omega\). In the special case \(m = h = N\), the function \(\nu_N\) coincides with the Vandermonde determinant. Moreover, $\nu_N$ is antisymmetric, meaning that
\begin{equation}
\label{eq: condizione di antisimmetria}
\nu_N(x_1, \dots, x_i, \dots, x_j, \dots, x_N) = -\nu_N(x_1, \dots, x_j, \dots, x_i, \dots, x_N), \quad \forall\, 1 \leq i \neq j \leq N.
\end{equation}
We also observe that when \(m = N\), the setting reduces to the classical Euclidean framework. Indeed, in this case, \(\nabla_\mathcal{L}\) coincides with the standard Euclidean gradient \(\nabla\), and \(\mathcal{L}\) reduces to the Laplacian \(\Delta\).  
The gradient of \(\nu_m\) is given by  
\(
\nabla_\mathcal{L} \nu_m = (\nabla_x \nu_m, 0_{h-m}),
\)
where \(\nabla_x\) denotes the Euclidean gradient in \(\mathbb{R}^m\), and \(0_{h-m}\) is the zero vector in \(\mathbb{R}^{h-m}\). Moreover, \(\nu_m\) satisfies    
\[
\nabla_\mathcal{L} \nu_m \cdot \left(\frac{x}{|x|}, 0_{h-m} \right) = \nabla_x \nu_m \cdot \frac{x}{|x|} = \frac{m(m-1)}{2} \frac{\nu_m}{|x|}.
\]
This follows from the fact that the restriction \(\nu_m|_{\mathbb{R}^m}\) is homogeneous of degree \(\frac{m(m-1)}{2}\) with respect to Euclidean dilations in \(\mathbb{R}^m\). Finally,
\[
\mathcal{L} \nu_m = \Delta_x \nu_m = 0,
\]
where \(\Delta_x\) is the standard Euclidean Laplacian in \(\mathbb{R}^m\). This follows directly from assumption~\eqref{eq: condizione su sigma_1 per vandermonde}.  

We can now state the main result of this section.
\begin{theorem}
\label{thm: hardy in special domain}
Let \(\theta \in \mathbb{R}\), then for all \(u \in C^\infty_c(\Omega)\) one has
\begin{equation}
\label{eq: hardy in special domain}
\int_\Omega \frac{|\nabla_\mathcal{L} u(x,y)|^2}{|x|^{2(\theta-1)}} \, dx \, dy \geq \left( \left(\frac{m^2 - 2\theta}{2}\right)^2 +  m (m-1) (\theta -1) \right) \int_\Omega \frac{|u(x,y)|^2}{|x|^{2\theta}} \, dx \, dy.
\end{equation}
Moreover, if \(m = N\) and \(N^2> 2\theta\), the inequality extends to the larger space \(C^\infty_0(\Omega)\) of functions with compact support in \(\overline{\Omega}\) and, in this case, the constant is sharp, with the candidate maximizing function given by  
\(
|x|^{-\frac{N^2 - 2\theta}{2}} \nu_N(x).
\)
\end{theorem}

\begin{remark}
Setting \(\theta = 1\) in inequality~\eqref{eq: hardy in special domain} yields the constant \(\big(\frac{m^2 - 2}{2}\big)^2\), which closely resembles the sharp constant of the cylindrical Hardy inequality~\eqref{eq: sharp_cilindrical_hardy_inequality}, except that it now involves \(m^2\) instead of \(m\). In particular, here the constant remains always strictly positive.  
\end{remark}

\begin{proof} 
We want to show~\eqref{eq: hardy in special domain} as a consequence of~\eqref{eq: hardy_completa} in Theorem~\ref{thm: Hardy_completa}. For the hypotheses of Theorem~\ref{thm: Hardy_completa} to be satisfied, it is sufficient to exhibit a nontrivial function \( h \in C^1(\Omega) \), with \( h \neq 0 \), such that $h$ is a solution of the equation
\begin{equation*}
	-\Div_\mathcal{L} \left(\frac{1}{|x|^{2(\theta-1)}} \nabla_\mathcal{L} h \right) = \lambda \frac{h}{|x|^{2\theta}}, 
	\qquad  
\lambda = \left(\frac{m^2 - 2\theta}{2}\right)^2 + m(m-1)(\theta -1).  
\end{equation*}  
Seeking a solution of the form \( h(x,y) = |x|^\alpha \nu_m(x,y) \) and substituting into the differential equation, we obtain  
\[
\lambda = -\alpha^2 - \alpha (m^2 - 2\theta) + m(m-1)(\theta -1),
\]
which gives  
\(
\alpha = -\frac{m^2 - 2\theta}{2}.
\)

In the case \(m = N\), we show that the inequality extends to functions in \(C^\infty_0(\Omega)\) with compact support in \(\overline{\Omega}\). 
In order to do that we introduce the following approximation pair 
\begin{equation*}
	V_\varepsilon(x) = \frac{1}{(|x|^2 + \varepsilon^2)^{\theta - 1}}, 
	\qquad 
	W_\varepsilon(x) := -\frac{1}{h_\varepsilon(x)} \Div \left(V_\varepsilon(x) \nabla h_\varepsilon \right),
\end{equation*}
where 
\begin{equation*}
	h_\varepsilon(x) = (|x|^2 + \varepsilon^2)^{-\frac{N^2 - 2\theta}{4}} \left(\nu_N^2(x) + \varepsilon^2\right)^{\frac{1}{2}}.
\end{equation*}

Reproducing the same argument as in Theorem~\ref{thm: Hardy_completa} for \(p=2\), \emph{i.e.} taking 
\begin{equation*}
	f_\varepsilon(x)= V_\varepsilon(x)^{1/2} \nabla h_\varepsilon,
	\qquad
	g_\varepsilon(x)=\frac{\nabla h_\varepsilon}{h_\varepsilon} V_\varepsilon(x)^{1/2},
\end{equation*}
and observing that 
\[
\left. \frac{|u|^2}{h_\varepsilon} V_\varepsilon \nabla h_\varepsilon \right|_{\partial \Omega} = 0, 
\qquad \forall\, u \in C^\infty_0(\Omega) \text{ with compact support in } \overline{\Omega},
\]
we obtain  
\[
\int_\Omega \frac{|\nabla u|^2}{ ( |x|^2+\varepsilon^2)^{\theta-1}} dx \geq \int_\Omega |u|^2 W_\varepsilon(x) dx, 
\qquad \forall\, u \in C^\infty_0(\Omega) \text{ with compact support in } \overline{\Omega}.
\]
Computing explicitly \(W_\varepsilon(x)\) one gets  
\begin{multline*}
W_\varepsilon(x) = -\left(\frac{N^2 - 2\theta}{2} \right) \left(\frac{N^2 + 2\theta}{2} \right) \frac{|x|^2}{(|x|^2+\varepsilon^2)^{\theta+1}}  
+ N \left(\frac{N^2 - 2\theta}{2}\right) \frac{1}{(|x|^2+\varepsilon^2)^\theta} \\
+ N(N-1) \left(\frac{N^2-2}{2}\right) \frac{1}{(|x|^2+\varepsilon^2)^\theta} \frac{\nu_N^2}{\nu_N^2+\varepsilon^2}  
- \frac{\varepsilon^2}{(|x|^2+\varepsilon^2)^{\theta-1}} \frac{|\nabla \nu_N|^2}{(\nu_N^2+\varepsilon^2)^2}.
\end{multline*}
Taking the limit as $\varepsilon$ goes to $0^+,$ we obtain
\[
\int_\Omega \frac{|\nabla u(x)|^2}{|x|^{2(\theta-1)}} \, dx  \geq \left( \left(\frac{N^2 - 2\theta}{2}\right)^2 + (\theta -1) N (N-1) \right) \int_\Omega \frac{|u(x)|^2}{|x|^{2\theta}} \, dx,
\]
for all \( u \in C^\infty_0(\Omega) \) with compact support in \(\overline{\Omega}\). 
This follows from Fatou's lemma and Beppo-Levi’s theorem if \(\theta \geq 1\), or otherwise from uniform convergence.  
Observe that condition \(N^2 > 2\theta\) ensures the integrability of the right-hand side, making it possible taking the limit. Indeed, since \(u\) vanishes on \(\partial \Omega\), the function \(\varphi = u/\nu_N\) is smooth up to the boundary, so one has 
\[
\int_\Omega \frac{|u|^2}{|x|^{2\theta}} dx
\lesssim
\int_K \frac{|\nu_N(x)|^2}{|x|^{2\theta}} dx = \int_K |x|^{N(N-1)- 2\theta} \left|\nu_N\left(\frac{x}{|x|}\right) \right|^2 dx < +\infty\quad\text{for some compact set \( K \subseteq \overline{\Omega} \).}
\]
To prove that the constant is sharp, we employ a standard approximation argument. Let  
\[
u_\varepsilon(x) = |x|^{-\frac{N^2-2\theta}{2}} \nu_N(x) g_\varepsilon(|x|) = |x|^{-\frac{N-2\theta}{2}} \nu_N\left(\frac{x}{|x|}\right) g_\varepsilon(|x|),
\]
where, as usual, \( g_\varepsilon(r) \) is defined in \eqref{eq: definizione g_veps}. For every \(\varepsilon > 0\), \( u_\varepsilon \) is an admissible function with compact support contained in \(\overline{\Omega}\).  
Recall that in polar coordinates \((r, \theta)\), where \(r = |x|\) and \(\theta = \frac{x}{r}\), the Laplace operator is given by  
\[
\Delta = \partial^2_{rr} + \frac{N-1}{r} \partial_r + \frac{1}{r^2} \Delta_{\mathbb{S}^{N-1}},
\]
where \(\partial_r\) and \(\partial^2_{rr}\) denote the first and second order derivatives with respect to the radial variable \(r\), while \(\Delta_{\mathbb{S}^{N-1}}\) is the Laplace-Beltrami operator associated with the standard metric on \(\mathbb{S}^{N-1}\), the \((N-1)\)-dimensional unit sphere in \(\mathbb{R}^N\).  
Moreover, if \(P_k\) is a homogeneous harmonic polynomial of degree \(k\), then the corresponding spherical harmonic  
\(
\phi_k(x) = P_k\left(\frac{x}{|x|}\right)
\)
is an eigenfunction of \(-\Delta_{\mathbb{S}^{N-1}}\) with eigenvalue \(k(k+N-2)\), \emph{i.e.} 
\[
-\Delta_{\mathbb{S}^{N-1}} \phi_k = k(k+N-2) \phi_k.
\]
In particular, the function \(\nu_S := \nu_N\left(\frac{x}{|x|}\right)\) satisfies  
\begin{equation}
\label{eq:Laplace-Beltrami-nu}
-\Delta_{\mathbb{S}^{N-1}} \nu_S = \frac{N(N-1)}{2} \left(\frac{N^2 + N - 4}{2}\right) \nu_S.
\end{equation}
Since \( u_\varepsilon \) is antysymmetric (see~\eqref{eq: condizione di antisimmetria}), then one has
\[
\int_\Omega \frac{|\nabla u_\varepsilon |^2 }{|x|^{2(\theta-1)}} dx = \frac{1}{N!} \int_{\mathbb{R}^N} \frac{|\nabla u_\varepsilon|^2}{|x|^{2(\theta-1)}}dx.
\]
To further simplify the left-hand side above, one employs~\eqref{eq:Laplace-Beltrami-nu}  together with the standard decomposition  
\[
|\nabla u|^2 = |\partial_r u|^2 + \frac{1}{r^2} |\nabla_{\mathbb{S}^{N-1}} u |^2,
\]
to obtain  
\[
\int_\Omega \frac{|\nabla u_\varepsilon |^2 }{|x|^{2(\theta-1)}} dx =
\frac{1}{N!} \int_0^{+\infty} \int_{\mathbb{S}^{N-1}} \left( \frac{|\partial_r (r^{-\frac{N-2\theta}{2}} g_\varepsilon(r))|^2}{r^{2(\theta-1)}} |\nu_S|^2 + \frac{1}{r^N} |g_\varepsilon(r)|^2 |\nabla_{\mathbb{S}^{N-1}} \nu_S|^2 \right) r^{N-1} d\mathcal{H}^{N-1} dr.
\]
Expanding the radial derivative inside the square and using~\eqref{eq:g_eps1}-\eqref{eq:g_eps3} in Proposition~\ref{prop:straighforward-g_eps} we get
\[
\int_\Omega \frac{|\nabla u_\varepsilon |^2 }{|x|^{2(\theta-1)}} dx 
\leq \frac{1}{N!} \left( \left(\frac{N^2 - 2\theta}{2}\right)^2 + N(N-1)(\theta-1) \right) \int_{\mathbb{S}^{N-1}} |\nu_S|^2 d\mathcal{H}^{N-1} \left\{ - \ln(4\varepsilon^2) + \mathcal{O}(1) \right\}.
\]
Similarly,   
\[
\int_\Omega \frac{|u_\varepsilon|^2}{|x|^{2\theta}} dx = \frac{1}{N!} \int_{\mathbb{S}^{N-1}} |\nu_S|^2 d\mathcal{H}^{N-1} \left\{ -\ln(4\varepsilon^2) + \mathcal{O}(1) \right\}.
\]
This concludes the proof of the result.
\end{proof}

As a corollary of the previous result Theorem~\ref{thm: hardy in special domain}, when \(m = N\) (\emph{i.e.} in the Euclidean setting), one gets Hardy-type inequalities for antisymmetric functions (see~\eqref{eq: condizione di antisimmetria}). We introduce the following notation for the space of smooth, compactly supported, antisymmetric functions
\begin{equation*}
	C^\infty_{c,A}(\mathbb{R}^N):=
\{ u \in C^\infty_c(\mathbb{R}^N) \colon u \text{ antisymmetric} \}.
\end{equation*}

\begin{corollary}
\label{corollario: hardy funzioni antisimmetriche}
Let \(\theta \in \mathbb{R}\) and assume \(N^2> 2\theta\), then for all $u\in C^\infty_{c,A}(\mathbb{R}^N)$ one has 
\begin{equation}
\label{eq: hardy funzioni antisimmetriche}
\int_{\mathbb{R}^N} \frac{|\nabla u|^2}{|x|^{2(\theta-1)}} dx \geq \left( \left(\frac{N^2 - 2\theta}{2} \right)^2 + N(N-1)(\theta-1) \right) \int_{\mathbb{R}^N} \frac{|u|^2}{|x|^{2\theta}} dx, \quad \forall u \in C^\infty_{c,A}(\mathbb{R}^N).
\end{equation}
The candidate maximizing function is given by  
\(
|x|^{-\frac{N^2 - 2\theta}{2}} \nu_N(x).
\) 
Moreover, the constant is sharp. 
\end{corollary}

\begin{remark}
The inequality~\eqref{eq: hardy funzioni antisimmetriche} was first established in the case \(\theta = 1\) by Hoffmann-Ostenhof, Laptev and Shcherbakov in~ \cite{HOL}, employing the decomposition into spherical harmonics. For an \(L^p\)-version of the inequality we refer to Kijaczko~\cite{K}.  
\end{remark}

\begin{proof}
It suffices to observe that, for any $\Omega$ of the form  
\[
\Omega = \{(x_1, \dots, x_N) \mid x_1 < x_2 < \dots < x_N\},
\]
then every antisymmetric function vanishes on \(\partial \Omega\). Moreover, one easily has
\[
\int_{\mathbb{R}^N} \frac{|\nabla u|^2}{|x|^{2(\theta-1)}} dx = N! \int_{\Omega} \frac{|\nabla u|^2}{|x|^{2(\theta-1)}} dx, 
\qquad \text{and} \qquad
\int_{\mathbb{R}^N} \frac{|u|^2}{|x|^{2\theta}} dx = N! \int_{\Omega} \frac{|u|^2}{|x|^{2\theta}} dx,
\]
for all $u \in C^\infty_{c,A}(\mathbb{R}^N).$
Then the thesis is given directly from applying Theorem~\ref{thm: hardy in special domain}.  
\end{proof}

\section{Application III: Hardy-type inequalities for some subelliptic operators}\label{section:examples-subelliptic}
In this section, we apply the previous general results to specific settings.

\subsection{Euclidean Setting}
In \(\mathbb{R}^N\), we assume $\sigma(x)=\mathcal{I}_N$ in~\eqref{eq: nabla_L come gradiente euclideo}, namely we consider the standard $N$ components-vector field \(\nabla_\mathcal{L} = \nabla,\). The associated differential operators $\mathcal{L}$ and $\mathcal{L}_p$ being the standard Laplaciand $\Delta$ and the \(p\)-Laplacian \(\Delta_p\), respectively. 

Let \(\delta_\lambda\) denote the family of dilations defined by \(\delta_\lambda(x) = (\lambda x_1, \dots, \lambda x_N)\), and set \(d(x) = |x|\).  

Theorem~\ref{thm: sharp_d_radial_weighted_hardy_inequaity} applied in this situation, assuming \(Q = N\), allows to recover the classical weighted Hardy-type inequality in the Euclidean setting. More precisely, we have the following result:  
\begin{theorem}
Let \(\theta \in \mathbb{R}\), and let \(\Omega \subseteq \mathbb{R}^N \setminus \{0\}\), then for all $u \in C^\infty_c(\Omega)$ one has
\begin{equation*}
\int_\Omega \left| \nabla u \cdot \frac{x}{|x|} \right|^p \frac{1}{|x|^{p(\theta - 1)}} \, dx 
\geq \left| \frac{N - p \theta}{p} \right|^p \int_\Omega \frac{|u|^p}{|x|^{p \theta}}  \, dx.
\end{equation*}
Moreover, if \(\Omega \cup \{0\}\) contains a neighborhood of the origin, then the constant is sharp.
\end{theorem}  
Theorem~\ref{thm: sharp_logaritmic_hardy_inequality} gives the following logarithmic inequality.  
\begin{theorem}
Let \(\theta \in \mathbb{R}\), and let \(\Omega \subseteq \mathbb{R}^N \setminus \{0\}\), then for all \(u \in C^\infty_c(\Omega)\) one has
\begin{equation*}
\int_\Omega \left| \nabla u \cdot \frac{x}{|x|} \right|^p \left( \ln \frac{R}{|x|} \right)^{\theta + p} \frac{1}{|x|^{N-p}} \, dx 
\geq \left| \frac{\theta + 1}{p} \right|^p \int_\Omega \frac{|u|^p}{|x|^N} \left( \ln \frac{R}{|x|} \right)^\theta \, dx.
\end{equation*}
Moreover, if \(\Omega \cup \{0\}\) contains a neighborhood of the origin, then the constant is sharp.
\end{theorem}  
By considering weight functions of Gaussian-type, from Theorem~\ref{thm: hardy gaussiana tipo A} and Theorem~\ref{thm: hardy gaussiana tipo B}, we get the two results stated below.
\begin{theorem}
Let \(\alpha\) and \(\beta\) be constants such that \(\alpha \geq 2\) and \(\beta > 0,\) then for all \(u \in C^\infty_c(\mathbb{R}^N)\) one has
\begin{multline*}
\int_{\mathbb{R}^N} \left| \nabla u \cdot \frac{x}{|x|} \right|^p \displaystyle{e^{-\tfrac{|x|^{ \alpha}}{ \beta}}}\, dx 
\geq
\left( \frac{\alpha}{p \beta} \right)^p \int_{\mathbb{R}^N} |u|^p |x|^{p( \alpha - 1)} \displaystyle{e^{-\tfrac{|x|^{ \alpha}}{ \beta}}}\, dx \\
- \left( \frac{\alpha}{p \beta} \right)^{p-1} \left(\alpha  (p - 1) + N - p\right) \int_{\mathbb{R}^N} |u|^p |x|^{ \alpha (p - 1) - p} \displaystyle{e^{-\tfrac{|x|^{ \alpha}}{ \beta}}}\, dx.
\end{multline*}
The candidate maximizing function is \(\displaystyle{e^{\frac{|x|^{ \alpha}}{p \beta}}}\).
Moreover, this inequality is sharp in the sense that it cannot hold if the right-hand side is multiplied by any constant \(c > 1.\)
\end{theorem}
\begin{theorem}
Let \(\theta \in \mathbb{R}\) and  let \(\alpha, \beta\) be two constants such that \(\alpha \geq 2\) and \(\beta > 0,\) then for all \(u \in C^\infty_c(\mathbb{R}^N \setminus \{0\})\) one has
\begin{multline*}
\int_{\mathbb{R}^N} \left| \nabla u \cdot \frac{x}{|x|} \right|^p \frac{1}{|x|^{p(\theta-1)}} \displaystyle{e^{-\tfrac{|x|^{ \alpha}}{ \beta}}} \, dx \geq
\left| \frac{N - p \theta}{p} \right|^p \int_{\mathbb{R}^N} \frac{|u|^p}{|x|^{p \theta}} \displaystyle{e^{-\tfrac{|x|^{ \alpha}}{ \beta}}} \, dx\\ 
- \frac{\alpha}{\beta} \left(\frac{N - p \theta}{p}\right) \left| \frac{N - p \theta}{p} \right|^{p-2} \int_{\mathbb{R}^N} \frac{|u|^p}{|x|^{p \theta - \alpha}} \displaystyle{e^{-\tfrac{|x|^{ \alpha}}{ \beta}}} \, dx.
\end{multline*}
The candidate maximizing function is \(|x|^{-\frac{N - p \theta}{p}}\).
Moreover, if \(N \neq p \theta\), the inequality is sharp in the sense that it does not hold if the right-hand side is multiplied by any constant \(c > 1.\) 
\end{theorem}  
In the setting of annular domains, Theorem~\ref{thm: sharp hardy inequality in an annulus} and Corollary~\ref{cor:sharp hardy inequality in an annulus} give
\begin{theorem}
Let \(\theta \in \mathbb{R},\) then for any \(u \in C^\infty_c(\Omega_{a,b})\),one has
\begin{equation}
\int_{\Omega_{a,b}} \left|\nabla u \cdot \frac{x}{|x|} \right|^p \frac{1}{|x|^{p(\theta-1)}} \, dx \geq \lambda_1 \int_{\Omega_{a,b}} \frac{|u|^p}{|x|^{p\theta}} \, dx,
\end{equation}
where \(\lambda_1\) denotes the first eigenvalue of~\eqref{eq: problema autovalori A}.
The candidate maximizing function is the corresponding eigenfunction \(\varphi_1\).
Moreover, the constant is sharp. 
\end{theorem}
\begin{remark}\label{rmk:first-eigen-prop}
Remember (ref. Lemma~\ref{lemma: lambda_1 e piu grande costante classica}) that \(\lambda_1 > \left| \frac{Q - p\theta}{p} \right|^p\): Moreover, for \( p = 2, \) $\lambda_1$ is explicit (ref. Corollary~\ref{cor:sharp hardy inequality in an annulus}) and it is given by $
\lambda_1 = \left( \frac{Q - 2\theta}{2} \right)^2 + \left( \frac{\pi}{\ln \frac{b}{a}} \right)^2.$
\end{remark}

\smallskip
Observe that one can rewrite $\sigma(x)=\mathcal{I}_N$ as 
\begin{equation*}
	\sigma 
	= \begin{pmatrix} \mathcal{I}_m & 0 \\ 0 & \mathcal{I}_{N-m} \end{pmatrix},
	\qquad \text{for} \quad 1 \leq m \leq N.
\end{equation*}
Thus, $\sigma$ has the form~\eqref{eq: matrice sigma cilindrica}, with $\sigma_1=0_{m,N-m}$ and $\sigma_2=\mathcal{I}_{N-m,N-m}.$

Since \(\sigma_1 \equiv 0\), hypothesis~\eqref{eq: condizione su sigma_1 per vandermonde0} is trivially satisfied. We shall adopt the notation \(\mathbb{R}^N = \mathbb{R}^m \times \mathbb{R}^{N-m} \ni (x,y)\), where \(x \in \mathbb{R}^m\) and \(y \in \mathbb{R}^{N-m}\). 

Since hypotheses of Theorem~\ref{thm: sharp_cilindrical_hardy_inequality} and Theorem~\ref{thm: Cylindrical logaritmic hardy inequality} are satisfied one has the following results. 
\begin{theorem}
Fix \(\theta \in \mathbb{R}\), and assume \(\Omega \subseteq \mathbb{R}^m \setminus \{0\} \times \mathbb{R}^{N-m}\), then for all \(u \in C^\infty_c(\Omega)\) one has
\begin{equation*}
\int_\Omega \left| \nabla u \cdot \left( \frac{x}{|x|}, 0_{N-m} \right) \right|^p \frac{1}{|x|^{p(\theta - 1)}} \, dx \, dy 
\geq \left| \frac{m - p \theta}{p} \right|^p \int_\Omega \frac{|u(x, y)|^p}{|x|^{p \theta}} \, dx \, dy.
\end{equation*}
Moreover, if \(B^m_r \setminus \{0\} \times B^{N-m}_r \subseteq \Omega\) for some \(r > 0\), then the constant is sharp.
\end{theorem}
\begin{theorem}
Let \(R > 0\), fix \(\,\theta \in \mathbb{R}\), and assume \(\Omega = B^m_R \setminus \{0\} \times \mathbb{R}^{N-m},\) then for all \(u \in C^\infty_c(\Omega)\) one has
\begin{equation*}
\int_\Omega \left| \nabla u \cdot \left(\frac{x}{|x|}, 0_{N-m}\right) \right|^p \frac{1}{|x|^{m-p}} \left(\ln\left(\frac{R}{|x|}\right)\right)^{\theta+p} \, dx \, dy 
\geq \left|\frac{\theta+1}{p}\right|^p \int_\Omega \frac{|u|^p}{|x|^m} \left(\ln\left(\frac{R}{|x|}\right)\right)^\theta \, dx \, dy.
\end{equation*}
Moreover, the constant is sharp.
\end{theorem}

\subsection{Heisenberg Greiner Operator}
Let \(\gamma \geq 1\). In \(\mathbb{R}^{2n+1}\), with coordinates \((z,t) = (x,y,t) \in \mathbb{R}^n \times \mathbb{R}^n \times \mathbb{R}\), we consider the vector fields  
\[
X_i = \frac{\partial}{\partial x_i} + 2 \gamma y_i |z|^{2\gamma-2} \frac{\partial}{\partial t},  
\qquad  
Y_i = \frac{\partial}{\partial y_i} - 2 \gamma x_i |z|^{2\gamma-2} \frac{\partial}{\partial t},  
\qquad i = 1, \dots, n.
\]
The associated horizontal gradient is given by \(\nabla_\mathcal{L} = \sigma \nabla\), where  
\[
\sigma (x,y,t) =  
\begin{pmatrix}  
\mathcal{I}_n & 0 & 2 \gamma y |z|^{2\gamma-2}\\  
0 & \mathcal{I}_n & -2 \gamma x |z|^{2\gamma-2}  
\end{pmatrix}
\in \mathcal{M}_{2n\times (2n+1)}.
\]
\begin{remark}
When \(p=2\) and \(\gamma=1\), the operator \(\mathcal{L}_p\) coincides with the standard sub-Laplacian on the Heisenberg group  
\[
\mathcal{L}_2 = \Delta_{\mathbb{H}} = \Delta_z + 4|z|^2 \frac{\partial^2}{\partial t^2} + 4 \frac{\partial}{\partial t} (T \cdot),
\qquad
\text{where}  
\qquad
T= \sum_{i=1}^n \left( y_i \frac{\partial }{\partial x_i} - x_i \frac{\partial }{\partial y_i} \right).
\]
For \(p=2\) and \(\gamma > 1\), the operator \(\mathcal{L}_p\) is referred to as the Greiner operator (see \cite{GR}).
\end{remark}
We introduce the function  
\[
\rho(z,t) = (|z|^{4\gamma} + t^2)^{\frac{1}{4\gamma}},
\]
which, together with the horizontal gradient \(\nabla_\mathcal{L}\), is homogeneous of degree one with respect to the dilation group  
\[
\delta_\lambda(x,y,t) = (\lambda x, \lambda y, \lambda^{2\gamma} t).
\]
Setting \(Q = 2n + 2\gamma\) as the homogeneous dimension, it is well known (see, e.g., D'Ambrosio~\cite{DA05} and Zhang, Niu~\cite{ZN}) that for \(p \geq 2\),
\[
\begin{cases}  
\displaystyle \mathcal{L}_p\, \rho^{\frac{p-Q}{p-1}} = 0 \quad &\text{in } \mathbb{R}^{2n+1} \setminus \{0\}, \quad \text{if } p \neq Q, \\[5pt]  
\displaystyle \mathcal{L}_Q \,(-\ln \rho) = 0 \quad &\text{in } \mathbb{R}^{2n+1} \setminus \{0\}, \quad \text{if } p = Q.  
\end{cases}
\]
Finally, we observe that  
\[
|\nabla_\mathcal{L} \rho| = \frac{|z|^{2\gamma-1}}{\rho^{2\gamma-1}} \neq 0 \qquad \text{a.e. } (z,t) \in \mathbb{R}^{2n+1}.
\]
We now derive various Hardy-type inequalities within this framework. 

By applying Theorem~\ref{thm: sharp_d_radial_weighted_hardy_inequaity}, we obtain the following weighted Hardy inequality.
\begin{theorem}
Let \(\theta \in \mathbb{R}\), and assume \(\Omega \subseteq \mathbb{R}^N \setminus \{0\}\), then for all $u \in C^\infty_c(\Omega)$ one has
\begin{equation}
\int_\Omega \left| \nabla_\mathcal{L} u \cdot \frac{\nabla_\mathcal{L} \rho}{|\nabla_\mathcal{L} \rho|} \right|^p \frac{1}{\rho^{p(\theta - 1)}} \, dx 
\geq \left| \frac{Q - p \theta}{p} \right|^p \int_\Omega \frac{|u|^p}{\rho^{p \theta}} |\nabla_\mathcal{L} \rho|^p \, dx.
\end{equation}
Moreover, if \(\Omega \cup \{0\}\) contains a neighborhood of the origin, then the constant \(\left| \frac{Q - p \theta}{p} \right|^p\) is sharp.
\end{theorem}
Using Theorem~\ref{thm: sharp_logaritmic_hardy_inequality}, we can also establish a logarithmic Hardy inequality. 
\begin{theorem}
Let \(\theta \in \mathbb{R}\), and let \(\Omega \subseteq \mathbb{R}^N \setminus \{0\},\) then for all \(u \in C^\infty_c(\Omega)\) one has
\begin{equation}
\int_\Omega \left| \nabla_\mathcal{L} u \cdot \frac{\nabla_\mathcal{L} \rho}{|\nabla_\mathcal{L} \rho|} \right|^p \left( \ln \frac{R}{\rho} \right)^{\theta + p} \frac{1}{\rho^{Q-p}} \, dx 
\geq \left| \frac{\theta + 1}{p} \right|^p \int_\Omega \frac{|u|^p}{\rho^Q} \left( \ln \frac{R}{\rho} \right)^\theta |\nabla_\mathcal{L} \rho|^p \, dx.
\end{equation}
Moreover, if \(\Omega \cup \{0\}\) contains a neighborhood of the origin, then the constant is sharp.
\end{theorem}
Further refinements follow from the results in Section \ref{sezione: Gaussian-type Hardy inequality}, leading to Hardy-type inequalities involving Gaussian weights. More precisely, using~\ref{thm: hardy gaussiana tipo A} one gets
\begin{theorem}
Let \(\alpha\) and \(\beta\) be constants such that \(\alpha \geq 2\) and \(\beta > 0,\) then for all \(u \in C^\infty_c(\mathbb{R}^N)\) one has
\begin{multline*}
\int_{\mathbb{R}^N} \left| \nabla_\mathcal{L} u \cdot \frac{\nabla_\mathcal{L} \rho}{|\nabla_\mathcal{L} \rho|} \right|^p \displaystyle{e^{-\tfrac{\rho^{ \alpha}}{ \beta}}}\, dx 
\geq
\left( \frac{\alpha}{p \beta} \right)^p \int_{\mathbb{R}^N} |u|^p \rho^{p( \alpha - 1)} |\nabla_\mathcal{L} \rho|^p \displaystyle{e^{-\tfrac{\rho^{ \alpha}}{ \beta}}}\, dx\\ 
- \left( \frac{\alpha}{p \beta} \right)^{p-1} \left(\alpha  (p - 1) + Q - p\right) \int_{\mathbb{R}^N} |u|^p \rho^{ \alpha (p - 1) - p} |\nabla_\mathcal{L} \rho|^p \displaystyle{e^{-\tfrac{\rho^{ \alpha}}{ \beta}}}\, dx.
\end{multline*}
The candidate maximizing function is \(\displaystyle{e^{\frac{\rho^{ \alpha}}{p \beta}}}\). Moreover, this inequality is sharp in the sense that it cannot hold if the right-hand side is multiplied by any constant \(c > 1\). 
\end{theorem}
Now, Theorem~\ref{thm: hardy gaussiana tipo B} gives
\begin{theorem}
Let \(\theta \in \mathbb{R}\) and let \(\alpha, \beta\) be constants such that \(\alpha \geq 2\) and \(\beta > 0,\) then for all \(u \in C^\infty_c(\mathbb{R}^N \setminus \{0\}),\) one has
\begin{multline*}
\int_{\mathbb{R}^N} \left| \nabla_\mathcal{L} u \cdot \frac{\nabla_\mathcal{L} \rho}{|\nabla_\mathcal{L} \rho|} \right|^p \frac{1}{\rho^{p(\theta-1)}} \displaystyle{e^{-\tfrac{\rho^{ \alpha}}{ \beta}}} \, dx 
\geq
\left| \frac{Q - p \theta}{p} \right|^p \int_{\mathbb{R}^N} \frac{|u|^p}{\rho^{p \theta}} |\nabla_\mathcal{L} \rho|^p \displaystyle{e^{-\tfrac{\rho^{ \alpha}}{ \beta}}} \, dx\\ 
- \frac{\alpha}{\beta} \left(\frac{Q - p \theta}{p}\right) \left| \frac{Q - p \theta}{p} \right|^{p-2} \int_{\mathbb{R}^N} \frac{|u|^p}{\rho^{p \theta - \alpha}} |\nabla_\mathcal{L} \rho|^p \displaystyle{e^{-\tfrac{\rho^{ \alpha}}{ \beta}}} \, dx.
\end{multline*}
The candidate maximizing function is \(\rho^{-\frac{Q - p \theta}{p}}\).
Furthermore, if \(Q \neq p \theta\), the inequality is sharp in the sense that it does not hold if the right-hand side is multiplied by any constant \(c > 1\). 
\end{theorem}
Finally, applying Theorem~\ref{thm: sharp hardy inequality in an annulus} from Section~\ref{sezione: Hardy inequality in an Annulus}, we obtain a Hardy inequality in the annular domain
\[
\Omega_{a,b} = \{ (z,t) \in \mathbb{R}^{2n}\times\mathbb{R} \mid a < \rho(z,t) < b \}.
\]
More precisely, one has
\begin{theorem}
Let \(\theta \in \mathbb{R},\) then for all \(u \in C^\infty_c(\Omega_{a,b})\)
one has
\begin{equation}
\int_{\Omega_{a,b}} \left|\nabla_\mathcal{L} u \cdot \frac{\nabla_\mathcal{L} \rho}{|\nabla_\mathcal{L} \rho|} \right|^p \frac{1}{\rho^{p(\theta-1)}} \, dx \geq \lambda_1 \int_{\Omega_{a,b}} \frac{|u|^p}{\rho^{p\theta}} |\nabla_\mathcal{L} \rho|^p \, dx,
\end{equation}
where \(\lambda_1\) denotes the first eigenvalue of~\eqref{eq: problema autovalori A}.
The candidate maximizing function is the corresponding eigenfunction \(\varphi_1\). Moreover, the constant is sharp. 
\end{theorem}
Remark~\ref{rmk:first-eigen-prop} holds also in this context.
\smallskip
We observe that the matrix \(\sigma\) defining \(\nabla_\mathcal{L}\) conforms to the structure given in \eqref{eq: matrice sigma cilindrica}, namely  one has
\[
\sigma = \begin{pmatrix} \mathcal{I}_{2n} & \sigma_1 \end{pmatrix},
\]
where  
\[
\sigma_1 = \gamma |z|^{2\gamma-2} \begin{pmatrix} 2y \\ 2x \end{pmatrix} \in \mathcal{M}_{2n \times 1}.
\]
Notice that \(\sigma_1\) is independent of the variable \(t\). Moreover
\[
z \cdot (\sigma_1 t) = (x,y) \cdot ( 2\gamma |z|^{2\gamma-2} y t, -2 \gamma |z|^{2\gamma-2} x t) = 0.
\]
This gives the following two results which are consequences of Theorem~\ref{thm: sharp_cilindrical_hardy_inequality} and Theorem~\ref{thm: Cylindrical logaritmic hardy inequality}, respectively.
\begin{theorem}
Fix \(\theta \in \mathbb{R},\) and assume \(\Omega \subseteq \mathbb{R}^{2n} \setminus \{0\} \times \mathbb{R},\) then for all \(u \in C^\infty_c(\Omega)\) one has
\begin{equation*}
\int_\Omega \left| \nabla_\mathcal{L} u \cdot  \frac{z}{|z|}  \right|^p \frac{1}{|z|^{p(\theta - 1)}} \, dz \, dt 
\geq \left| \frac{2n - p \theta}{p} \right|^p \int_\Omega \frac{|u(z, t)|^p}{|z|^{p \theta}} \, dz \, dt.
\end{equation*}
Moreover, if \(B^{2n}_r \setminus \{0\} \times B^{1}_r \subseteq \Omega\) for some \(r > 0\), then the constant is sharp.
\end{theorem}
\begin{theorem}
Let \(R > 0\), fix \(\,\theta \in \mathbb{R}\), and assume \(\Omega = B^{2n}_R \setminus \{0\} \times \mathbb{R},\) then for all \(u \in C^\infty_c(\Omega)\) one has
\begin{equation*}
\int_\Omega \left| \nabla_\mathcal{L} u \cdot \frac{z}{|z|} \right|^p \frac{1}{|z|^{2n-p}} \left(\ln\left(\frac{R}{|z|}\right)\right)^{\theta+p} \, dz \, dt 
\geq \left|\frac{\theta+1}{p}\right|^p \int_\Omega \frac{|u|^p}{|z|^{2n}} \left(\ln\left(\frac{R}{|z|}\right)\right)^\theta \, dz \, dt.
\end{equation*}
Moreover, the constant is sharp.
\end{theorem}

\subsection{Baouendi-Grushin Operator}
We split \(\mathbb{R}^N\) into \(\mathbb{R}^n \times \mathbb{R}^k\), introducing coordinates \(x \in \mathbb{R}^n\) and \(y \in \mathbb{R}^k\). Given \(\gamma \geq 0\), we define the vector field  
\[
\nabla_\gamma = (\nabla_x, (1+\gamma) |x|^\gamma \nabla_y),
\]
which induces the second-order operator  
\[
\mathcal{L}_2 = \Delta_x + (1+\gamma)^2 |x|^{2\gamma} \Delta_y,
\]
known as the Baouendi–Grushin operator. The matrix associated with \(\nabla_\gamma\) is  
\[
\sigma = \begin{pmatrix} \mathcal{I}_n & 0 \\ 0 & (1+\gamma) |x|^\gamma \mathcal{I}_k \end{pmatrix}.
\]
\begin{remark}
If \( k=0 \) or \(\gamma=0\), the operator \(\mathcal{L}_2\) and the corresponding \(\mathcal{L}_p\) coincide with the standard Euclidean Laplacian and \(p\)-Laplacian.
\end{remark}
We define the function \(\rho(x,y) = ( |x|^{2(1+\gamma)} + |y|^2 )^{\tfrac{1 }{ 2(1+\gamma)}}\) and introduce the dilation family 
\[\delta_\lambda(x,y) = (\lambda x, \lambda^{1+\gamma} y),
\qquad
\lambda > 0.
\]
Both \(\nabla_\gamma\) and \(\rho(x,y)\) are homogeneous of degree one with respect to \(\delta_\lambda\), and the homogeneous dimension is given by \( Q = n + (1+\gamma) k \). 
A straightforward computation shows that  
\[
\begin{cases} 
\displaystyle \mathcal{L}_p \, \rho^{\frac{p-Q}{p-1}} = 0 \quad \text{in } \mathbb{R}^N \setminus \{0\}, & \text{if } p \neq Q, \\[5pt]
\displaystyle \mathcal{L}_Q (-\ln \rho) = 0 \quad \text{in } \mathbb{R}^N \setminus \{0\}, & \text{if } p = Q.
\end{cases}
\]
Finally, we observe that \( |\nabla_\gamma \rho| = |x|^\gamma / \rho^\gamma \).

The results from the previous sections allow us to derive several Hardy-type inequalities. In particular, Theorem~\ref{thm: sharp_d_radial_weighted_hardy_inequaity} yields a classical Hardy inequality.  
\begin{theorem}
Let \(\theta \in \mathbb{R}\), and let \(\Omega \subseteq \mathbb{R}^N \setminus \{0\}\), then for all $u \in C^\infty_c(\Omega)$ one has
\begin{equation}
\int_\Omega \left| \nabla_\gamma u \cdot \frac{\nabla_\gamma \rho}{|\nabla_\gamma \rho|} \right|^p \frac{1}{\rho^{p(\theta - 1)}} \, dx \, dy 
\geq \left| \frac{Q - p \theta}{p} \right|^p \int_\Omega \frac{|u|^p}{\rho^{p \theta}} |\nabla_\gamma \rho|^p \, dx \, dy.
\end{equation}
Moreover, if \(\Omega \cup \{0\}\) contains a neighborhood of the origin, then the constant is sharp.
\end{theorem}
Similarly, Theorem~\ref{thm: sharp_logaritmic_hardy_inequality} provides a logarithmic Hardy-type inequality.  
\begin{theorem}
Let \(\theta \in \mathbb{R}\), and let \(\Omega \subseteq \mathbb{R}^N \setminus \{0\},\) then for all \(u \in C^\infty_c(\Omega)\) one as
\begin{equation}
\int_\Omega \left| \nabla_\gamma u \cdot \frac{\nabla_\gamma \rho}{|\nabla_\gamma \rho|} \right|^p \left( \ln \frac{R}{\rho} \right)^{\theta + p} \frac{1}{\rho^{Q-p}} \, dx \, dy
\geq \left| \frac{\theta + 1}{p} \right|^p \int_\Omega \frac{|u|^p}{\rho^Q} \left( \ln \frac{R}{\rho} \right)^\theta |\nabla_\gamma \rho|^p \, dx \, dy.
\end{equation}
Moreover, if \(\Omega \cup \{0\}\) contains a neighborhood of the origin, then the constant is sharp.
\end{theorem}
By applying the results from Section~\ref{sezione: Gaussian-type Hardy inequality}, namely Theorem~\ref{thm: hardy gaussiana tipo A} and Theorem~\ref{thm: hardy gaussiana tipo B}, we obtain Hardy-type inequalities involving Gaussian weights.  
\begin{theorem}
Let \(\alpha\) and \(\beta\) be constants such that $\alpha\geq 2$ and \(\beta > 0,\) then for all \(u \in C^\infty_c(\mathbb{R}^N)\) one has
\begin{multline*}
\int_{\mathbb{R}^N} \left| \nabla_\gamma u \cdot \frac{\nabla_\gamma \rho}{|\nabla_\gamma \rho|} \right|^p \displaystyle{e^{-\tfrac{\rho^{ \alpha}}{ \beta}}}\, dx \, dy 
\geq 
\left( \frac{\alpha}{p \beta} \right)^p \int_{\mathbb{R}^N} |u|^p \rho^{p( \alpha - 1)} |\nabla_\gamma \rho|^p \displaystyle{e^{-\tfrac{\rho^{ \alpha}}{ \beta}}}\, dx \, dy\\ 
- \left( \frac{\alpha}{p \beta} \right)^{p-1} \left(\alpha  (p - 1) + Q - p\right) \int_{\mathbb{R}^N} |u|^p \rho^{ \alpha (p - 1) - p} |\nabla_\gamma \rho|^p \displaystyle{e^{-\tfrac{\rho^{ \alpha}}{ \beta}}}\, dx \, dy.
\end{multline*}
The candidate maximizing function is \(\displaystyle{e^{\frac{\rho^{ \alpha}}{p \beta}}}\).
Moreover, this inequality is sharp in the sense that it cannot hold if the right-hand side is multiplied by any constant \(c > 1\). 
\end{theorem}
\begin{theorem}
Let \(\theta \in \mathbb{R}\) and let \(\alpha, \beta\) be constants such that \(\alpha \geq 2\) and \(\beta > 0,\) then for all \(u \in C^\infty_c(\mathbb{R}^N \setminus \{0\})\) one has
\begin{multline*}
\int_{\mathbb{R}^N} \left| \nabla_\gamma u \cdot \frac{\nabla_\gamma \rho}{|\nabla_\gamma \rho|} \right|^p \frac{1}{\rho^{p(\theta-1)}} \displaystyle{e^{-\tfrac{\rho^{ \alpha}}{ \beta}}} \, dx \, dy \geq
\left| \frac{Q - p \theta}{p} \right|^p \int_{\mathbb{R}^N} \frac{|u|^p}{\rho^{p \theta}} |\nabla_\gamma \rho|^p \displaystyle{e^{-\tfrac{\rho^{ \alpha}}{ \beta}}} \, dx \, dy\\ 
- \frac{\alpha}{\beta} \left(\frac{Q - p \theta}{p}\right) \left| \frac{Q - p \theta}{p} \right|^{p-2} \int_{\mathbb{R}^N} \frac{|u|^p}{\rho^{p \theta - \alpha}} |\nabla_\gamma \rho|^p \displaystyle{e^{-\tfrac{\rho^{ \alpha}}{ \beta}}} \, dx \, dy.\nonumber
\end{multline*}
The candidate maximizing function is \(\rho^{-\frac{Q - p \theta}{p}}\).
Moreover, if \(Q \neq p \theta\), the inequality is sharp in the sense that it does not hold if the right-hand side is multiplied by any constant \(c > 1\). 
\end{theorem}
Moreover, we establish a Hardy inequality with an optimal constant in the annular region  
\[
\Omega_{a,b} = \{ (x,y) \in \mathbb{R}^n\times\mathbb{R}^k \mid a < \rho(x,y) < b \}.
\]  
This will be direct consequence of Theorem~\ref{thm: sharp hardy inequality in an annulus}.
\begin{theorem}
Let \(\theta \in \mathbb{R},\) then for all \(u \in C^\infty_c(\Omega_{a,b})\) one has
\begin{equation}
\int_{\Omega_{a,b}} \left|\nabla_\gamma u \cdot \frac{\nabla_\gamma \rho}{|\nabla_\gamma \rho|} \right|^p \frac{1}{\rho^{p(\theta-1)}} \, dx \, dy \geq \lambda_1 \int_{\Omega_{a,b}} \frac{|u|^p}{\rho^{p\theta}} |\nabla_\gamma \rho|^p \, dx \, dy,
\end{equation}
where \(\lambda_1\) denotes the first eigenvalue of \eqref{eq: problema autovalori A}.
The candidate maximizing function is the corresponding eigenfunction \(\varphi_1\). Moreover, the constant is sharp.
\end{theorem}
Notice that Remark~\ref{rmk:first-eigen-prop} holds also in this context.

\medskip
In order to apply the results from Section \ref{sezione: sharp_cylindrical_hardy_inequality} and Section \ref{sezione: Cylindrical Logarithmic Hardy Inequality}, let us observe that if we fix \( m \) such that \( 1 \leq m \leq n \), we can express \(\sigma\) as  
\[
\sigma = \begin{pmatrix} \mathcal{I}_m & 0 & 0 \\ 0 & \mathcal{I}_{n-m} & 0 \\ 0 & 0 & (1+\gamma)|x|^\gamma \mathcal{I}_k \end{pmatrix}  
= \begin{pmatrix} \mathcal{I}_m & \sigma_1 \\ 0 & \sigma_2 \end{pmatrix},
\]
where  
\[
\sigma_1 = 0 \in \mathcal{M}_{m \times (N-m)}, \quad \text{and} \quad \sigma_2 = \begin{pmatrix} \mathcal{I}_{n-m} & 0 \\ 0 & (1+\gamma) |x|^\gamma \mathcal{I}_k \end{pmatrix} \in \mathcal{M}_{(N-m) \times (N-m)}.
\]
Introducing the notation \( \xi = (x_1, \dots, x_m) \) and \( \eta = (x_{m+1}, \dots, x_n, y_1, \dots, y_k) \), we obtain the following results which are consequence of Theorem~\ref{thm: sharp_cilindrical_hardy_inequality} and Theorem~\ref{thm: Cylindrical logaritmic hardy inequality}.  
\begin{theorem}
Fix  \(\theta \in \mathbb{R}\) and assume \(\Omega \subseteq \mathbb{R}^m \setminus \{0\} \times \mathbb{R}^{N-m}\), then for all \(u \in C^\infty_c(\Omega)\) one has
\begin{equation*}
\int_\Omega \left| \nabla_\gamma u \cdot \left( \frac{\xi}{|\xi|}, 0_{N-m} \right) \right|^p \frac{1}{|\xi|^{p(\theta - 1)}} \, d\xi \, d\eta 
\geq \left| \frac{m - p \theta}{p} \right|^p \int_\Omega \frac{|u|^p}{|\xi|^{p \theta}} \, d\xi \, d\eta.
\end{equation*}
Moreover, if \(B^m_r \setminus \{0\} \times B^{N-m}_r \subseteq \Omega\) for some \(r > 0\), then the constant is sharp.
\end{theorem}
\begin{theorem}
Let \(R > 0\), fix \(\,\theta \in \mathbb{R}\), and assume \(\Omega = B^m_R \setminus \{0\} \times \mathbb{R}^{N-m},\) then for all \(u \in C^\infty_c(\Omega)\) one has
\begin{equation}
\int_\Omega \left| \nabla_\gamma u \cdot \left(\frac{\xi}{|\xi|}, 0_{N-m}\right) \right|^p \frac{1}{|\xi|^{m-p}} \left(\ln\left(\frac{R}{|\xi|}\right)\right)^{\theta+p} \, d\xi \, d\eta 
\geq \left|\frac{\theta+1}{p}\right|^p \int_\Omega \frac{|u|^p}{|\xi|^m} \left(\ln\left(\frac{R}{|\xi|}\right)\right)^\theta \, d\xi \, d\eta.
\end{equation}
Moreover, the constant is sharp.
\end{theorem}

\subsection{Hardy-Type Inequalities on Carnot Groups}
In this section, we establish Hardy-type inequalities within the framework of Carnot groups. We first recall some preliminary results and refer the interested reader to \cite{BLU} for a more comprehensive treatment.
\begin{definition}  
A Lie group \( G = (\mathbb{R}^N, \circ) \) is called a \textit{homogeneous Carnot group} if the following properties hold: 
\begin{enumerate}  
\item \(\mathbb{R}^N\) decomposes as \( \mathbb{R}^N = \mathbb{R}^{N_1} \times \dots \times \mathbb{R}^{N_r} \), and the family of dilations  
\[
\delta_\lambda (x^{(1)}, \dots, x^{(r)}) = (\lambda x^{(1)}, \lambda^2 x^{(2)}, \dots, \lambda^r x^{(r)}), \quad x^{(i)} \in \mathbb{R}^{N_i},
\]
is an automorphism of \( G \) for every \( \lambda > 0 \).  

\item Letting \( N_1 \) be as above, consider the left-invariant vector fields \( X_1, \dots, X_{N_1} \) on \( G \) such that  
\(\displaystyle
X_j(0) = \frac{\partial}{\partial x_j} \Big|_0,\)  \(j = 1, \dots, N_1.
\)  
Then, 
\[
\operatorname{rank} (\operatorname{Lie} \{X_1, \dots, X_{N_1} \} (x)) = N, \quad \text{for every } x \in \mathbb{R}^N.
\] 
\end{enumerate}  
\end{definition}  
We denote by \(\nabla_\mathcal{L} = (X_1, \dots, X_{N_1})\). Among the fundamental properties of Carnot groups, the vector field \(\nabla_\mathcal{L}\) is homogeneous of degree one with respect to the dilation family \(\delta_\lambda\).  
The canonical sub-Laplacian on \(G\) is the second-order differential operator defined as
\[
\mathcal{L}_2 = \sum_{i=1}^{N_1} X_i^2.
\]  
For \( p \geq 2 \), the \( p \)-sub-Laplacian is given by  
\[
\mathcal{L}_p u = \sum_{i=1}^{N_1} X_i \big( |\nabla_\mathcal{L} u |^{p-2} X_i u \big) = \operatorname{div}_\mathcal{L} ( |\nabla_\mathcal{L} u |^{p-2} \nabla_\mathcal{L} u).
\]  
%
Examples of Carnot groups include the Euclidean group \((\mathbb{R}^N, +)\) and the Heisenberg group \((\mathbb{H}_n, \circ_\mathbb{H})\) (see~\cite{BLU}). Furthermore, denoting by  
\(
Q = N_1 + 2N_2 + \dots + rN_r
\)
the homogeneous dimension, it follows that if \( Q < 3 \), the Carnot group reduces to the Euclidean group \(\mathbb{R}^Q\). Thus, from now on, we may assume \( Q \geq 3 \).
\vspace{0.15cm}

We call a \textit{homogeneous norm} on \( G \) any continuous function \( \rho \colon G \to [0, +\infty) \) satisfying the following properties
\begin{itemize}  
\item \( \rho(\delta_\lambda(x)) = \lambda \rho(x) \) for all \( \lambda > 0 \) and \( x \in G \);  
\item \( \rho(x) > 0 \) if and only if \( x \neq 0 \).  
\end{itemize}  
It is well known that there always exists a homogeneous norm \( \rho_2 \in C^\infty(G \setminus \{0\}) \cap C(G) \) related to the fundamental solution of \( \mathcal{L}_2 \), satisfying 
\[
\mathcal{L}_2 \rho_2^{2-Q} = 0 \quad \text{in } G \setminus \{0\}.
\]
However, there is no general guarantee that \( \rho_2 \) is also related to the fundamental solution of \( \mathcal{L}_p \), as is the case for the Heisenberg group.
We therefore introduce the following definition (see~\cite{DA05}).  
\begin{definition}  
A Carnot group \( G \) is said to be \textit{polarizable} if, \textit{for every} \( p \geq 2 \), the function \( \rho_2 \) satisfies  
\[
\begin{cases}  
\displaystyle \mathcal{L}_p \,\rho_2^{\frac{p-Q}{p-1}} = 0 \quad \text{in } G \setminus \{0\}, & \text{if } p \neq Q, \\[5pt]  
\displaystyle \mathcal{L}_Q  \, (-\ln \rho_2) = 0 \quad \text{in } G \setminus \{0\}, & \text{if } p = Q.  
\end{cases}  
\]
\end{definition}  
Among the examples of polarizable Carnot groups we can find Euclidean spaces, the Heisenberg group, and H-type groups (see~\cite{BT}).

\medskip
We are now in position to establish several Hardy-type inequalities in the framework of Carnot groups. 

The classical Hardy inequality (see \cite{DA05}) follows from Theorem~\ref{thm: sharp_d_radial_weighted_hardy_inequaity}, while a logarithmic counterpart can be derived from Theorem~\ref{thm: sharp_logaritmic_hardy_inequality}. Additionally, Hardy inequalities involving Gaussian weights are obtained from Theorem~\ref{thm: hardy gaussiana tipo A} and Theorem~\ref{thm: hardy gaussiana tipo B}. 

A further refinement is given by the Hardy inequality in the annular domain  
\[
\Omega_{a,b}=\{x\in\mathbb{R}^N \mid a<\rho_2(x)<b\}.
\]  
We shall distinguish between the case \( p=2 \), where these inequalities hold for \textit{any} Carnot group, and the case \( p>2 \), where our techniques allow to establish them only in \textit{polarizable} Carnot groups. 
\begin{theorem}
Let \(\theta \in \mathbb{R}\), and assume \(\Omega \subseteq \mathbb{R}^N \setminus \{0\}\). Let \( p \geq 2.\) If $p=2$ assume \( G \) to be any Carnot group,  if \( p > 2 \) assume \( G \)  to be a polarizable Carnot group,
then for all $u \in C^\infty_c(\Omega)$ one has
\begin{equation*}
\int_\Omega \left| \nabla_\mathcal{L} u \cdot \frac{\nabla_\mathcal{L} \rho_2}{|\nabla_\mathcal{L} \rho_2|} \right|^p \frac{1}{d^{p(\theta - 1)}} \, dx 
\geq \left| \frac{Q - p \theta}{p} \right|^p \int_\Omega \frac{|u|^p}{\rho_2^{p \theta}} |\nabla_\mathcal{L} \rho_2|^p \, dx.
\end{equation*}
Moreover, if \(\Omega \cup \{0\}\) contains a neighborhood of the origin, then the constant is sharp.
\end{theorem}
\begin{theorem}
Let \(\theta \in \mathbb{R}\) and \(\Omega \subseteq \mathbb{R}^N \setminus \{0\}\). Let \( p \geq 2.\) If $p=2$ assume \( G \) to be any Carnot group,  if \( p > 2 \) assume \( G \)  to be a polarizable Carnot group,
then for all $u \in C^\infty_c(\Omega)$ one has
\begin{equation*}
\int_\Omega \left| \nabla_\mathcal{L} u \cdot \frac{\nabla_\mathcal{L} \rho_2}{|\nabla_\mathcal{L} \rho_2|} \right|^p \left( \ln \frac{R}{\rho_2} \right)^{\theta + p} \frac{1}{\rho_2^{Q-p}} \, dx 
\geq \left| \frac{\theta + 1}{p} \right|^p \int_\Omega \frac{|u|^p}{\rho_2^Q} \left( \ln \frac{R}{\rho_2} \right)^\theta |\nabla_\mathcal{L} \rho_2|^p \, dx.
\end{equation*}
Moreover, if \(\Omega \cup \{0\}\) contains a neighborhood of the origin, then the constant is sharp.
\end{theorem}
\begin{theorem}
Let \(\alpha\) and \(\beta\) be constants such that  \(\alpha \geq 2\) and \(\beta > 0.\) Let \( p \geq 2.\) If $p=2$ assume \( G \) to be any Carnot group,  if \( p > 2 \) assume \( G \)  to be a polarizable Carnot group,
then for all $u \in C^\infty_c(\R^N)$ one has
\begin{multline*}
\int_{\mathbb{R}^N} \left| \nabla_\mathcal{L} u \cdot \frac{\nabla_\mathcal{L} \rho_2}{|\nabla_\mathcal{L} \rho_2|} \right|^p \displaystyle{e^{-\tfrac{\rho_2^{ \alpha}}{ \beta}}}\, dx 
\geq
\left( \frac{\alpha}{p \beta} \right)^p \int_{\mathbb{R}^N} |u|^p \rho_2^{p( \alpha - 1)} |\nabla_\mathcal{L} \rho_2|^p \displaystyle{e^{-\tfrac{\rho_2^{ \alpha}}{ \beta}}}\, dx \\
- \left( \frac{\alpha}{p \beta} \right)^{p-1} \left(\alpha  (p - 1) + Q - p\right) \int_{\mathbb{R}^N} |u|^p \rho_2^{ \alpha (p - 1) - p} |\nabla_\mathcal{L} \rho_2|^p \displaystyle{e^{-\tfrac{\rho_2^{ \alpha}}{ \beta}}}\, dx.
\end{multline*}
The candidate maximizing function is \(\displaystyle{e^{\frac{\rho_2^{ \alpha}}{p \beta}}}.\)
Moreover, this inequality is sharp in the sense that it cannot hold if the right-hand side is multiplied by any constant \(c > 1.\)
\end{theorem}
\begin{theorem}
Let \(\theta \in \mathbb{R}\) and let \(\alpha, \beta\) be two constants with \(\alpha \geq 2\) and \(\beta > 0\). Let \( p \geq 2.\) If $p=2$ assume \( G \) to be any Carnot group,  if \( p > 2 \) assume \( G \)  to be a polarizable Carnot group,
then for all $u \in C^\infty_c(\mathbb{R}^N \setminus \{0\})$ one has
\begin{multline*}
\int_{\mathbb{R}^N} \left| \nabla_\mathcal{L} u \cdot \frac{\nabla_\mathcal{L} \rho_2}{|\nabla_\mathcal{L} \rho_2|} \right|^p \frac{1}{\rho_2^{p(\theta-1)}} \displaystyle{e^{-\tfrac{\rho_2^{ \alpha}}{ \beta}}} \, dx \geq
\left| \frac{Q - p \theta}{p} \right|^p \int_{\mathbb{R}^N} \frac{|u|^p}{\rho_2^{p \theta}} |\nabla_\mathcal{L} \rho_2|^p \displaystyle{e^{-\tfrac{\rho_2^{ \alpha}}{ \beta}}} \, dx\\ 
- \frac{\alpha}{\beta} \left(\frac{Q - p \theta}{p}\right) \left| \frac{Q - p \theta}{p} \right|^{p-2} \int_{\mathbb{R}^N} \frac{|u|^p}{\rho_2^{p \theta - \alpha}} |\nabla_\mathcal{L} \rho_2|^p \displaystyle{e^{-\tfrac{\rho_2^{ \alpha}}{ \beta}}} \, dx.\nonumber
\end{multline*}
The candidate maximizing function is \(\rho_2^{-\frac{Q - p \theta}{p}}\). Moreover, if \(Q \neq p \theta\), the inequality is sharp in the sense that it does not hold if the right-hand side is multiplied by any constant \(c > 1\). 
\end{theorem}
\begin{theorem}
Let $\theta\in \R.$ Let \( p \geq 2.\) If $p=2$ assume \( G \) to be any Carnot group,  if \( p > 2 \) assume \( G \)  to be a polarizable Carnot group,
then for all $u \in C^\infty_c(\Omega_{a,b})$ one has
\begin{equation*}
\int_{\Omega_{a,b}} \left|\nabla_\mathcal{L} u \cdot \frac{\nabla_\mathcal{L} \rho_2}{|\nabla_\mathcal{L} \rho_2|} \right|^p \frac{1}{\rho_2^{p(\theta-1)}} \, dx \geq \lambda_1 \int_{\Omega_{a,b}} \frac{|u|^p}{\rho_2^{p\theta}} |\nabla_\mathcal{L} \rho_2|^p \, dx,
\end{equation*}
where \(\lambda_1\) denotes the first eigenvalue of \eqref{eq: problema autovalori A}.
The candidate maximizing function is the corresponding eigenfunction \(\varphi_1\). Moreover, the constant is sharp.
\end{theorem}
Notice that Remark~\ref{rmk:first-eigen-prop} holds in this setting too.

\medskip
In a homogeneous Carnot group, the first-layer vector fields \( X_1, \dots, X_{N_1} \) are of the form 
\[
X_i = \frac{\partial}{\partial x_i} + \sum_{h=2}^r \sum_{k=1}^{N_h} a^h_{i,k}(x^{(1)}, \dots, x^{(h-1)}) \frac{\partial}{\partial x_k^{(h)}}, \qquad i = 1, \dots, N_1, \qquad x^{(j)} \in \mathbb{R}^{N_j},
\]  
where \( a^h_{i,k} \) is a homogeneous polynomial of degree \( h-1 \) under the dilations \(\delta_\lambda\) (see~ \cite[Rmk.1.4.6]{BLU}). As a consequence, the matrix \(\sigma\) that defines \(\nabla_\mathcal{L}\) is given by
\[
\sigma = \begin{pmatrix} \mathcal{I}_{N_1} & \sigma_1 \end{pmatrix}, \quad \sigma_1 \in \mathcal{M}_{N_1 \times (N - N_1)}.
\]
From Theorem~\ref{thm: sharp_cilindrical_hardy_inequality} and Theorem~\ref{thm: Cylindrical logaritmic hardy inequality}, we obtain the following results:  
\begin{theorem}  
Let \( G \) be any Carnot group and fix \(\theta \in \mathbb{R}\). If \(\Omega \subseteq \mathbb{R}^{N_1} \setminus \{0\} \times \mathbb{R}^{N-N_1}\), then for all \(u \in C^\infty_c(\Omega)\) one has
\begin{equation}
\label{eq: hardy cilindrica in Carnot}
\int_\Omega \left| \nabla_\mathcal{L} u \cdot  \frac{x^{(1)}}{|x^{(1)}|}  \right|^p \frac{1}{|x^{(1)}|^{p(\theta - 1)}} \, dx
\geq \left| \frac{N_1 - p \theta}{p} \right|^p \int_\Omega \frac{|u|^p}{|x^{(1)}|^{p \theta}} \, dx
\end{equation}
\end{theorem}  
\begin{theorem}  
Let \( G \) be any Carnot group let \(R > 0\), \(\,\theta \in \mathbb{R}\), and \(\Omega = B^{N_1}_R \setminus \{0\} \times \mathbb{R}^{N-N_1},\) then for all \(u \in C^\infty_c(\Omega)\) one has
\begin{equation}
\label{eq: hardy logarotmica cilindrica Carnot}
\int_\Omega \left| \nabla_\mathcal{L} u \cdot \frac{x^{(1)}}{|x^{(1)}|} \right|^p \frac{1}{|x^{(1)}|^{N_1-p}} \left(\ln\left(\frac{R}{|x^{(1)}|}\right)\right)^{\theta+p} \, dx
\geq \left|\frac{\theta+1}{p}\right|^p \int_\Omega \frac{|u|^p}{|x^{(1)}|^{N_1}} \left(\ln\left(\frac{R}{|x^{(1)}|}\right)\right)^\theta \, dx
\end{equation}
\end{theorem}  
\vspace{0.2cm}
In this general context, we are unable to determine the optimality of the constant due to the lack of information regarding the matrix \(\sigma_1\). However, as previously noted, when \(G\) is either the Euclidean group or the Heisenberg group, the constants in \eqref{eq: hardy cilindrica in Carnot} and \eqref{eq: hardy logarotmica cilindrica Carnot} are sharp. 

\medskip
In what follows, we extend this result to a class of step-two Carnot groups, including, for instance, H-type groups.

Let us recall a general characterization of step-two Carnot groups. Let \(\mathbb{R}^N = \mathbb{R}^{N_1} \times \mathbb{R}^{N_2}\), with coordinates \((x, t)\), with \(x \in \mathbb{R}^{N_1}\) and \(t \in \mathbb{R}^{N_2}\). Given an \(N_2\)-tuple of \(N_1 \times N_1\) matrices \( B^{(1)}, \dots, B^{(N_2)} \), the group law is defined as  
\[
(x,t) \circ (\xi, \tau) = (x + \xi, t + \tau + \frac{1}{2} \langle Bx, \xi \rangle),
\]  
where \( \langle Bx, \xi \rangle \) represents the \(N_2\)-tuple \( \left( \langle B^{(1)} x, \xi \rangle, \dots, \langle B^{(N_2)} x, \xi \rangle \right) \).

It is straightforward to verify that the family of dilations  
\[
 \delta_\lambda(x,t) = (\lambda x, \lambda^2 t) 
\]
is an automorphism of the group \((\mathbb{R}^N, \circ)\) for all \(\lambda > 0\). The first-layer vector fields are explicitly given by  
\[
X_i = \frac{\partial}{\partial x_i} + \frac{1}{2} \langle (Bx)_i, \nabla_t \rangle, \quad i = 1, \dots, N_1,
\]  
where we denote by \((Bx)_i = ((B^{(1)}x)_i, \dots, (B^{(N_2)}x)_i)\).
\begin{proposition}[Proposition 3.2.1 in\cite{BLU}]
A characterization of homogeneous Carnot groups of step two and \( N_1 \) generators is given by the above \( G = (\mathbb{R}^{N_1 + N_2}, \circ, \delta_\lambda) \), where the skew-symmetric parts of the matrices \( B^{(k)} \) are linearly independent. 
\end{proposition}
Examples of step-two Carnot groups include the Heisenberg group and H-type groups, also known as Heisenberg-type groups. Below, we recall their definition.
\begin{definition}
Consider the homogeneous Lie group \( H = (\mathbb{R}^{N_1 + N_2}, \circ, \delta_\lambda) \) with the group law 
\[
(x,t) \circ (\xi, \tau) = \left( x + \xi, \, t_1 + \tau_1 + \frac{1}{2} \langle B^{(1)} x, \xi \rangle, \dots, \, t_{N_2} + \tau_{N_2} + \frac{1}{2} \langle B^{(N_2)} x, \xi \rangle \right),
\]  
where \( B^{(1)}, \dots, B^{(N_2)} \) are fixed \( N_1 \times N_1 \) matrices and the dilations are given by 
\(
\delta_\lambda(x,t) = (\lambda x, \lambda^2 t).
\) 
We also assume the following conditions
\begin{enumerate}  
\item Each \( B^{(j)} \) is a skew-symmetric orthogonal matrix of size \( N_1 \times N_1 \), for \( j \leq N_2 \);
\item \( B^{(i)} B^{(j)} = - B^{(j)} B^{(i)} \) for all \( i, j \in \{1, \dots, N_2\} \), with \( i \neq j \).  
\end{enumerate}  
If these conditions are satisfied, the group \( H \) is called a Heisenberg-type group, or more briefly, an H-type group.
\end{definition}
The following theorem can now be stated.
\begin{theorem}
Let \( G \) be a step-two Carnot group, and assume that the matrices \( B^{(1)}, \dots, B^{(N_2)} \) are skew-symmetric. Then, the constant in the inequality \eqref{eq: hardy logarotmica cilindrica Carnot} is sharp. Moreover, if we denote by \( B^{l}_r \) the Euclidean ball centered at the origin with radius \( r \) in \( \mathbb{R}^l \), if \( B^{N_1}_r \setminus \{0\} \times B^{N_2}_r \subseteq \Omega \) for some \( r > 0 \), the constant in the Hardy inequality \eqref{eq: hardy cilindrica in Carnot} is sharp.
\end{theorem}
\begin{remark}
This result applies, in particular, to H-type groups.
\end{remark}
\begin{proof}
In a homogeneous step-two Carnot group, due to the specific form of the vector fields \( X_1, \dots, X_{N_1} \), the matrix \( \sigma \) takes the form  
\[
\sigma = \begin{pmatrix} \mathcal{I}_{N_1} & \sigma_1 \end{pmatrix}, \quad \text{where} \quad (\sigma_1)_{i,j} = \frac{1}{2}(B^{(j)} x)_i, \qquad i = 1, \dots, N_1, \quad  j = 1, \dots, N_2.
\]  
In particular, we have  
\[
x \cdot (\sigma_1 t) =\frac{1}{2} t_1 \langle B^{(1)} x, x \rangle + \dots +\frac{1}{2} t_{N_2} \langle B^{(N_2)} x, x \rangle = 0,
\]  
since the matrices \( B^{(1)}, \dots, B^{(N_2)} \) are assumed to be skew-symmetric. This condition is sufficient to establish the optimality of the constants from Theorem~\ref{thm: sharp_cilindrical_hardy_inequality} and Theorem~\ref{thm: Cylindrical logaritmic hardy inequality}.
\end{proof}

\appendix
\section{Proof of auxiliary results}
\label{app:auxiliary}
In this appendix we collect the proofs of several auxiliary results that have been used to obtain our main results. We shall start with the proof of Proposition~\ref{prop:identita_fondamentale_campivettoriali}.
\begin{proof}[Proof of Prop.~\ref{prop:identita_fondamentale_campivettoriali}]
In order to prove~\eqref{eq: identita_fondamentale_campivettoriali} in Proposition~\ref{prop:identita_fondamentale_campivettoriali} it is enough to show the following identity
\begin{multline}
\label{eq:concise_identity}
|\zeta|^p + (p-1)|\xi|^p - p |\xi|^{p-2} \Re \langle \xi,\zeta \rangle_{\C^h}
=p|\zeta-\xi|^2 \int_0^1 s\, |s \xi + (1-s) \zeta|^{p-2}\, ds\\ 
+ p(p-2) \int_0^1 s\, |s \xi + (1-s) \zeta|^{p-4} \big[\Re\langle \zeta-\xi,s \xi + (1-s) \zeta\rangle_{\C^h}\big]^2 ds,
\end{multline}
for any \(\zeta,\xi \in \mathbb{C}^h\). Indeed, then~\eqref{eq: identita_fondamentale_campivettoriali} follows from~\eqref{eq:concise_identity} replacing $\zeta$ with $f(x)$ and $\xi$ with $g(x)$ and then integrating over \(\Omega\).

Using the standard identification $\C^h \backsimeq \R^{2h},$ that is defining the $2h$-components vectors $\mu$ and $\nu$ associated to $\zeta$ and $\xi$ as
\begin{equation*}
\mu=(\Re \zeta_1, \Im \zeta_1, \cdots,\Re \zeta_{h}, \Im \zeta_h),
\qquad 
\nu=(\Re \xi_1, \Im \xi_1, \cdots, \Re \xi_{h}, \Im \xi_h),
\end{equation*}
identity~\eqref{eq:concise_identity} follows from the following identity for $\mu, \nu\in \R^{2h}$
\begin{multline}
\label{eq:concise_identity_real}
|\mu|^p + (p-1)|\nu|^p - p |\nu|^{p-2}  \langle \nu,\mu \rangle_{\R^{2h}}
=p|\mu-\nu|^2 \int_0^1 s\, |s \nu + (1-s) \mu|^{p-2}\, ds\\ 
+ p(p-2) \int_0^1 s\, |s \nu + (1-s) \mu|^{p-4} \big[\langle \mu-\nu,s \nu + (1-s) \mu\rangle_{\R^{2h}}\big]^2 ds,
\end{multline}
where we have used the trivial facts 
\begin{equation*}
	|\mu|_{\R^{2h}}^2
	=|\zeta|_{\C^h}^2,
	\qquad 
	|\nu|_{\R^{2h}}^2
	=|\xi|_{\C^h}^2,
	\qquad
	\langle \nu, \mu \rangle_{\R^{2h}}=\langle \xi, \zeta \rangle_{\C^{h}}. 
\end{equation*}
In order to prove~\eqref{eq:concise_identity_real} we define
\begin{equation*}
z(\mu_1, \dots, \mu_{2h}) := |\mu|^p + (p-1)|\nu|^p - p|\nu|^{p-2}\langle \nu\cdot\mu \rangle_{\R^{2h}}.
\end{equation*}
Applying the multivariate Taylor expansion of $z$ around the point \(\nu\), we have
\begin{equation*}  
z(\mu_1, \dots, \mu_{2h}) = \sum_{|\alpha|=2} R_\alpha(\mu) (\mu - \nu)^\alpha,
\qquad \text{where} \qquad
R_\alpha(\mu) := \frac{2}{\alpha!} \int_0^1 (1 - t) D^\alpha z\,(\nu + t(\mu - \nu)) \, dt.  
\end{equation*}
Here we have used that $\sum_{|\alpha|\leq 1} \frac{D^\alpha z(v)}{\alpha !} (\mu-\nu)^\alpha=0.$
A direct computation shows that $D^\alpha z(\mu),$ for $|\alpha|=2$ is as follows
\begin{equation*}
	D^{\alpha_i} D^{\alpha_j} z(\mu)
	= p(p-2) |\mu|^{p-4} \mu_i \mu_j + p |\mu|^{p-2}\delta_{ij},	
	\qquad i,j=1,2, \dots, 2h,
\end{equation*}
where $\delta_{ij}$ represents the Kronecker symbol.
This gives 
\begin{multline*}
z(\mu_1, \dots, \mu_{2h})
=2 \sum_{\substack{i,j=1\\ i< j}}^{2h}
\int_0^1 (1-t) p(p-2) |\nu-t(\mu-\nu)|^{p-4} (\nu_i-t(\mu_i-\nu_i)) (\mu_i-\nu_i) (\nu_j-t(\mu_j-\nu_j)) (\mu_j-\nu_j)\, dt\\
+ \sum_{i=1}^{2h}
\int_0^1 (1-t) \big[ p(p-2) |\nu-t(\mu-\nu)|^{p-4} (\nu_i-t(\mu_i-\nu_i))^2(\mu_i-\nu_i)^2
+ p |\nu-t(\mu-\nu)|^{p-2} (\mu_i-\nu_i)^2 \big]\, dt,
\end{multline*}
which, rearranging all the terms, become
\begin{multline*}
z(\mu) =
p(p-2) \int_0^1 (1-t) |\nu-t(\mu-\nu)|^{p-4} [(\nu-t(\mu-\nu))\cdot (\mu-\nu)]^2\, dt
+p \int_0^1 (1-t) |\nu-t(\mu-\nu)|^{p-2} |\mu-\nu|^2\, dt.
\end{multline*}
A final change of variable \(s=1-t\) yields the desired identity \eqref{eq:concise_identity_real}.
\end{proof}

\subsection*{A cut-off function}
In order to prove sharpness of the constant in various inequalities of this manuscript we will make use of a suitable cut-off function. Let $\varepsilon>0,$ consider a family of functions $g_\varepsilon \in C^\infty_c(\mathbb{R}_+),$ $0 \leq g_\varepsilon \leq 1,$ such that   
\begin{equation}
\label{eq: definizione g_veps}
g_\varepsilon(r) =
\begin{cases} 
0 & \text{if } 0 \leq r \leq \varepsilon \text{ or } r \geq \frac{1}{\varepsilon}, \\[5pt]
1 & \text{if } 2\varepsilon \leq r \leq \frac{1}{2\varepsilon},
\end{cases}
\qquad 
\text{and}
\qquad 
|g'_\varepsilon(r)| \leq
\begin{cases} 
 \frac{c}{\varepsilon} & \text{for } \varepsilon \leq r \leq 2\varepsilon, \\[5pt]
 c\varepsilon & \text{for } \frac{1}{2\varepsilon} \leq r \leq \frac{1}{\varepsilon},
\end{cases}
\end{equation}
for some constant \(c > 0\).

It is straightforward to show that the following properties holds true.
\begin{proposition}\label{prop:straighforward-g_eps}
Given $\varepsilon>0,$ let $g_\varepsilon \in C^\infty_c(\R_+),$ $0\leq g_\varepsilon \leq 1$ and assume~\eqref{eq: definizione g_veps}. Then one has 
\begin{align}
\label{eq:g_eps1}  
\int_0^{+\infty} r^{-1} |g_\varepsilon(r)|^p \, dr &= -\ln(4\varepsilon^2) + \mathcal{O}(1),\\
\int_0^{+\infty} |g_\varepsilon(r)|^{p-1} |g'_\varepsilon(r)| \, dr &= \mathcal{O}(1),
\label{eq:g_eps2}\\
\int_0^{+\infty} r^{p-1} |g'_\varepsilon(r)|^p \, dr &= \mathcal{O}(1).
\label{eq:g_eps3}
\end{align}
\end{proposition}
We additionally have the following result whose proof also follows easily from the definition of $g_\varepsilon$ given in~\eqref{eq: definizione g_veps}. Nevertheless, for reader's convenience, we provide the complete proof below.  
\begin{proposition}\label{prop: asintotico_g_veps_logaritmico}
	Given $\varepsilon>0,$ define $g_\varepsilon$ as in Proposition~\ref{prop:straighforward-g_eps}. Let $R>0,$ then one has 
\begin{align}\label{eq:g_eps-log1}
\int_0^{R} r^{-1} \frac{1}{\ln\left(\frac{R}{r}\right)} \left|g_\varepsilon\left(\ln\left(\frac{R}{r}\right)\right)\right|^p \, dr &= -\ln(4 \varepsilon^2) + \mathcal{O}(1),\\[7pt]
\int_0^{R} r^{-1} \left|g_\varepsilon\left(\ln\left(\frac{R}{r}\right)\right)\right|^{p-1} \left|g'_\varepsilon\left(\ln\left(\frac{R}{r}\right)\right)\right| \, dr &= \mathcal{O}(1),
\label{eq:g_eps-log2}\\[7pt]
\int_0^{R} r^{-1} \ln\left(\frac{R}{r}\right)^{p-1} \left|g'_\varepsilon\left(\ln\left(\frac{R}{r}\right)\right)\right|^p \, dr &= \mathcal{O}(1),
\label{eq:g_eps-log3}
\end{align}
as $\varepsilon$ goes to $0^+.$
\end{proposition}

\begin{proof}
We start with the proof of~\eqref{eq:g_eps-log1}.
From the definition of \(g_\varepsilon\) given in~\eqref{eq: definizione g_veps}, the integral in the ledt-hand-side of~\eqref{eq:g_eps-log1} can be split as follows
\[
\int_0^{R} r^{-1} \frac{1}{\ln\left(\frac{R}{r}\right)} \left|g_\varepsilon\left(\ln\left(\frac{R}{r}\right)\right)\right|^p \, dr = I_1 + I_2 + I_3,
\]
where
\begin{equation*}
I_1 := \int_{R e^{-\frac{1}{\varepsilon}}}^{R e^{-\frac{1}{2 \varepsilon}}} r^{-1} \frac{1}{\ln\left(\frac{R}{r}\right)} \left|g_\varepsilon\left(\ln\left(\frac{R}{r}\right)\right)\right|^p \, dr,
\qquad  \qquad
I_2 := \int_{R e^{-2\varepsilon}}^{R e^{-\varepsilon}} r^{-1} \frac{1}{\ln\left(\frac{R}{r}\right)} \left|g_\varepsilon\left(\ln\left(\frac{R}{r}\right)\right)\right|^p \, dr, 
\end{equation*}
and 
\[
I_3 := \int_{R e^{-\frac{1}{2\varepsilon}}}^{R e^{-2\varepsilon}} r^{-1} \frac{1}{\ln\left(\frac{R}{r}\right)} \, dr. 
\]
It is easy to show, then, that $I_1=\mathcal{O}(1)$ and that $I_2=\mathcal{O}(1).$

Moreover,
\begin{equation*}
	I_3= \left[-\ln\left(\ln\left(\frac{R}{r}\right)\right)\right]_{R e^{-2\varepsilon}}^{R e^{-\frac{1}{2\varepsilon}}}=-\ln(4 \varepsilon^2).
\end{equation*}
Combining those together gives~\eqref{eq:g_eps-log1}. As for~\eqref{eq:g_eps-log2} one has
\[
\int_0^{R} r^{-1} \left|g_\varepsilon\left(\ln\left(\frac{R}{r}\right)\right)\right|^{p-1} \left|g'_\varepsilon\left(\ln\left(\frac{R}{r}\right)\right)\right| \, dr = J_1 + J_2,
\]
where we have defined
\begin{multline*}
J_1 := \int_{R e^{-\frac{1}{\varepsilon}}}^{R e^{-\frac{1}{2 \varepsilon}}} r^{-1} \left|g_\varepsilon\left(\ln\left(\frac{R}{r}\right)\right)\right|^{p-1} \left|g'_\varepsilon\left(\ln\left(\frac{R}{r}\right)\right)\right| \, dr, \\
J_2 := \int_{R e^{-2 \varepsilon}}^{R e^{-\varepsilon}} r^{-1} \left|g_\varepsilon\left(\ln\left(\frac{R}{r}\right)\right)\right|^{p-1} \left|g'_\varepsilon\left(\ln\left(\frac{R}{r}\right)\right)\right| \, dr. 
\end{multline*}
Now it is easy to show that
\begin{equation*}
	J_1\leq c \varepsilon \int_{R e^{-\frac{1}{\varepsilon}}}^{R e^{-\frac{1}{2 \varepsilon}}} \frac{1}{r} \, dr
= \mathcal{O}(1),
\qquad \text{and} \qquad
J_2\leq \frac{c}{\varepsilon} \ln\left(e^{\varepsilon}\right) = \mathcal{O}(1).
\end{equation*}
For~\eqref{eq:g_eps-log3}, one similarly has
\[
\int_0^{R} r^{-1} \ln\left(\frac{R}{r}\right)^{p-1} \left|g'_\varepsilon\left(\ln\left(\frac{R}{r}\right)\right)\right|^p \, dr = K_1 + K_2,
\]
where we have defined 
\begin{equation*}
K_1 := \int_{R e^{-\frac{1}{\varepsilon}}}^{R e^{-\frac{1}{2 \varepsilon}}} r^{-1} \ln\left(\frac{R}{r}\right)^{p-1} \left|g'_\varepsilon\left(\ln\left(\frac{R}{r}\right)\right)\right|^p \, dr, 
\qquad 
K_2 = \int_{R e^{-2 \varepsilon}}^{R e^{-\varepsilon}} r^{-1} \ln\left(\frac{R}{r}\right)^{p-1} \left|g'_\varepsilon\left(\ln\left(\frac{R}{r}\right)\right)\right|^p \, dr 
\end{equation*}
Since $K_1=\mathcal{O}(1)$ and $K_2=\mathcal{O}(1)$ one gets~\eqref{eq:g_eps-log3}.
\end{proof}

\section*{Acknowledgements}

L.C was supported by the grant Ramón y Cajal RYC2021-032803-I  funded by MCIN/AEI/10.13039/50110
\noindent
0011033 and by the European Union NextGenerationEU/PRTR and by Ikerbasque.

\nocite{*}
\bibliographystyle{abbrv} 
\bibliography{Bibliografia03}

\end{document}